\documentclass{article}
\usepackage[margin=1in]{geometry}

\usepackage{graphicx, dsfont, amssymb, natbib, color} 

\usepackage{amsmath,amsthm, bbm, multicol,array,mathtools,enumerate, hyperref} 
\newtheorem{theorem}{Theorem}
\newtheorem{lemma}{Lemma}
\newtheorem{prop}{Proposition}

\newtheorem{rem}{Remark}

\usepackage[normalem]{ulem}

\newcommand{\tH}{\widetilde{H}}
\newcommand{\tY}{\widetilde{Y}}
\newcommand{\tX}{\widetilde{X}}
\newcommand{\norm}[1]{\lVert#1\rVert}

\newenvironment{keywords}
  {\par\smallskip\noindent\small\textbf{Keywords:}\ }
  {\par}

\newcommand{\bluesout}{%
  \bgroup
  \markoverwith{\textcolor{blue}{\rule[0.5ex]{2pt}{0.5pt}}}%
  \ULon
}

\newcommand{\bluewave}{%
  \bgroup
  \markoverwith{\lower3.5pt\hbox{\color{blue}\sixly\char58}}%
  \ULon
}

\title{Unmatched data augmentation for classification with Linkit regression}
\author{Fadoua Balabdaoui and Jinyu Chen}
\date{\today}

\begin{document}

\maketitle

\begin{abstract}
We consider binary regression based on a small sample of matched covariate-response pairs and much larger, independent samples of covariates and binary responses whose pairing is unavailable. Unlike classical semi-supervised classification, the auxiliary responses are observed but unpaired. We assume that the matched and unmatched populations share a common regression vector, while allowing their covariate distributions and known link functions to differ. We propose a semi-supervised maximum likelihood estimator that combines the conditional likelihood of the matched sample with a marginal Bernoulli likelihood obtained from the unmatched samples. We establish existence, consistency, and asymptotic normality of the estimator. We then derive explicit measures of the asymptotic gain over the maximum likelihood estimator based only on the matched observations, for both parameter estimation and prediction of conditional success probabilities.  We find that the gain depends on the abundance of unmatched data, the information supplied by the unmatched population, and some aspects of the underlying information geometry. We characterize settings in which the gain vanishes and study special cases involving common and differently scaled link functions. Simulations and an application to the Dry Bean data illustrate the resulting improvements in parameter estimation and out-of-sample classification.
\end{abstract}


\section{Introduction}

Binary regression is one of the basic tools for modeling the relationship between a predictor vector and a two-valued outcome. In the classical setting, one observes matched pairs $(X,Y)$, where $X\in\mathbb R^p$ is a vector of covariates and $Y\in\{0,1\}$ is the response. A common formulation assumes that
$\mathbb P(Y=1\mid X=x)=H(\beta_0^\top x)$, 
where $H:\mathbb R\to(0,1)$ is a known increasing link function and $\beta_0\in\mathbb R^p$ is the regression parameter of interest. This formulation includes logistic and probit regression as special cases and more generally belongs to the family of CDF-linked binary regression models. Standard likelihood methods rely crucially on observing $X$ and $Y$ together. When this link is missing, partially or completely, the usual likelihood is no longer available in its ordinary paired form, and even basic questions of identifiability and efficiency become substantially more delicate.

The loss of correspondence between covariates and responses has been studied under several related names, including permuted regression, shuffled regression, unlabeled sensing, unlinked regression, and unmatched regression. In the shuffled or permuted regression formulation, the covariates and responses originate from the same experimental units, but the pairing has been obscured by an unknown permutation. This line of work has revealed sharp differences between ordinary regression and regression under unknown correspondence, both statistically and computationally; see, for instance, \citet{Pananjady:2018hd}, \citet{hsu2017linear}, and \citet{Unnikrishnan:2018gp}. In contrast, unmatched or unlinked regression allows the covariate sample and the response sample to be collected separately, possibly from different sources, so that there may be no latent finite-sample permutation connecting the two observed lists. This perspective was developed for linear regression in \citet{monafadoua2024}, where the equality of paired observations is replaced by an equality in distribution and the estimation problem is studied through a deconvolution viewpoint.

A further distinction concerns the meaning of the word \lq\lq unlabeled\rq\rq. We would like to note that the setting studied here is different from standard semi-supervised classification problems. In much of the semi-supervised learning and semi-supervised inference literature, a labeled observation is a full pair $(X,Y)$, whereas an unlabeled observation consists only of covariates $X$ with the outcome missing. In that framework, the additional unlabeled sample informs the marginal distribution of the covariates and can improve estimation or inference for certain targets; see, for example, \citet{ZhangBrownCai2019}. The usefulness of the unlabeled sample depends on assumptions linking this distribution to the classification rule. Generative mixture models, entropy minimization, and manifold regularization exploit such connections in different ways \cite{CozmanCohenCirelo2003,GrandvaletBengio2005,BelkinNiyogiSindhwani2006}.

Our observation scheme is different. In addition to a small matched sample which carries direct information about the conditional law of $Y$ given $X$, we observe separate samples from the marginal covariate and response distributions without knowing which response corresponds to which covariate. The auxiliary binary responses are therefore not missing, but their links to the auxiliary covariates are unknown. This data structure arises naturally when a limited amount of linked data can be collected or recovered, while larger unlinked samples are available from administrative records, anonymized databases, surveys, or distinct data repositories. In such settings, discarding the unmatched data may waste substantial information, whereas treating the unmatched observations as if they were paired would introduce a fundamental error.


The statistical information in the unmatched response sample enters through the marginal law of $Y$, while the unmatched covariate sample enters through the marginal law of $X$. This provides aggregate information about the class probability and places the problem closer to learning from label proportions than to conventional semi-supervised classification, although our independent marginal samples, small matched component, and inferential objectives require a distinct analysis. The likelihood-based approach studied in this paper presents one way of conducting such an analysis.


For continuous outcomes, \citet{balabdaoui2026SSLEMLE} considered estimation in a semi-supervised unmatched linear regression model in which a small matched sample is combined with larger unmatched samples. Their analysis is based on an empirical likelihood approach and quantifies the extent to which the additional unmatched data can improve statistical precision relative to using only the matched sample. The present binary-response setting is motivated by the same principle but requires a different likelihood structure. The fact that the response is binary means that the unmatched outcome sample contributes through the empirical proportion of successes, while the unmatched covariate sample contributes through the empirical distribution of the covariates inside the link-integrated marginal probability. 
The resulting likelihood therefore combines a matched component, corresponding to ordinary link-based binary regression, with an unmatched component, corresponding to the marginal Bernoulli law induced by the covariate distribution and the link function. Another difference from the work in \citet{balabdaoui2026SSLEMLE} is that we will allow for \textit{different} link functions in the matched and unmatched samples and also \textit{different} distributions of the covariates. One of the main results of this paper is that no matter how different the link functions and the distributions of the covariates are, one can only do at least as well as the MLE based solely on the matched sample as long as the true regression vector is the same.


\paragraph{Our contributions.} The main goal of this paper is to determine whether, and by how much, a large unmatched sample can improve estimation and prediction in CDF-linked binary regression when only a small sample of matched covariate-response pairs is available. To this end, we introduce a semi-supervised maximum likelihood estimator that combines the ordinary binary-regression likelihood from the matched sample with a marginal Bernoulli likelihood constructed from the unmatched covariates and responses. The latter contribution depends on the empirical proportion of successes in the unmatched response sample and on the empirical average of the link-transformed linear predictors in the unmatched covariate sample. Our framework allows the matched and unmatched populations to have different covariate distributions and different known link functions. The two samples are connected only through the common regression vector $\beta_0$. This flexibility is important in applications in which the matched and unmatched observations are collected from different populations, instruments, or data repositories. Despite these differences, we show that incorporating the unmatched sample cannot decrease asymptotic precision relative to the maximum likelihood estimator based only on the matched observations. The paper makes four main contributions. First, we formulate a likelihood-based framework for augmenting matched binary-regression data with independent samples from the marginal covariate and binary-response distributions. This observation scheme is distinct from classical semi-supervised classification, where the auxiliary sample contains covariates with missing labels. In our setting, the auxiliary binary responses are observed, but their correspondence with the auxiliary covariates is unavailable. Second, we establish existence and consistency of the resulting semi-supervised maximum likelihood estimator under general conditions on the link functions and covariate distributions. The existence theory accounts explicitly for the possibility of separation in the matched binary-regression sample and shows that, with probability tending to one, maximization may be restricted to a fixed bounded parameter set. Third, we derive the asymptotic distribution of the proposed estimator and obtain explicit expressions for its asymptotic covariance matrix. We use these expressions to quantify the statistical gain over the matched-sample maximum likelihood estimator, both for estimation of the regression vector and for prediction of conditional success probabilities. The gain is characterized in terms of the relative amount of unmatched data, the signal carried by the unmatched sample, and the alignment between the prediction direction and the direction in which the unmatched data provide information. We further identify settings in which the gain vanishes and settings in which it can be substantial. Fourth, we investigate important special cases, including settings in which the matched and unmatched samples share the same link function and covariate distribution, and settings in which the two links belong to the same family but have different known scales. Simulation studies confirm the theoretical efficiency comparisons. An application to the Dry Bean data illustrates that the proposed estimator can improve both parameter estimation and out-of-sample classification, and can produce a more concentrated estimate of the decision boundary than the matched-sample estimator. 

\paragraph{Organization of the paper.} Section 2 introduces the matched and unmatched binary-regression models, defines the combined empirical log-likelihood, and states the assumptions used throughout the paper. Section 2.1 studies existence of the semi-supervised maximum likelihood estimator, first for a fixed sample and then asymptotically, while Section 2.2 establishes consistency. Section 3 derives the weak limit of the estimator. After obtaining the asymptotic distribution of the score and Hessian, we establish asymptotic normality and derive closed-form measures of the statistical gain over the maximum likelihood estimator based solely on the matched sample. The comparison covers both estimation of $\beta_0$ and prediction of $H(\beta_0^\top x)$ at a new covariate value $x$. Section 4.1 specializes the theory to the case in which the matched and unmatched samples have the same link function and covariate distribution. The section identifies configurations with no asymptotic gain and examines illustrative examples, including logistic regression with a uniform covariate distribution. Section 4.2 considers matched and unmatched samples whose link functions belong to the same family but have different known scale parameters. Section 5 presents an application to the Dry Bean data. The empirical study compares the matched-sample and semi-supervised estimators in terms of mean squared error, test-set misclassification rate, and estimation of the classification boundary. Section 6  concludes the paper and discusses extensions to partially specified and shape-constrained link functions, including links known only up to scale and links satisfying the log-concavity or bi-log-concavity constraint.

\begin{figure}
    \centering
    \includegraphics[width=0.8\linewidth]{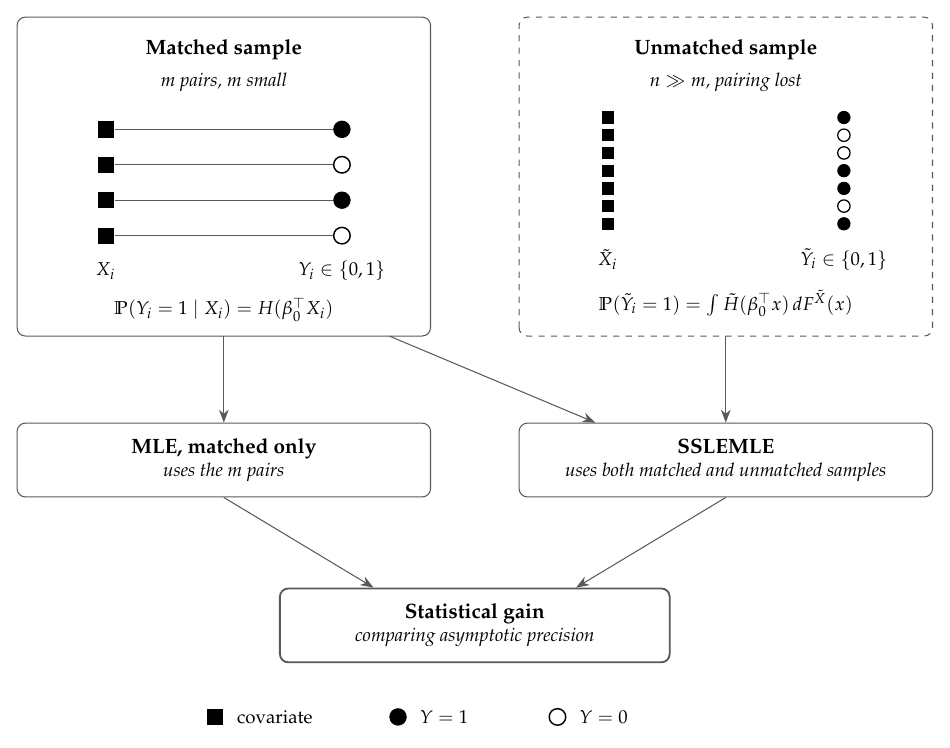}
    \caption{A small matched sample of size $m$ is augmented with a large unmatched sample of size $n\gg m$, in which covariates and binary responses are observed without their pairing. The statistical gain is assessed by comparing the MLE based on the matched sample only with SSLEMLE based on both samples.}
    \label{fig:placeholder}
\end{figure}

\section{Preliminaries and assumptions}

Let $H: \mathbb R \mapsto (0,1) $ and $\tH: \mathbb R \mapsto (0,1) $ be known strictly increasing  and continuous functions. Let $X \in \mathbb R^p$ be a predictor and $Y \in \{0,1\}$ the outcome of interest. In the case where $X$ and $Y$ are linked, that is, when they can be put together as a pair $(X, Y)$ defined on the same probability space, we will assume that the conditional distribution of $Y$ is determined through the relationship 
\begin{align*}
\mathbb P(Y = \delta | X)  =  H(\beta_0^\top X)^\delta (1-H(\beta_0^\top X))^{1-\delta}
\end{align*}
for $\delta \in \{0,1\}$ and some unknown regression vector $\beta_0 \in \mathbb R^p$. In the case where $\tX$ and $\tY$ are unmatched, that is, when we lose the link between $\tX$ and $\tY$. We assume that their marginal distributions are compatible with an underlying regression model: there exists a pair of random variables $(\widetilde{x}_i,\widetilde{y}_i)$, with
\begin{align*}
    \widetilde{x}_i \sim \tX \quad \text{and} \quad \widetilde{y}_i\sim \tY,
\end{align*}
such that
\begin{align*}
     \mathbb P(\widetilde{y}_i = \delta | \widetilde{x}_i)  =  \tH(\beta_0^\top \widetilde{x}_i)^\delta (1-\tH(\beta_0^\top \widetilde{x}_i))^{1-\delta}.
\end{align*}
Consequently, the marginal distribution of \(\widetilde Y\) satisfies: 
\begin{align*}
\mathbb P(\tY = \delta)  =  \left(\int  \tH (\beta_0^\top x) dF^{\widetilde{X}}(x)\right)^\delta  \left(1- \int \tH (\beta_0^\top x) dF^{\widetilde{X}}(x) \right)^{1-\delta},
\end{align*}
where $F^{\tX}$ denotes the distribution function of $\tX$.

Let $\widetilde X_i \in \mathbb{R}^p,\, i=1, \ldots, n$ and  $\widetilde Y_j \in 
\{0,1\},\, j =1, \ldots, n$ be the unmatched predictors and responses, and $(X_k, Y_k) \in \mathbb R^p \times \{0,1\}, \, k=1, \ldots, m$ the pairs in the matched sample.  We denote by $\mathbb F^{\widetilde X}_n$ the empirical distribution of the unmatched covariates $\widetilde X_i, \ i =1, \ldots, n$.  Let $\overline{\widetilde Y}_n$ denote the empirical average of the unmatched responses. 

Then, the empirical log-likelihood based on the matched and unmatched samples is given by   
\begin{eqnarray}\label{loglik}
    \ell_{n,m}  =  \frac{n}{n+m} \ell_n +  \frac{m}{n+m} \ell_m
\end{eqnarray}    
where
\begin{eqnarray}\label{elln}
\ell_n(\beta) =  \overline{\widetilde Y}_n \log\left(  \int \tH(\beta^\top x) d\mathbb F^{\widetilde X}_n(x)\right) + (1-\overline{\widetilde Y}_n) \log\left(1- \int \tH(\beta^\top x) d\mathbb F^{\widetilde X}_n(x)\right) 
\end{eqnarray}
and
\begin{eqnarray}\label{ellm}
\ell_m(\beta) = \frac{1}{m} \sum_{k=1}^m \Big \{ Y_k \log\left( H(\beta^\top X_k)\right)  + (1-Y_k) \log\left(1 - H(\beta^\top X_k)\right) \Big\}.
\end{eqnarray}
In what follows, we specify the assumptions under which we will study the asymptotics of a maximizer of (\ref{loglik}).

\begin{itemize}
    \item (A0):  $\lim_{m, n \to \infty}  m/n = \lambda \in (0,1)$. 
    \item (A1): The functions $H: \mathbb R \to (0,1)$ and $\tH: \mathbb R \mapsto (0,1) $ are continuous and strictly monotone on $\mathbb R$,  and  $H(\mathbb R) = \tH(\mathbb R) = (0,1)$.

    \item (A1'): $H$ is bi-log-concave: i.e., $H$ and $1-H$ are log-concave. 
    
    \item (A2): The distribution of $X$ satisfies the property: For any given $R > 0$ and all $\beta \in \overline{\mathcal B}(0, R)$, there exist $0 < c_R \le C_R < 1$ such that 
     \begin{eqnarray*}
     c_R  \le \int H(\beta^\top x)\, dF^X(x) \le C_R .
     \end{eqnarray*}
     
      \item (A3): For any given $R > 0$, the distribution of $X$ satisfies 
      $$ \int \left(\log{H(-R\lVert x\rVert)}\right)^2 d F^X(x) < \infty$$ and 
      $$\int \left(\log{\left(1-H(R\lVert x\rVert)\right)}\right)^2 d F^X(x) < \infty.$$
      \item (A4): $\operatorname{supp}(X) $ has a non-empty interior and $X$ admits a density with respect to the Lebesgue measure on $\mathbb R^p$.
      \item (A5): Assumptions (A1') (A2), (A3), and (A4) should also hold with $X$ and $H$ replaced by $\tX$ and $\tH$, respectively. We denote the corresponding assumptions by (A5-1'),  (A5-2), (A5-3), and (A5-4), respectively.
      
\end{itemize}
Note that Assumption (A1') holds if $H$ admits a log-concave density with respect to the Lebesgue measure. Both Assumptions (A2) and (A3) are satisfied for a compactly supported predictor $X$. In fact, if $L > 0$ is such that $\mathbb P(\Vert X \Vert \le L) =1$, then we can take
\begin{eqnarray*}
c_R =  H(-R L), \ \text{and} \  C_R =  H(R L).
\end{eqnarray*}
In addition, $H(- R \Vert x \Vert ) \ge H(- R L)$ with probability 1, and hence 
$$\int (\log H(- R \Vert x \Vert ))^2 dF^X(x) \le  (\log H(- RL))^2 < \infty.$$
Similarly, $\int \log(1 - H(R \Vert x \Vert))^2 dF^X(x)   \le  \log{(1- H( RL))}^2 < \infty$.

Assumption (A2) is also satisfied for $X \sim \mathcal{N}(\mu, \Sigma)$ for some $\mu \in \mathbb R^p$ and $\Sigma \in \mathbb R^{p \times p}$ such that $\mathbb P(X = \mu) = 0$; i.e., $\Sigma \ne 0$. In fact, if $\varphi$ denotes the probability density function of $\mathcal N(0,1)$, we can write
\begin{eqnarray*}
\int H(\beta^\top x)  dF^X(x)  =  \int H( m_\beta + \sigma_\beta z) \varphi(z) dz  
\end{eqnarray*}
with $m_\beta =  \beta^\top \mu$ and $\sigma_\beta =  \sqrt{\beta^\top \Sigma \beta}$. Then,
\begin{eqnarray*}
\int H( m_\beta + \sigma_\beta z) \varphi(z) dz  \ge   \int_0^\infty H( m_\beta + \sigma_\beta z) \varphi(z) dz    \ge    \frac{H(m_\beta)}{2} \ge \frac{H(- \Vert \mu \Vert R)}{2} > 0
\end{eqnarray*}
using the Cauchy-Schwarz inequality and Assumption (A1). On the other hand,  if $\lambda_{\max} > 0$ is the largest eigenvalue of $\Sigma $, we have $\beta^\top \Sigma \beta \le  \lambda_{\max} R^2$.   Hence, 
\begin{eqnarray*}
\int H( m_\beta + \sigma_\beta z) \varphi(z) dz  \le   \int H( R \Vert \mu \Vert+ \sqrt{\lambda_{\max}} R \ \vert z \vert) \varphi(z) dz < 1
\end{eqnarray*}
since otherwise $\int (1- H(  R \Vert \mu \Vert+ \sqrt{\lambda_{\max}} R \ \vert z \vert)) \varphi(z) dz = 0$, and this can only occur if $H( R \Vert \mu \Vert+ \sqrt{\lambda_{\max}} R \ \vert z \vert) = 1$ for a.e. $z$, which is impossible by Assumption (A1). Therefore, for $X \sim \mathcal N(\mu, \Sigma)$ 
\begin{eqnarray*}
c_R  =  \frac{H(- R \Vert \mu \Vert)}{2}, \ \textrm{and} \ \ C_R = \int H( R \Vert \mu \Vert+ \sqrt{\lambda_{\max}} R \ \vert z \vert) \varphi(z) dz. 
\end{eqnarray*}
Assumption (A3) is satisfied for the logit and probit links and $X \sim \mathcal N(\mu, \Sigma)$. We start with logit link, i.e., $H(x) = \exp{(x)}/(1+\exp{(x)})$. We have that
\begin{align*}
    \log{H(x)} &= x - \log{(1+e^x)} \\
    &\geq x - \left(\log{2} \cdot \mathbbm{1}_{\{x \leq 0\}} + (x + \log{2}) \cdot \mathbbm{1}_{\{x > 0\}} \right)\\
    &= (x-\log{2}) \cdot \mathbbm{1}_{\{x \leq 0\}} + (-\log{2}) \cdot \mathbbm{1}_{\{x > 0\}} 
\end{align*}
for all $x \in \mathbb{R}$. Since $\log{H(x)} < 0$ for all $x\in \mathbb{R}$, then is follows that 
\begin{align*}
     \int \left(\log{H(-R\lVert x\rVert)}\right)^2 d\mathbb{P}^X(x) 
     &\leq \int \left(-R\lVert x\rVert- \log{2}\right)^2 d\mathbb{P}^X(x)\\
     &= R^2 \ \mathbb{E}_X[\lVert X\rVert^2] + R\log{(4)}\ \mathbb{E}_X[\lVert X\rVert] +  (\log{2})^2 \\
     &< \infty.
\end{align*}
For the probit link, we have $H(x) = \Phi(x)$. By the Mill's ratio we have for all $t <0$
\begin{align*}
    \Phi(t) > \frac{-t}{t^2+1} \varphi(t)
\end{align*}
which implies that for all $t \in \mathbb{R}$ we have
\begin{align*}
   0 &> \log{H(t)}\\
   &\geq (\log{(-t)} - \log{(t^2+1)} - \frac{1}{2}\log{(2\pi)} - \frac{t^2}{2}) \cdot \mathbbm{1}_{\{t < -1\}} + (\log \Phi(-1))  \cdot \mathbbm{1}_{\{t \geq -1 \}}  \\
    &\geq (\log{(-t)} -\frac{1}{2}\log{(2\pi)} - \frac{3}{2} t^2) \cdot \mathbbm{1}_{\{t < -1\}}+ (\log \Phi(-1))  \cdot \mathbbm{1}_{\{t \geq -1 \}}  \\
    & \geq (-\frac{1}{2}\log{(2\pi)} - \frac{3}{2} t^2) \cdot \mathbbm{1}_{\{t < -1\}} + (\log \Phi(-1))  \cdot \mathbbm{1}_{\{t \geq -1 \}}.
\end{align*}
Then it follows that 
\begin{align*}
     &\int \left(\log{H(-R\lVert x\rVert)}\right)^2 d\mathbb{P}^X(x) \\
     &\leq \int \left(\frac{3}{2}R^2\lVert x\rVert^2 + \frac{1}{2}\log{(2\pi)}\right)^2 d\mathbb{P}^X(x) + (\log \Phi(-1))^2\\
     &\leq \frac{9}{4}R^4\ \mathbb{E}_X[\lVert X\rVert^4] + \frac{3}{2}\log{(2\pi)} R^2\ \mathbb{E}_X[\lVert X\rVert^2] + \frac{1}{4} (\log{(2\pi)})^2 + (\log \Phi(-1))^2 \\
     &< \infty.
\end{align*}

\subsection{Existence of the SSLEMLE}

In the following, we show that maximizing the log-likelihood in (\ref{loglik}) can be restricted on a bounded Euclidean ball. In the sequel, $\operatorname{conv}(z_1, \ldots, z_d)$, for some integer $d \ge 1$, denotes the convex hull of $z_1, \ldots, z_d \in \mathbb R^p$ and $\operatorname{int}(S)$ the interior of some set $S \subseteq \mathbb R^p$.  The following lemma and its proof are  inspired by the results of \cite{silvapulle1981}.
\subsubsection{The finite sample case}
\begin{lemma}\label{lemma: existence-auxiliary}
    For $z_1, \cdots, z_d \in \mathbb{R}^p$, if $0\in\operatorname{int}(\operatorname{conv}(z_1,\cdots,z_d))$, then for any $v \in \mathcal{S}^{p-1}$ there exists at least one $k \in \{1,\cdots,d\}$ with $v^\top z_k <0$ and hence
    $$ \inf_{v \in \mathcal{S}^{p-1}}\frac{1}{d}\sum_{i=1}^d \mathbbm{1}_{\{v^\top z_i < 0\}}  >0.$$
\end{lemma}
\begin{proof}
    If there exists $v' \in \mathcal{S}^{p-1}$ such that $v'^\top z_i \geq 0$ for all $i =1,\cdots,d$, then for any $y = \sum_{i=1}^d \alpha_i z_i \in \operatorname{conv}(z_1,\cdots,z_d)$ with $\alpha_i \geq 0$ for all $i=1,\cdots,d$ and $\sum_{i=1}^d \alpha_i=1$ we have 
    \begin{align}\label{nonneg}
        v'^\top y = v'^\top \sum_{i=1}^d \alpha_i z_i = \sum_{i=1}^d \alpha_i \left(v'^\top z_i\right) \geq 0.
    \end{align}
    Since $0\in\operatorname{int}(\operatorname{conv}(z_1,\cdots,z_d))$, there exists $\epsilon > 0$, such that $\mathcal{B}(0, \epsilon) \subseteq \operatorname{conv}(z_1,\cdots,z_d)$. Therefore, $-\frac{\epsilon}{2}v' \in \mathcal{B}(0, \epsilon) \subseteq \operatorname{conv}(z_1,\cdots,z_d)$. But we have 
    \begin{align*}
         v'^\top \left(-\frac{\epsilon}{2}v'\right) = -\frac{\epsilon}{2} < 0
    \end{align*}
    which is in contradiction with (\ref{nonneg}).  We conclude that for any $v \in \mathcal{S}^{p-1}$ there exists at least one $k \in \{1,\cdots,d\}$ with $v^\top z_k <0$. This means that for all $v \in \mathcal{S}^{p-1}$
    \begin{align*}
        \frac{1}{d}\sum_{i=1}^d \mathbbm{1}_{\{v^\top z_i < 0\}} > 0.
    \end{align*}
    Since $\mathcal{S}^{p-1}$ is a compact set and $v \mapsto \frac{1}{d}\sum_{i=1}^d \mathbbm{1}_{\{v^\top z_i < 0\}}$ is lower semi-continuous, there exists $v^\ast  \in \mathcal{S}^{p-1}$ with
    \begin{align*}
        \inf_{v \in \mathcal{S}^{p-1}}\frac{1}{d}\sum_{i=1}^d \mathbbm{1}_{\{v^\top z_i < 0\}} = \frac{1}{d}\sum_{i=1}^d \mathbbm{1}_{\{{v^\ast}^\top z_i < 0\}} >0.
    \end{align*}
\end{proof}

\begin{lemma}\label{lemma: bi log concavity left derivative}
Assumption (A1') implies that for all $t \in \mathbb{R}$ we have
\begin{align*}
    \log{H}(t) \leq \log{H}(0) + \lim_{s\to 0^-}\frac{\log{H}(s)-\log{H}(0)}{s} \cdot t \leq \frac{H_-'(0)}{H(0)} t
\end{align*}
and 
\begin{align*}
    \log{(1-H(t))} &\leq \log{(1-H(0))}  + \lim_{s\to 0^-}
    \frac{\log{(1-H(s))} - \log{(1-H(0))} }{s}  \cdot t \leq -\frac{H_-'(0)}{1-H(0)} t
\end{align*}
where $H_-'(0)$ denotes the derivative of $H$ to the left at 0.
\end{lemma}

\bigskip

\begin{proof}
 Put $g = \log H$. Then, $g$ is concave on $\mathbb R$, and hence 
$$
g(t)  \le g(t_0)  +  (t-t_0)  g'_{-}(t_0)
$$
where $g'_{-}(t_0)$ is the derivative of $g$ to the left of $t_0$. Applying this inequality to $t_0 =0$ and using the fact that 
$$
g'_{-}(0)  = \frac{H'_{-}(0)}{H(0)}
$$
the first inequality follows.  The second one follows from taking $g = \log(1-H)$. 
\end{proof}
\medskip

Note that since $\log H$ and $\log(1-H)$ are negative on $\mathbb R$, the first and second inequalities of Lemma \ref{lemma: bi log concavity left derivative} are trivial for $t \ge 0 $ and $t \le 0$ respectively. However, they will be used in the proof of Lemma \ref{lemma: Existence finite sample case} below on the \lq\lq non-trivial\rq\rq \  domains.    \\

Let $r_m$ denote the number of matched observations $(X_i, Y_i), i \in \{1, \ldots, m \}$ such that 
$Y_i = 1$. We re-order the matched data so that the first $r_m$ matched observations correspond to the successful responses. In the next lemma we establish the the existence of a solution to the likelihood maximization if certain conditions hold. 

\medskip

\begin{lemma}\label{lemma: Existence finite sample case}
    Let $H$ be a function that satisfies Assumptions (A1) and (A1').
    If any of the following two conditions holds, then the empirical log-likelihood function $\ell_{n,m}(\beta)$ admits at least one maximizer. \\
    
    \par \noindent $\bullet$ \textbf{Condition 1:} $0 < r_m < m$ and \ 
         $0\in\operatorname{int}(\operatorname{conv}(X_1,\cdots,X_{r_m}, -X_{r_m+1}, -X_{r_m+2},\cdots, -X_m))$. \\
    \par \noindent $\bullet$ \textbf{Condition 2:} $r_m \in\{0,m\}$ and \ $0\in\operatorname{int}(\operatorname{conv}(X_1,\cdots,X_m))$.
\end{lemma}
\begin{proof}
    The result follows from Lemmas \ref{lemma: existence-auxiliary} and \ref{lemma: bi log concavity left derivative}. The details are provided in the Appendix \ref{Appendix: Existence}. The idea is to show that $\beta$ cannot maximize $\ell_{n,m}$ if $\norm{\beta}$ is sufficiently large.
\end{proof}

\subsubsection{The asymptotic case}

We next show that, for sufficiently large $n$ and $m$, there exists $R^\ast>0$ such that the maximizer of $\ell_{n,m}$ lies in $\overline{\mathcal{B}}(0,R^\ast)$ almost surely. We need four lemmas to establish the proof.

\begin{lemma} \label{Lemma: Case I holds a.s.}
Suppose that (A2) holds. Then, there exists almost surely an integer $\tilde{m} \ge 1$ such that for all $m\geq \tilde{m}$ there is at least one $i \leq m$ with $Y_i = 1$ and at least one $j \leq m$ with $Y_j = 0$.
\end{lemma}

\begin{lemma} \label{lemma: Case I holds a.s. 2}
    Suppose that Assumptions (A1) and (A4) hold. Then,  there exists  almost surely an integer $\tilde{m} \ge 1$ such that for all $m \geq \tilde{m}$ we have $$0\in\operatorname{int}(\operatorname{conv}(X_1,\cdots,X_{r_m}, -X_{r_m+1},\cdots, -X_m)).$$
\end{lemma}

\begin{lemma} \label{lemma: ULLN existence asymptotic}
    Suppose that Assumption (A3) holds. Then, as $ m \to  \infty$
    \begin{align*}
    &\sup_{v \in \mathcal{S}^{p-1}} \frac{1}{m} \left( \sum_{k=1}^{m} v^\top X_k  \mathbbm{1}_{\{v^\top X_k<0\}} \mathbbm{1}_{\{Y_k=1\}}   +\sum_{s=1}^m (-v^\top X_s ) \mathbbm{1}_{\{v^\top X_s > 0\}} \mathbbm{1}_{\{Y_s=0\}}\right) \\
    & \overset{\text{a.s.}}{\longrightarrow} \sup_{v\in \mathcal{S}^{p-1}} \mathbb{E}\left[v^\top X  \mathbbm{1}_{\{v^\top X<0\}}  \mathbbm{1}_{\{Y=1\}} -v^\top X  \mathbbm{1}_{\{v^\top X>0\}}  \mathbbm{1}_{\{Y=0\}} \right]. 
\end{align*}
\end{lemma}

\begin{lemma} \label{lemma: asymptotic non-zero denominator}
    Suppose that Assumptions (A1), (A3) and (A4) hold. Then,
    \begin{align*}
        \sup_{v\in \mathcal{S}^{p-1}} \left( \mathbb{E}\left[v^\top X  \mathbbm{1}_{\{v^\top X<0\}}  \mathbbm{1}_{\{Y=1\}} \right] + \mathbb{E}\left[-v^\top X  \mathbbm{1}_{\{v^\top X>0\}}  \mathbbm{1}_{\{Y=0\}} \right] \right) < 0.
    \end{align*}
\end{lemma}

The proofs of Lemmas $\ref{Lemma: Case I holds a.s.}$, $\ref{lemma: Case I holds a.s. 2}$, $\ref{lemma: ULLN existence asymptotic}$ and $\ref{lemma: asymptotic non-zero denominator}$ are provided in Appendix \ref{Appendix: Existence auxiliary lemma}.
\medskip
\begin{theorem}
    Suppose that (A0)-(A4) hold. Then,  with probability 1, there exists an integer $m^\ast \ge 1$ such that for all $n,m \ge m^\ast$ a maximizer of $\beta \mapsto \ell_{n,m}(\beta)$ belongs to $\overline {\mathcal{B}} (0, R^\ast)$. 
\end{theorem}

\medskip

\begin{proof}
By Lemma \ref{Lemma: Case I holds a.s.} and \ref{lemma: Case I holds a.s. 2},  there exists with probability 1 an integer $\tilde{m} \ge 1$ such that \textbf{Condition 1} holds for all $m \geq \tilde{m}$. Thus, $\beta \mapsto \ell_{n,m}(\beta)$ admits at least one maximizer. In what follows, we present an explicit closed-form expression for the ball containing the maximizer under Assumptions (A1') as $n,m\to \infty$. \\

\par \noindent By Assumption (A0) we have $\lim_{m, n \to \infty} \frac{m}{n} = \lambda \in (0,1)$. Then, there exists an integer $m_0 \ge 1$, such that for all $n,m \geq m_0$:
    \begin{align*}
        \frac{m+n}{m} \leq 1+\frac{2}{\lambda}.
    \end{align*}
By Lemma \ref{lemma: ULLN existence asymptotic} and \ref{lemma: asymptotic non-zero denominator}, with probability 1, there exists an integer $m_1\geq 1$ such that for all $m \geq m_1$:
\begin{align*}
     &\sup_{v \in \mathcal{S}^{p-1}} \frac{1}{m} \left( \sum_{k=1}^{m} v^\top X_k  \mathbbm{1}_{\{v^\top X_k<0\}} \mathbbm{1}_{\{Y_k=1\}}   +\sum_{s=1}^m (-v^\top X_s ) \mathbbm{1}_{\{v^\top X_s > 0\}} \mathbbm{1}_{\{Y_s=0\}}\right) \\
     &\leq \frac{1}{2} \sup_{v\in \mathcal{S}^{p-1}} \left( \mathbb{E}\left[v^\top X  \mathbbm{1}_{\{v^\top X<0\}}  \mathbbm{1}_{\{Y=1\}} \right] + \mathbb{E}\left[-v^\top X  \mathbbm{1}_{\{v^\top X>0\}}  \mathbbm{1}_{\{Y=0\}} \right] \right)
     <0.
\end{align*}
Then for all $n,m \geq \max{(\tilde{m},m_0,m_1)}$, we have with probability 1
\begin{align*}
    R_{n,m} & \leq  \frac{\frac{ML}{H_-'(0)}\cdot(2+\frac{4}{\lambda})}{\sup_{v\in \mathcal{S}^{p-1}} \left( \mathbb{E}\left[v^\top X  \mathbbm{1}_{\{v^\top X<0\}}  \mathbbm{1}_{\{Y=1\}} \right] + \mathbb{E}\left[-v^\top X  \mathbbm{1}_{\{v^\top X>0\}}  \mathbbm{1}_{\{Y=0\}} \right] \right)} := R^\ast.
\end{align*}
This means that with probability 1 there exists an integer $m^{\ast} \ge 1$ such that for all $n,m \geq m^\ast$ we have that 
$$\{\beta \in\mathbb{R}^p:\left\lVert \beta\right\rVert > R^\ast\} \subseteq  \{\beta\in\mathbb{R}^p:\left\lVert\beta\right\rVert > R_{n,m}\}. $$
Now, from Lemma \ref{lemma: Existence finite sample case} it follows that for all $n, m \geq m^\ast$ $\ell_{n,m}(\beta) <\ell_{n,m}(0)$ for all $\beta$ such that $\Vert \beta\Vert > R^\ast$. Using again continuity of the function $\beta \mapsto \ell_{n,m}(\beta)$ and compactness of closed Euclidean balls, we conclude that with probability 1, there exists an integer $m^\ast$ such that for all $n,m \ge m^\ast$ a maximizer of $\beta \mapsto \ell_{n,m}(\beta)$ belongs to the closed ball $\overline{\mathcal{B}}(0, R^\ast)$.

\end{proof}

\medskip

\subsection{Consistency of the SSLEMLE}
For $\beta \in \mathbb R^p$.  Define 
\begin{eqnarray*}
 \ell_{U}(\beta) &= & \left( \int \tH(\beta_0^\top x) dF^{\tX}(x) \right) \log\left( \int \tH(\beta^\top x) dF^{\tX}(x)  \right) \\
 && \ + \left( 1- \int \tH(\beta_0^\top x) dF^{\tX}(x) \right) \log\left( 1- \int \tH(\beta^\top x) dF^{\tX}(x)  \right)
\end{eqnarray*}
and
\begin{eqnarray*}
\ell_{M}(\beta) = \int \Big (H( \beta_0^\top x) \log H(\beta^\top x) + (1- H(\beta_0^\top x)) \log(1 -H(\beta^\top x)) \Big) dF^X(x)
\end{eqnarray*}
the population entropy for the unmatched and matched parts, respectively. Then, the combined population entropy is given by
\begin{eqnarray*}
\ell =  \frac{1}{\lambda + 1}   \ell_U  +  \frac{\lambda}{\lambda +1} \ell_M.
\end{eqnarray*}
In the sequel, we shall use the following notation.

\begin{itemize}

\item For any measure $Q$ and function $f$ which is integrable with respect to $Q$, $Qf:= \int f dQ$,

\item $\mathbb P_m$ the empirical probability measure based on the matched sample $(X_k, Y_k), k=1, \ldots, m$,

\item $\mathbb P$ the probability measure of a matched pair $(X, Y)$ such that $Y | X \sim \text{Bernoulli}(H(\beta^\top_0 X))$, 

\item $\mathbb P^{\tX}$ the probability measure of $\tX$, which is not necessarily equal to $\mathbb P^{X}$.

\item $\mathbb P^{\tY}$ the (marginal) probability measure of the binary response $\tY$. Note that  $\mathbb P^{\tY}$ is the probability measure of a Bernoulli($p_0$) where  $p_0=  \int  \tH(\beta_0^\top x) dF^{\tX}(x)$,

\item $\mathbb G_m:=  \sqrt m (\mathbb P_m - \mathbb P)$ and $\mathbb G_n = \sqrt n (\mathbb P_n - \mathbb P)$.
\end{itemize}

\medskip

\begin{prop}\label{unifcons}
Under Assumptions (A0) -- (A5), for any given $R > 0$, 
\begin{eqnarray}
\sup_{\beta \in \overline{\mathcal{B}}(0, R)}  \vert \ell_{n, m}(\beta) - \ell(\beta) \vert   = o_{\mathbb P \otimes \mathbb P^{\tX} \otimes \mathbb P^{\tY}}(1).
\end{eqnarray}
And for any $r>0$ it holds that
$$
\sup_{\beta \in \mathcal O_r} \ell(\beta) < \ell(\beta_0),
$$
where $ \mathcal O_r = \mathcal{B}(0,R^\ast) \cap \{\beta: \norm{\beta-\beta_0}> r\}$. This implies that $\widehat \beta_{n,m}  \to_{\mathbb P \otimes \mathbb P^{\tX} \otimes \mathbb P^{\tY}} \beta_0$ as $n,m\to\infty$, where $\widehat \beta_{n,m}  \in \arg\max \ell_{n,m}(\beta)$.
\end{prop}

\begin{proof}
    We show separately that the unmatched and matched parts converge uniformly to their respective population criterion, with errors of the order $ O_{\mathbb P^{\tX} \otimes \mathbb P^{\tY}}(n^{-1/2})$ and $ O_{\mathbb P}(m^{-1/2})$, respectively. Identifiability is ensured by the matched data. The proof is provided in the Appendix \ref{Appendix: Consistency}.
\end{proof}

\medskip

\section{Weak convergence of the SSLEMLE and the statistical gain}
By \cite{DUMBGEN20171}[Theorem 1 (iv)], Assumption (A1), (A1') and (A5-1') imply that $H$ and $\tH$ are continuously
differentiable on $\mathbb{R}$, and $H'(\beta_0^\top X) > 0$ almost surely. Define
\begin{align} \label{eq: Gamma1}
    \Gamma_1 &:= \frac{\left(\int \tH'(\beta_0^\top x) x \, d\mathbb{P}^{\tX}(x) \right)\left(\int \tH'(\beta_0^\top x) x^\top \, d\mathbb{P}^{\tX}(x) \right)}{\int \tH(\beta_0^\top x) d\mathbb{P}^{\tX}(x) (1-\int \tH(\beta_0^\top x) d\mathbb{P}^{\tX}(x))},
\end{align}
\begin{align} \label{eq: Gamma2}
    \Gamma_2 &:= \left(\int \tH'(\beta_0^\top x) x \, d\mathbb{P}^{\tX}(x) \right)\left(\int \tH'(\beta_0^\top x) x^\top \, d\mathbb{P}^{\tX}(x) \right)   \nonumber\\
    &\phantom{=} \times\frac{ \int \tH(\beta_0^\top x)^2\, d\mathbb{P}^{\tX}(x)  - \left(\int \tH(\beta_0^\top x)\, d\mathbb{P}^{\tX}(x)\right)^2 }{\left(\int \tH(\beta_0^\top x) d\mathbb{P}^{\tX}(x)\right)^2 \left( 1- \int \tH(\beta_0^\top x) d\mathbb{P}^{\tX}(x)\right)^2},
\end{align}
and
\begin{align} \label{eq: Sigma2}
    \Sigma_2 := \int xx^\top \frac{\left(H'(\beta_0^\top x)\right)^2}{H(\beta_0^\top x)\left(1-H(\beta_0^\top x)\right)} \, d\mathbb{P}^X(x).
\end{align}
These three matrices will be used in establishing the asymptotic normality result.
\begin{lemma} \label{lemma: Sigma2 finite}
    Under Assumptions (A1), (A1'),  and (A3) we have
    $\Sigma_2 < \infty$.
\end{lemma}

\begin{lemma} \label{lemma: Gamma1 and Gamma2 finite}
    Under Assumptions (A1), (A5-1'), (A5-3) we have $\Gamma_1 < \infty$ and $\Gamma_2 < \infty$.
\end{lemma}

\begin{lemma} \label{lemma: Sigma2 positive definite}
    Under Assumptions (A1),(A1') and (A4), $\Sigma_2$ is positive definite.
\end{lemma}
 The proofs of Lemmas \ref{lemma: Sigma2 finite}, \ref{lemma: Gamma1 and Gamma2 finite} and \ref{lemma: Sigma2 positive definite} are provided in the Appendix \ref{Appendix: convergence auxialary lemma}.

\subsection{Asymptotic normality}
Using Lemmas \ref{lemma: Sigma2 finite}, \ref{lemma: Gamma1 and Gamma2 finite}, and \ref{lemma: Sigma2 positive definite}, we can show the following two lemmas, which together imply asymptotic normality. Their proofs are provided in the Appendix \ref{Appendix: Convergence}.
\begin{lemma} \label{lemma: dot ell convergence in distribution}
Under Assumptions (A0) -- (A5) it holds that
    \begin{align*}
    \sqrt{m+n}\ \dot{\ell}_{n,m}(\beta_0) \overset{d}{\to} \mathcal{N}(0,\,\frac{1}{1+\lambda} (\Gamma_1+\Gamma_2) +\frac{\lambda}{1+\lambda}\Sigma_2),
\end{align*}
where $\Gamma_1$, $\Gamma_2$, and $\Sigma_2$ are given in (\ref{eq: Gamma1}), (\ref{eq: Gamma2}), and \ref{eq: Sigma2}, respectively.
\end{lemma}

\begin{lemma}\label{lemma: ddot ell convergence in probability}
    Under Assumptions (A0) -- (A5) it holds that
    \begin{align*}
      \ddot{\ell}_{n,m}(\beta_0) \to_{\mathbb{P}^{\tX} \otimes \mathbb{P}^{\tY}  \otimes \mathbb{P} } -\frac{1}{1+\lambda} \Gamma_1 - \frac{\lambda}{1+\lambda} \Sigma_2,
      \end{align*}
where $\Gamma_1$ and $\Sigma_2$ are given in (\ref{eq: Gamma1}) and \ref{eq: Sigma2}, respectively.
\end{lemma}

\begin{theorem}
Under Assumptions (A0) -- (A5) it holds that
\begin{align*}
    \sqrt{m+n}\, (\widehat{\beta}_{n,m} - \beta_0)  \to_d \mathcal{N}(0,\, \Sigma_{SSL})
\end{align*}
with 
\begin{align*}
     \Sigma_{SSL} = \left(\frac{1}{1+\lambda} \Gamma_1 + \frac{\lambda}{1+\lambda} \Sigma_2\right)^{-1} \left(\frac{1}{1+\lambda} (\Gamma_1+\Gamma_2) +\frac{\lambda}{1+\lambda}\Sigma_2\right) \left(\frac{1}{1+\lambda} \Gamma_1 + \frac{\lambda}{1+\lambda} \Sigma_2\right)^{-1}.
\end{align*}
$\Gamma_1$, $\Gamma_2$, and $\Sigma_2$ are given in (\ref{eq: Gamma1}), (\ref{eq: Gamma2}), and \ref{eq: Sigma2}, respectively.
\end{theorem}
\begin{proof}
    The result follows by applying Slutsky’s theorem to Lemmas \ref{lemma: dot ell convergence in distribution} and \ref{lemma: ddot ell convergence in probability}. The claim follows from Lemmas \ref{lemma: Sigma2 finite}, \ref{lemma: Gamma1 and Gamma2 finite} and \ref{lemma: Sigma2 positive definite}.
\end{proof}

\subsection{The statistical gain}
In this section, we will use $\operatorname{AVar}$ to denote the asymptotic variance of some asymptotically normal estimator. In the following result, we will derive the expressions of the statistical gain of using the SLLEMLE $\widehat \beta_{n, m}$ over the matched MLE $\widehat  \beta_m$ to estimate $\beta_0$. In addition, we will derive the gain for estimating the predictive probability $H(\beta_0^\top x)$, for a new cavariate $x \in \mathbb{R}^p$, using $\widehat \beta_{n, m}$ instead of $\widehat \beta_m$. We will denote those gains by $G^{\text{par}}$ and $G^{\text{pred}}(x)$ respectively. 

\medskip

\begin{theorem}\label{theorem: Gain}
     Define the quantities
     \begin{align*}
        & a : = \mathbb{E}\left[\tH'(\beta_0^\top \tX)\tX\right],\\
        & b : =  \mathbb{E}\left[\tH(\beta_0^\top \tX)\right] \left(1-  \mathbb{E}\left[\tH(\beta_0^\top \tX)\right]\right), \\
        & c:= \mathbb{E}\left[\tH(\beta_0^\top \tX) \left(1- \tH(\beta_0^\top \tX)\right)\right], \\
        & \Sigma_2 := \mathbb{E} \left[XX^\top \frac{\left(H'(\beta_0^\top X)\right)^2}{H(\beta_0^\top X)\left(1-H(\beta_0^\top X)\right)} \right],\\
        & d := a^\top \Sigma_2^{-1}a, \\
        & \kappa := \frac{c}{b},\\
        & \tau := \frac{d}{b},\\
        & \cos^2{\theta_x} := \frac{\langle x,a\rangle^2_{\Sigma_2^{-1}}}{\lVert x\rVert^2_{\Sigma_2^{-1}}\lVert a\rVert^2_{\Sigma_2^{-1}}} = \frac{\left(x^\top \Sigma_2^{-1} a\right)^2}{\left(x^\top \Sigma_2^{-1} x\right)\left(a^\top \Sigma_2^{-1} a\right)}.
     \end{align*}
   Then,
    \begin{align*}
     &G^{\text{pred}}(x) :=\sqrt{\frac{\operatorname{AVar}\left(H(\widehat{\beta}_{m}^\top x)\right)}{\operatorname{AVar}\left(H(\widehat{\beta}_{n,m}^\top x)\right)}} =  \frac{1}{\sqrt{1 - \frac{(\tau+ \kappa \lambda)\tau}{(\tau+\lambda)^2}\cos^2{\theta_x}}}
     \end{align*}
      for any realization $x \in \mathbb{R}^p$, and 
     \begin{align*}
     &G^{\text{par}} := \sqrt{\frac{\det{\Sigma_{\text{mMLE}}}}{\det{\tilde{\Sigma}_{\text{SSL}}}}}= \frac{1+\frac{1}{\lambda} \tau}{\sqrt{1+\frac{1}{\lambda}\tau (2-\kappa)}}= \max_x G^{\text{pred}}(x).
    \end{align*}
\end{theorem}

\begin{proof}
    We assume that the new realization $x \in \mathbb{R}^p$  is fixed.  By the delta-method we have
    \begin{align*}
        \sqrt{m}\left[H(\widehat{\beta}_{n,m}^\top x) - H({\beta}_0^\top x)\right] \to_d \mathcal{N}\left(0, \left(H'(\beta_0^\top x)\right)^2x^\top\tilde{\Sigma}_{\text{SSL}}\,x\right)
    \end{align*}
    with
    \begin{align*}
        \tilde{\Sigma}_{\text{SSL}} = \left(\frac{1}{\lambda} \Gamma_1 + \Sigma_2\right)^{-1} \left(\frac{1}{\lambda} (\Gamma_1+\Gamma_2) +\Sigma_2\right) \left(\frac{1}{\lambda} \Gamma_1 +  \Sigma_2\right)^{-1}.
    \end{align*}
    We also have
        \begin{align*}
        \sqrt{m}\left[H(\widehat{\beta}_{m}^\top x) - H\left({\beta}_0^\top x\right)\right] \to_d \mathcal{N}\left(0, \left(H'\left(\beta_0^\top x\right)\right)^2x^\top\Sigma_{\text{mMLE}} \,x\right)
    \end{align*}
    with $\Sigma_{\text{mMLE}} = \Sigma_2^{-1}$. Recall $\Gamma_1$ in (\ref{eq: Gamma1}) and $\Gamma_2$ in (\ref{eq: Gamma2}):
    \begin{align*}
        \Gamma_1 = \frac{1}{b} aa^\top \quad \text{and} \quad \Gamma_2 = \frac{b-c}{b^2} aa^\top.
    \end{align*}
    We have
    \begin{align*}
        \tilde{\Sigma}_{\text{SSL}} &= \left(\frac{1}{\lambda} \Gamma_1 + \Sigma_2\right)^{-1} \left(\frac{1}{\lambda} (\Gamma_1+\Gamma_2) +\Sigma_2\right) \left(\frac{1}{\lambda} \Gamma_1 +  \Sigma_2\right)^{-1} \\
        &= \left(\frac{1}{\lambda} \Gamma_1 + \Sigma_2\right)^{-1} + \frac{1}{\lambda}  \left(\frac{1}{\lambda} \Gamma_1 + \Sigma_2\right)^{-1} \Gamma_2  \left(\frac{1}{\lambda} \Gamma_1 + \Sigma_2\right)^{-1}.
    \end{align*}
    By the Sherman-Morrison formula, we can write
    \begin{align*}
        \left(\frac{1}{\lambda} \Gamma_1 + \Sigma_2\right)^{-1} = \left(\frac{1}{\lambda b} aa^\top + \Sigma_2\right)^{-1} = \Sigma_2^{-1} - \frac{\Sigma_2^{-1} aa^\top \Sigma_2^{-1}}{\lambda b+a^\top \Sigma_2^{-1} a}.
    \end{align*}
    It follows that
    \begin{align*}
          \left(\frac{1}{\lambda} \Gamma_1 + \Sigma_2\right)^{-1} \Gamma_2  \left(\frac{1}{\lambda} \Gamma_1 + \Sigma_2\right)^{-1} = \left(\frac{\lambda b}{\lambda b + a^\top \Sigma_2^{-1}a}\right)^2 \left(\frac{b-c}{b^2}\right)  \Sigma_2^{-1} aa^\top \Sigma_2^{-1}.
    \end{align*}
    Thus,
    \begin{align*}
        \tilde{\Sigma}_{\text{SSL}} &= \Sigma_2^{-1} - \frac{a^\top \Sigma_2^{-1}a + c\lambda}{\left(a^\top \Sigma_2^{-1}a + b\lambda\right)^2} \Sigma_2^{-1} aa^\top \Sigma_2^{-1} \\
        &= \Sigma_2^{-1} - \frac{d + c\lambda}{\left(d + b\lambda\right)^2} \Sigma_2^{-1} aa^\top \Sigma_2^{-1}. 
     \end{align*}
     For fixed $x$ we have
     \begin{align*}
         G^{\text{pred}}(x) &= \sqrt{\frac{\operatorname{AVar}\left(H\left(\hat{\beta}_{m}^\top x\right)\right)}{\operatorname{AVar}\left(H\left(\hat{\beta}_{n,m}^\top x\right)\right)}} \\
         &= \left(\frac{x^\top\Sigma_2^{-1}x}{x^\top \Sigma_2^{-1}x - \frac{d + c\lambda}{\left(d + b\lambda\right)^2} x^\top \Sigma_2^{-1} aa^\top \Sigma_2^{-1}x}\right)^{\frac{1}{2}} \\
         &=  \sqrt{\frac{1}{1 - \frac{(\tau+ \kappa \lambda)\tau}{(\tau+\lambda)^2} \cos^2{\theta_x}}} \\
         &\leq \frac{1+\frac{1}{\lambda} \tau}{\sqrt{1+\frac{1}{\lambda} \tau (2-\kappa)}}.
     \end{align*}
     Note that the last inequality holds with equality if and only if $x \in \operatorname{span}\{a\}$. As for $G^{\text{par}}$, we can write
     \begin{align*}
     G^{\text{par}} &= \sqrt{\frac{\det{\Sigma_{\text{mMLE}}}}{\det{\tilde{\Sigma}_{\text{SSL}}}}} \\
         &= \sqrt{\frac{\det{\Sigma_2^{-1}}}{\det{\left(\Sigma_2^{-1} - \frac{d + c\lambda}{\left(d + b\lambda\right)^2} \Sigma_2^{-1} aa^\top \Sigma_2^{-1}\right)}}} \\
         &=\frac{1}{\sqrt{\det{\left( \mathbbm{1}_{p\times p} -  \frac{d + c\lambda}{\left(d + b\lambda\right)^2} \Sigma_2^{-\frac{1}{2}}aa^\top \Sigma_2^{-\frac{1}{2}} \right)}}} \\
         &=\frac{1}{\sqrt{1- \frac{\left(d+c\lambda\right)d}{\left(d + b\lambda\right)^2}}} \\
         &= \frac{1+\frac{1}{\lambda} \tau}{\sqrt{1+\frac{1}{\lambda} \tau (2-\kappa)}}.
     \end{align*}
\end{proof}

\begin{rem}
Note that $G^{\text{par}} \geq 1$. In fact 
\begin{align*}
    G^{\text{par}} = \max_x G^{\text{pred}}(x) \geq \sqrt{\frac{1}{1 - \frac{(\tau+ \kappa \lambda)\tau}{(\tau+\lambda)^2} \cos^2{\frac{\pi}{2}}}} = 1.
\end{align*}
\end{rem}

\begin{rem}
 Note that the parameter $\kappa$ is closely related to the signal-to-noise ratio in the unmatched sample. In fact, if this signal-to-noise ratio is denoted by $\widetilde{\text{SNR}}$ then
\begin{eqnarray*}
    \widetilde{\text{SNR}} = \frac{\operatorname{Var}\big(\tH(\beta_0^\top \tX)\big)}{\mathbb{E}\left[\tH(\beta_0^\top \tX) (1-\tH(\beta_0^\top \tX))\right]} = \frac{\operatorname{Var}(\tY) -  \mathbb{E}\Big[\tH(\beta_0^\top \tX) (1-\tH(\beta_0^\top \tX))\Big]}{\mathbb{E}\left[\tH(\beta_0^\top \tX) (1-\tH(\beta_0^\top \tX))\right]} =   \frac{1}{\kappa} - 1,
\end{eqnarray*}
since 
\begin{eqnarray*}
\operatorname{Var}(\tY)  & =   &     \mathbb{E}\Big(\tH(\beta_0^\top \tX) \Big) \Big( 1-   \mathbb{E}\Big(\tH(\beta_0^\top \tX) \Big) \Big)  \\
& =  &   \mathbb{E}\Big(\tH(\beta_0^\top \tX) \Big) - \mathbb{E}\Big(\tH(\beta_0^\top \tX) \Big)^2 \\
& = &  \mathbb{E}\Big(\tH(\beta_0^\top \tX) \Big) - \mathbb{E}\Big(\tH(\beta_0^\top \tX)^2 \Big) + \operatorname{Var}\big(\tH(\beta_0^\top \tX)\big).
\end{eqnarray*}

\end{rem}

\section{The gain in some specific settings} 

\subsection{A common link function and covariate distribution} \label{Section: Matched and Unmatched Data with a Common Link Function and Covariate Distribution}
In this section, we will focus on the special case where 
\begin{eqnarray}\label{common}
\mathbb{P}^{X} = \mathbb{P}^{\tX}, \ \  \text{and} \ \   \tH = H
\end{eqnarray}
where $H$ is some known link function that satisfies the assumptions under which the asymptotic gain was derived. 

Recall the quantities $\tau$ and $\kappa$ defined in Theorem \ref{theorem: Gain}. Under (\ref{common}) it holds that $\tau = \kappa \,\psi$, 
where 
\begin{align*}
    \psi = \frac{\mathbb{E}[UZ]^\top \mathbb{E}[UU^\top]^{-1}\mathbb{E}[UZ]}{\mathbb{E}[Z^2]}
\end{align*}
with
\begin{align*}
    Z := \sqrt{H(\beta_0^\top X)(1-H(\beta_0^\top X)) }\quad \text{and} \quad
    U := \frac{H'(\beta_0^\top X)}{\sqrt{H(\beta_0^\top X) (1-H(\beta_0^\top X))}} X.
\end{align*}
If $H$ is the logistic link, then using the identity $H' =  H (1-H)$ we can show that $UZ = Z^2X$ and $UU^\top = Z^2 XX^\top$. Now, define $\mathbb Q$ to be the probability measure that satisfies
\begin{align*}
    d\mathbb{Q}(x) = \frac{Z(x)^2}{\mathbb{E}(Z^2)}d\mathbb{P}^X(x).
\end{align*}
Then,
\begin{align*}
    \mathbb{E}(UZ) = \mathbb{E}(Z^2)\mathbb{E}\left(\frac{Z^2}{\mathbb{E}(Z^2)}X\right) = \mathbb{E}(Z^2)\mathbb{E}_{\mathbb{Q}}(X), \ \ \text{and}  \ \ \  \mathbb{E}(UU^\top) = \mathbb{E}(Z^2)\mathbb{E}_{\mathbb{Q}}\left(XX^\top\right).
\end{align*}
Consequently,
\begin{align*}
  \psi = \frac{\mathbb{E}(UZ)^\top \mathbb{E}(UU^\top)^{-1}\mathbb{E}(UZ)}{\mathbb{E}(Z^2)} = \mathbb{E}_{\mathbb{Q}}(X)^\top \left(\mathbb{E}_{\mathbb{Q}}\left(XX^\top\right)\right)^{-1} \mathbb{E}_{\mathbb{Q}}(X) = \frac{d_{\mathbb Q}}{1+d_{\mathbb Q}},
\end{align*}
with $d_{\mathbb Q} = \mu_{\mathbb Q}^\top \Sigma_{\mathbb Q}^{-1}\mu_{\mathbb Q}$ the squared Mahalanobis distance, with
   $ \mu_{\mathbb Q} : = \mathbb{E}_{\mathbb{Q}}(X)$ and $\Sigma_{\mathbb Q} := \operatorname{Var}_{\mathbb{Q}}\left(X\right)$.
 Thus, the gain is determined by the abundance of unmatched data $\lambda$, the unmatched signal-to-noise-ratio-related quantity $\kappa$, and $\tau$ which encodes aspects of the underlying information geometry.

\begin{figure}[htbp]
    \centering
    \includegraphics[width=0.8\textwidth]{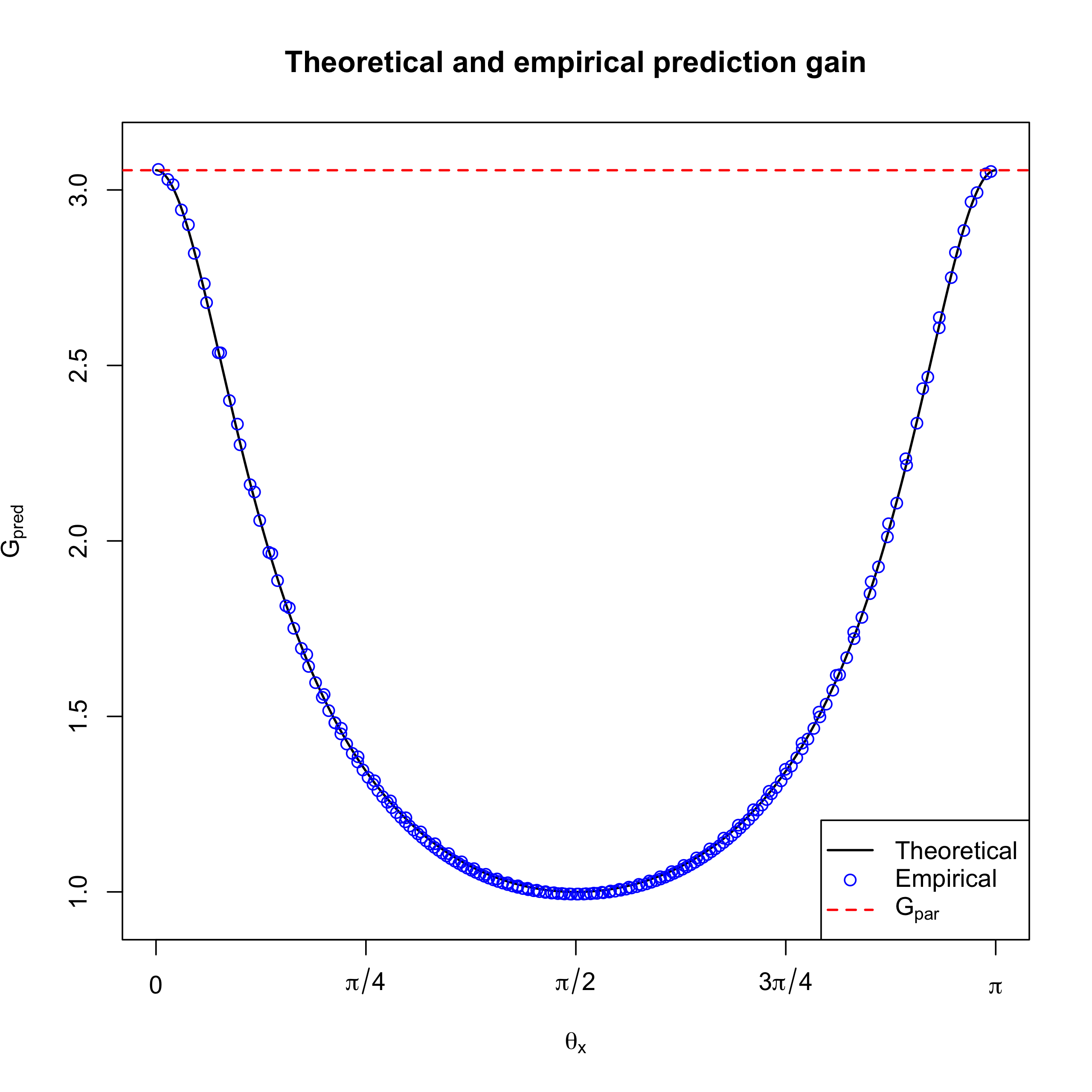}
    \caption{Theoretical and empirical prediction gains $G^{\mathrm{pred}}(x)$ for the logistic link function, with $ \beta_0= \begin{pmatrix}
    1\\
    -1/2
    \end{pmatrix}$, $X \overset{d}{=} \tX$, $ X\sim\mathcal{N}(\mu,\Sigma)$, where $\mu= (1, 1/2)^\top $ and $ \Sigma= \begin{pmatrix}
    1 & 1/2\\
    1/2 & 3/2 \end{pmatrix}$.
    We take $m=500$ and $n=12500$ and $\lambda=0.04$ . The theoretical prediction gain curve is approximated by Monte Carlo method. The empirical prediction gain is evaluated at 200 points on the unit circle $ x_i=    (\cos\phi_i, \sin\phi_i)^\top$,  $\phi_i\in[0,2\pi],$ with the $\phi_i$ equally spaced over $[0,2\pi]$.}
    \label{plot: prediction gain X Gaussian}
\end{figure}

In Figure \ref{plot: prediction gain X Gaussian}, We give a first illustration of the gain in a setting where we take the link function to be logistic and both $X$ and $\tilde X$ are distributed as $\mathcal N(\mu, \Sigma)$ with $\mu = (1, 1/2)^\top$ and $\displaystyle \Sigma= \begin{pmatrix}
    1 & 1/2\\
    1/2 & 3/2 \end{pmatrix}$. The solid curve in black  shows the theoretical predictive gain $G^{\text{pred}}(x)$ as a function of the angle $\theta_x$ as defined   Theorem \ref{theorem: Gain} for $\lambda = 0.04$. The blue dots depict the empirical gain obtained for 200 points on the unit circle. They show very good agreement with the theoretical values. Using a Monte Carlo approach, the empirical gain values were computed using independent matched samples of size $m = 500$ and unmatched samples of size $n = 12500$ drawn from the Gaussian distribution with the aforementioned parameters $\mu$ and $\Sigma$.

\subsubsection{Example with no gain}
If $ H(t) = 1-H(-t)$, as in the case for both the logistic and probit link functions, then $H'$ is an even function. In this case, there is no statistical gain if the covariate $X$ satisfies $X \overset{d}{=} -X$.  To see this, note that
\begin{align*}
     &a = \mathbb{E}\left[ H'(\beta_0^\top X) X  \right] = \mathbb{E}\left[ H'(-\beta_0^\top X) (-X)  \right] = -\mathbb{E}\left[ H'(\beta_0^\top X) X  \right] = 0.
\end{align*}
Thus $\tau = d = 0$, which implies that $G^{\text{par}}=1$. 

\subsubsection{Example with \texorpdfstring{$\beta_0 = 0$}{beta0 = 0} }
\begin{figure}[htbp]
    \centering
    \includegraphics[width=0.8\textwidth]{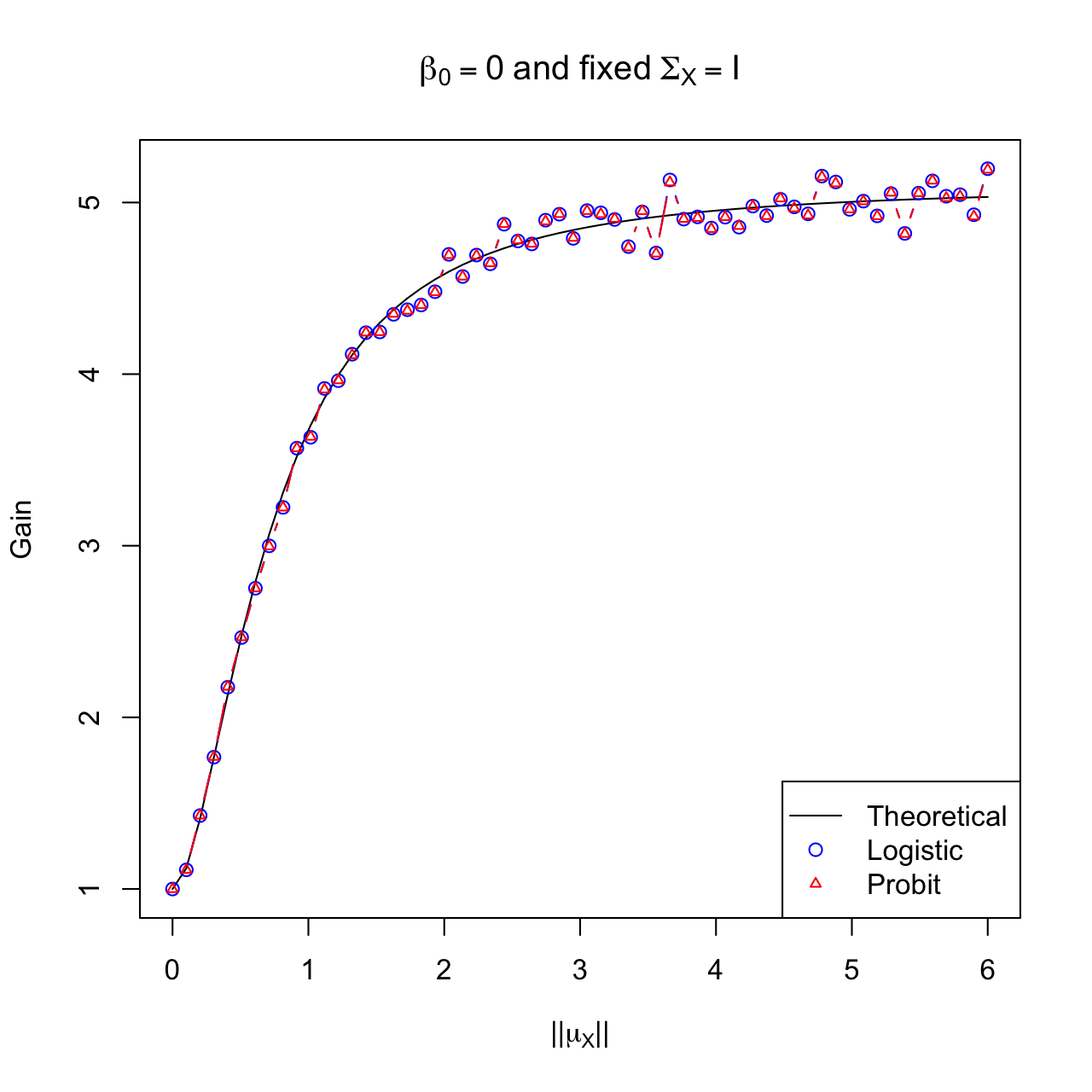}
    \caption{Theoretical and empirical parameter estimation gains $G^{\mathrm{par}}$ for logistic and probit link functions, with $\beta_0=(0,0)^\top$, $X \overset{d}{=} \tX$ and $X\sim \mathcal{N}(\mu_X,\mathbbm{1}_{2\times 2})$. The covariate mean $\mu_X$ varies along a spiral, $\mu_X=r(\cos{\theta},\sin{\theta})^\top$, where $r$ increases from 0 to 6 and $\theta$ increases from 0 to $6\pi$ over 60 equally spaced configurations. The sample sizes are $m=5000$ and $n=125000$, with $4000$ Monte Carlo replications per configuration. The solid curve represents the theoretical gain, while he blue circles and red triangles represent the empirical gains for logistic and probit links, respectively.}
    \label{plot: Gain_beta=0}
\end{figure}

Suppose that $\beta_0 = 0$. Using the notation of Theorem \ref{theorem: Gain}, we have 
\begin{align*}
    &a = H'(0)\mathbb{E}(X),\\
    &b = c =  H(0)\left(1-H(0)\right),\\
    &\Sigma_2 = \frac{\left(H'(0)\right)^2}{H(0)\left(1-H(0)\right)} \mathbb{E}(XX^\top),\\
    &d = a^\top \Sigma_2^{-1}a = H(0)\left(1-H(0)\right) \mathbb{E}(X)^\top \mathbb{E}(XX^\top)^{-1}\mathbb{E}(X),\\
    &\kappa = \frac{c}{b} = 1,\\
    &\tau = \frac{d}{b} =  \mathbb{E}(X)^\top \mathbb{E}(XX^\top)^{-1}\mathbb{E}(X).
\end{align*}
Put $\mu_X = \mathbb{E}(X)$ and $\Sigma_X = \operatorname{Cov}(X)$. Then by the Sherman–Morrison formula it follows that
\begin{align*}
     \mathbb{E}(X)^\top \mathbb{E}(XX^\top)^{-1}\mathbb{E}(X) &= \mu_X^\top \left(\Sigma_X+\mu_X\mu_X^\top\right)^{-1}\mu_X \\
     &= \mu_X^\top \left(\Sigma_X^{-1} - \frac{\Sigma_X^{-1}\mu_X\mu_X^\top \Sigma_X^{-1} }{1+ \mu_X^\top \Sigma_X^{-1}\mu_X}\right)\mu_X   \\
     &= \frac{\mu_X^\top\Sigma_X^{-1}\mu_X}{1+\mu_X^\top\Sigma_X^{-1}\mu_X}.
\end{align*}
Thus
\begin{align}\label{eq: Gain with beta = 0}
    G^{\text{par}} = \sqrt{1+\frac{1}{\lambda} \tau} =\sqrt{1+\frac{1}{\lambda}\frac{1}{1+\rho}} \quad \text{with} \quad \rho = \frac{1}{\mu_X^\top\Sigma_X^{-1}\mu_X}.
\end{align}
This result is illustrated in Figure \ref{plot: Gain_beta=0}.  Note that if $\Sigma_X  = \mathbbm{1}_{p \times p}$, then
\begin{eqnarray*}
G^{\text{par}} =\sqrt{1+\frac{1}{\lambda} \frac{\Vert \mu_X \Vert^2}{\Vert \mu_X \Vert^2 + 1}}.
\end{eqnarray*}
In this case,  $G^{\text{par}}$ converges to its maximal value $\displaystyle \sqrt{1 + \frac{1}{\lambda}}$ as $\Vert \mu_X \Vert$ gets larger. For smaller values of $\Vert \mu_X \Vert$, the gain tends to $1$. \\ 

It might be surprising that the gain depends on $\Vert \mu_X \Vert$ even though $\beta_0 = 0$. Mathematically, the gain is expected to depend on this quantity since when $\beta_0 = 0$ the only relevant quantity is $\tau$, which involves a scaled version of $\Vert \mu_X \Vert^2$. Here, we would like to give an intuitive explanation for what might be seen as a non-intuitive finding. Recall that we are assuming that $X \overset{d}{=} \tX$. Imagine now that the true $\beta_0$ results from making a small perturbation around $0$ and let us denote this perturbation by $\delta$. By Taylor expansion, we have
\begin{align*}
    \mathbb{E}(\tY) = \mathbb{P}(\tY = 1) = \int H(\epsilon^\top x) d\mathbb{P}^{X}(x) = H(0) + H'(0) \delta^\top \mu_X + o(\norm{\epsilon}).
\end{align*}
So we can write 
\begin{align*}
    \delta^\top \mu_X \approx \frac{1}{H'(0)} \left(\mathbb{E}(\tY) -H(0)\right).
\end{align*}
If $\delta \not\perp \mu_X$, then  $\delta^\top \mu_X = \norm{\delta}\norm{\mu_{X}} \cos{\alpha}$ with $\cos{\alpha} \neq 0$. It follows that
\begin{align*}
    \norm{\delta} \approx \frac{1}{H'(0) \cos{\alpha}} \frac{\mathbb{E}(\tY) -H(0)}{\norm{\mu_{X}}}.
\end{align*}
Heuristically, if everything is fixed and only $\Vert \mu_X \Vert$ is allowed to change, then one sees that 
large values of $\norm{\mu_{X}}$ give  stronger evidence that $\delta \neq 0$ and vice versa. Hence, $\norm{\mu_{X}}$ is expected to have an impact on the statistical precision when estimating $\beta_0$ in the setting where it is equal to $0$. In other words, when $\Vert \mu_X \Vert$ is large, we should be more confident that $\beta_0 =0$ since the confidence sets are smaller. In fact, this is just another way to say that we have a greater statistical gain.

\bigskip

\subsubsection{Example with Logistic link and a uniform covariate}
\begin{figure}[htbp]
    \centering
    \includegraphics[width=0.8\textwidth]{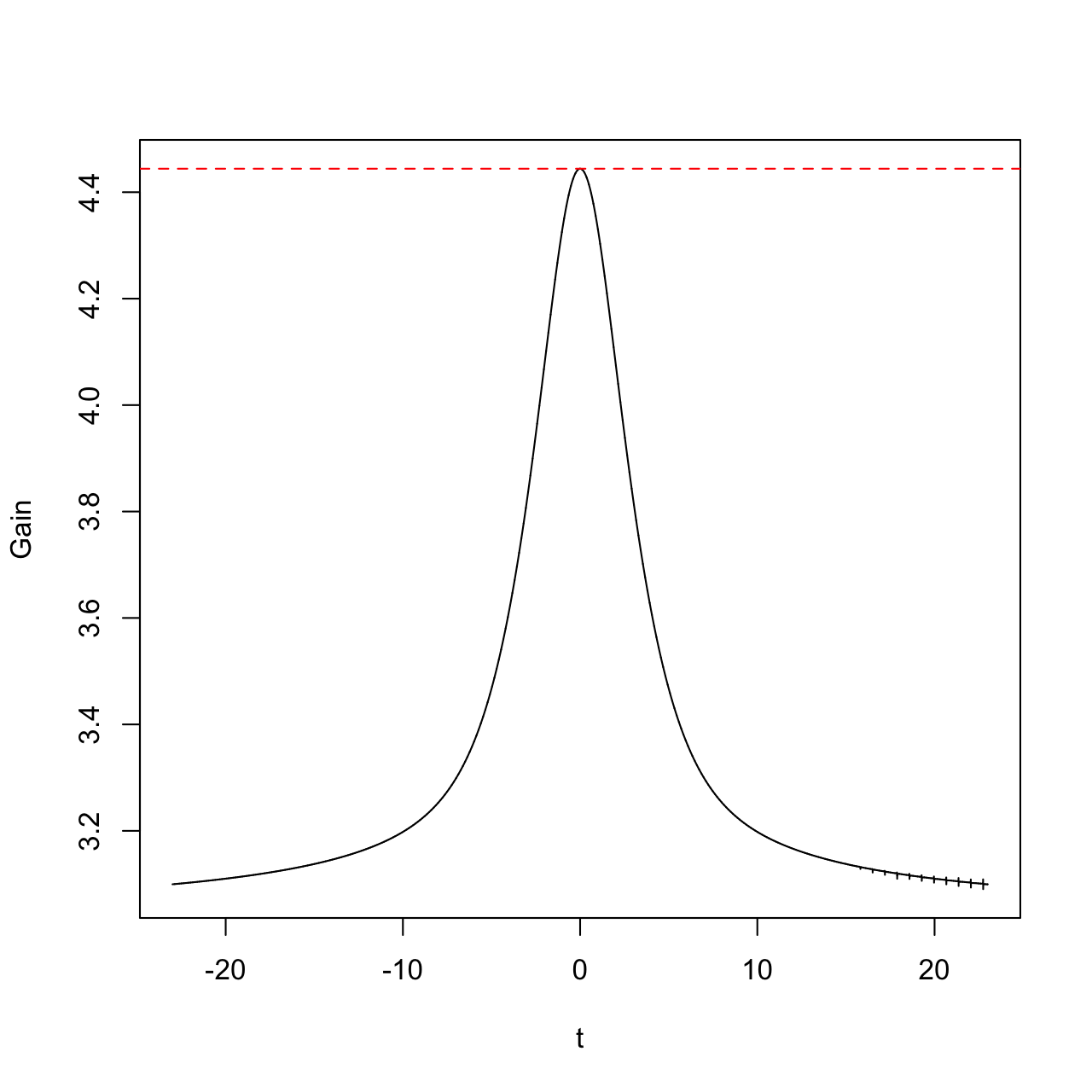}
    \caption{Theoretical gain for $H$ being the logistic link function and $X \sim \mathrm{U}[0,r]$, with $\lambda=0.04$. $t = \beta_0 r$. The gain is maximized at $t=0$, and its maximum value can be obtained from (\ref{eq: Gain with beta = 0}).}
    \label{plot: Gain Uniform X}
\end{figure}

Take again the link function $H$ to be the Logistic function, 
and consider the case $ X \sim \mathrm{U}[0,r]$.  Then, it is clear that the model remains unchanged if we replace $X$ by $U \sim  \mathrm{U}[0,1]$ and $\beta_0$ by $\beta_0 r$.  Put $t = \beta_0 r$, and define
$$
 I_1(t)=\int_0^1 H(tu)\,du, \ \  I_2(t)=\int_0^1 u H'(tu)\,du, \ \  I_3(t)=\int_0^1 H(tu)^2\,du
 $$
 and $I_4(t)=\int_0^1 u^2H'(tu)\,du$.  Using the expression of $H$, it can be shown that
\begin{eqnarray*}
        I_1(t)  =   
      \begin{cases}
      \frac{1}{t} \log{\left(\frac{1+\exp{\left(t\right)}}{2}\right)}, \ \text{if $t \ne 0$} \\
      \frac{1}{2}, \ \ \ \hspace{2cm} \text{if $t = 0$},
      \end{cases} 
\end{eqnarray*}
\begin{eqnarray*}
   I_2(t)  =     
    \begin{cases}
    \frac{1}{t} H(t) - \frac{1}{t^2}\log{\left(\frac{1+\exp{\left(t\right)}}{2}\right)}, \ \ \text{if $t \ne 0$} \\
   \frac{1}{8}, \ \ \  \ \hspace{3.5cm} \text{if $t = 0$},
      \end{cases} 
   \end{eqnarray*}
\begin{eqnarray*}
     I_3(t) =   
     \begin{cases} 
     \frac{1}{t} \left( \log{\left(\frac{1+\exp{\left(t\right)}}{2}\right)} -H(t) + \frac{1}{2}\right),\ \ \text{if $t \ne 0$} \\
   \frac{1}{4}, \ \ \  \ \hspace{4.2cm} \text{if $t = 0$}, 
     \end{cases}
\end{eqnarray*}
and
\begin{eqnarray*}
       I_4(t) =
       \begin{cases}
       \frac{1}{t} H(t) - \frac{2}{t^2} \log{\left(1+\exp{\left(t\right)}\right)} - \frac{2}{t^3} \operatorname{Li}_2(-\exp{\left(t\right)}) - \frac{\pi^2}{6 t^3}, \  \text{if $t \ne 0$} \\
        \frac{1}{12}, \ \ \  \ \hspace{7.2cm} \text{if $t = 0$} \\ 
        \end{cases} 
\end{eqnarray*}
where
 $\operatorname{Li}_2(x) = -\int_0^x \frac{\log{\left(1-t\right)}}{t}\,dt$ is the di-logarithmic function.
 Furthermore, using a similar notation as in Theorem \ref{theorem: Gain}, we have
\begin{align*}
    &a(t) = r  I_2(t), \\
    &b(t) = I_1(t) - I_1(t)^2, \\
    &c(t) = I_1(t) - I_3(t) ,\\
    &\Sigma_2(t) = r^2I_4(t) ,\\
    & d(t) = a(t)^\top \Sigma_2^{-1}(t) a(t) = \frac{I_2(t)^2}{I_4(t)},\\
    & \kappa(t) = \frac{I_1(t) -I_3(t)}{ I_1(t) - I_1(t)^2} ,\\
    & \tau(t) =  \frac{I_2(t)^2}{ I_4(t)\left( I_1(t) - I_1(t)^2 \right)}.
\end{align*}
By Theorem \ref{theorem: Gain}, the gain of using the SLLEMLE $\widehat \beta_{n, m}$ over the matched MLE $\widehat  \beta_m$  is 
\begin{align*}
    G_t^{\text{par}} &= \frac{1+\frac{1}{\lambda} \tau(t)}{\sqrt{1+\frac{1}{\lambda} \tau(t) (2-\kappa(t))}}.
\end{align*}
Figure \ref{plot: Gain Uniform X} plots the gain $\lambda = 0.04$, that is, the function 
$$
t \mapsto G_{t}^{\text{par}}.
$$
As the plot shows, the gain is unimodal with mode $G_0^{\text{par}}$ and stabilizes around  the asymptotic value $3.03$. Although it is hard to theoretically prove unimodality, we can show that $G$ has to decrease in the neighborhood of $0$. First, we can use the identities $H(-x) = 1 - H(x)$ and $H'(-x) = H'(x)$ to show that $G_t^{\text{par}} = G_{-t}^{\text{par}}$ and hence the derivative $\displaystyle \frac{\partial G_t^{\text{par}}}{\partial t}|_{t = 0} = 0 $. Furthermore,  computing the odd degree derivatives of $H$ at $0$ up to the third degree yields 
\begin{eqnarray*}
        H(x) =
        \frac{1}{2}+\frac{x}{4}-\frac{x^3}{48}+O(x^5),
        \qquad x \to 0,
\end{eqnarray*}
and hence
\begin{eqnarray*}
        H'(x) =
        \frac{1}{4}-\frac{x^2}{16}+\frac{x^4}{96}+O(x^6),
        \qquad x\to 0.
\end{eqnarray*}
Then, we obtain
We obtain
\begin{eqnarray*}
        I_1(t) =
    \frac{1}{2}+\frac{t}{8}-\frac{t^3}{192}+O(t^5),  \ I_2(t)=
 \frac{1}{8}-\frac{t^2}{64}+\frac{t^4}{576}+O(t^6),
\end{eqnarray*}
and 
\begin{eqnarray*}
        I_4(t) = 
        \frac{1}{12}-\frac{t^2}{80}+\frac{t^4}{672}+O(t^6).
\end{eqnarray*}
For $I_3(t)$, the logistic identity $ H'(x)=H(x) (1-H(x))$ implies that  $H(x)^2=H(x)-H'(x)$, and therefore
\begin{eqnarray*}
    I_3(t)  =  I_1(t)-\int_0^1 H'(tu)\,du =    I_1(t)-\frac{H(t)-H(0)}{t}.
\end{eqnarray*}
Using the Taylor expansion of $H$, this gives
\begin{eqnarray*}
        I_3(t)=
        \frac{1}{4}+\frac{t}{8}+\frac{t^2}{48}-\frac{t^3}{192}+O(t^4).
\end{eqnarray*}
Consequently,
\begin{eqnarray*}
I_1(t)(1-I_1(t))
=  \frac{1}{4}-\frac{t^2}{64}+O(t^4), \ \  \text{and} \ \   I_3(t)-I_1(t)^2=
        \frac{t^2}{192}+O(t^4).
\end{eqnarray*}
It follows that
\begin{eqnarray*}
\tau (t)=\frac{I_2(t)^2}{I_1(t)(1-I_1(t))I_4(t)} =\frac{3}{4}-\frac{9}{320}t^2+O(t^4),
\end{eqnarray*}
and
\begin{eqnarray*}
   \tau (t) \left(2-\kappa(t)\right)
        &=&
        \frac{I_2(t)^2}{I_4(t)}
        \left[
        \frac{1}{I_1(t)(1-I_1(t))}
        +
        \frac{I_3(t)-I_1(t)^2}{I_1(t)^2 (1-I_1(t))^2}
        \right] \\
        &=  &  \frac{3}{4}-\frac{1}{80}t^2+O(t^4).
\end{eqnarray*}
Now, note that
\begin{eqnarray*}
        (G_t^{\text{par}})^2= \frac{(1+\lambda^{-1}\tau(t))^2}{1+\lambda^{-1}\tau(t) \left(2-\kappa(t)\right)},
\end{eqnarray*}
and hence when $t \sim 0$
\begin{eqnarray*}
 (G_t^{\text{par}})^2= 1+\frac{3}{4\lambda} - \frac{7}{160\lambda}t^2 +O(t^4).
\end{eqnarray*}
Therefore,
\begin{eqnarray*}
G_t^{\text{par}} = \sqrt{1+\frac{3}{4\lambda}}  -  \frac{7}{320\lambda \sqrt{1+\frac{3}{4\lambda}}}t^2
        +O(t^4).
\end{eqnarray*}
This shows that for a fixed $\lambda$, the gain is unimodal in the neighborhood of the origin with mode  $G_0^{\text{par}} =  \sqrt{1+\frac{3}{4\lambda}}$.  
\begin{lemma} \label{lemma: G_infinity X uniform H logis} 
If $H$ is logistic link function, i.e., $H(z) = 1/(1+\exp{(-z)})$ and $X \sim \mathrm{U}[0,r]$. Write $t := \beta_0r$, then we have
    \begin{eqnarray}\label{Ginfty}
\lim_{ \vert t \vert \to \infty}  G_t^{\text{par}} =   \frac{
        1+\lambda^{-1}\dfrac{6\log 2}{\pi^2}
        }
        {\sqrt{
        1+\lambda^{-1}\dfrac{12\log 2-3}{\pi^2}
        }} := G_\infty(\lambda).
\end{eqnarray}
\end{lemma}
The proof and calculation are provided in Appendix \ref{appendix: Uniform covariate example}. Note that $G_\infty(0.04) \approx 3.03$, which can be obtained from the general formula above. 


\subsection{Matched and unmatched data under differently scaled \texorpdfstring{$H$}{H}}
\begin{figure}[htbp]
    \centering
    \includegraphics[width=0.8\textwidth]{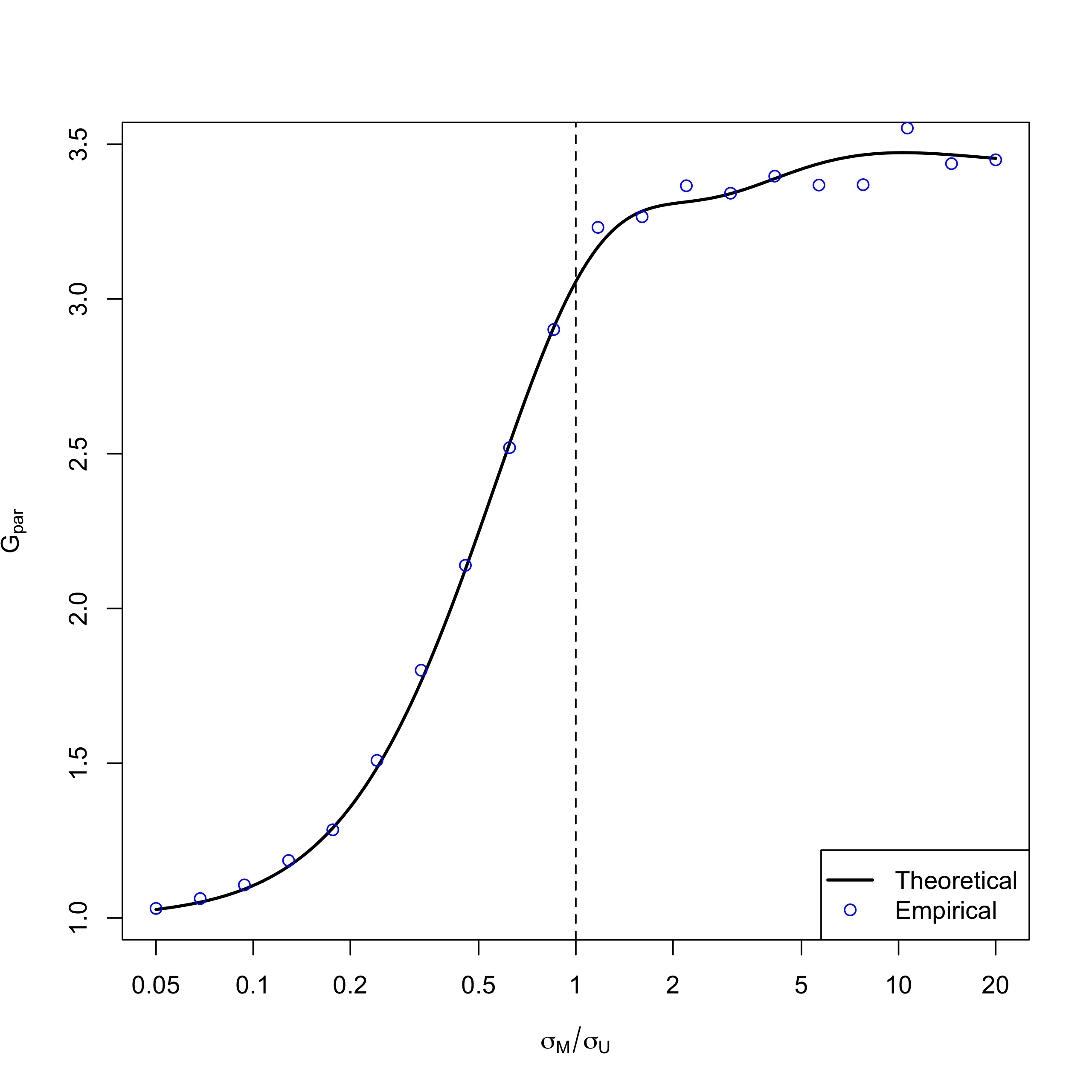}
    \caption{Theoretical and empirical parameter gains $G^{\mathrm{par}}$ for the logistic link function, with $ \beta_0= \begin{pmatrix}
    1\\
    -1/2
    \end{pmatrix}$, $ X\sim\mathcal{N}(\mu,\Sigma)$, where $\mu= \begin{pmatrix}
    1\\
    1/2
    \end{pmatrix}$ and $ \Sigma= \begin{pmatrix}
    1 & 1/2\\
    1/2 & 3/2 \end{pmatrix}$.
    We set $\lambda=0.04$ and $m=500$, so that $n=12500$. The theoretical gain curve is approximated by Monte Carlo method. The model is given in (\ref{eq: model with different scale}) and we use $\sigma_M = 1$ in the empirical simulation.}
    \label{plot: parameter gain X Gaussian}
\end{figure}

\begin{figure}[h]
    \includegraphics[width=0.48\textwidth]{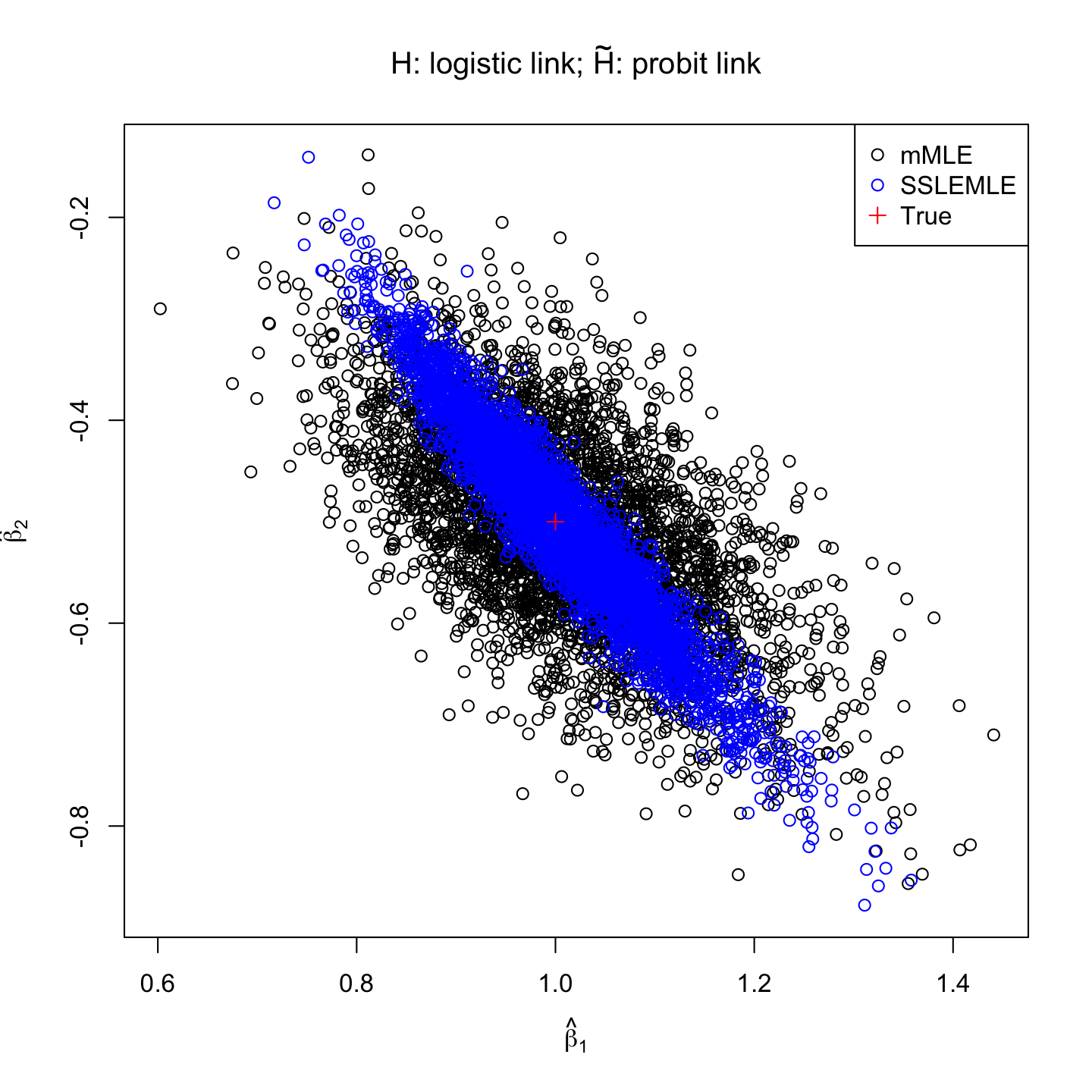}
    \includegraphics[width=0.48\textwidth]{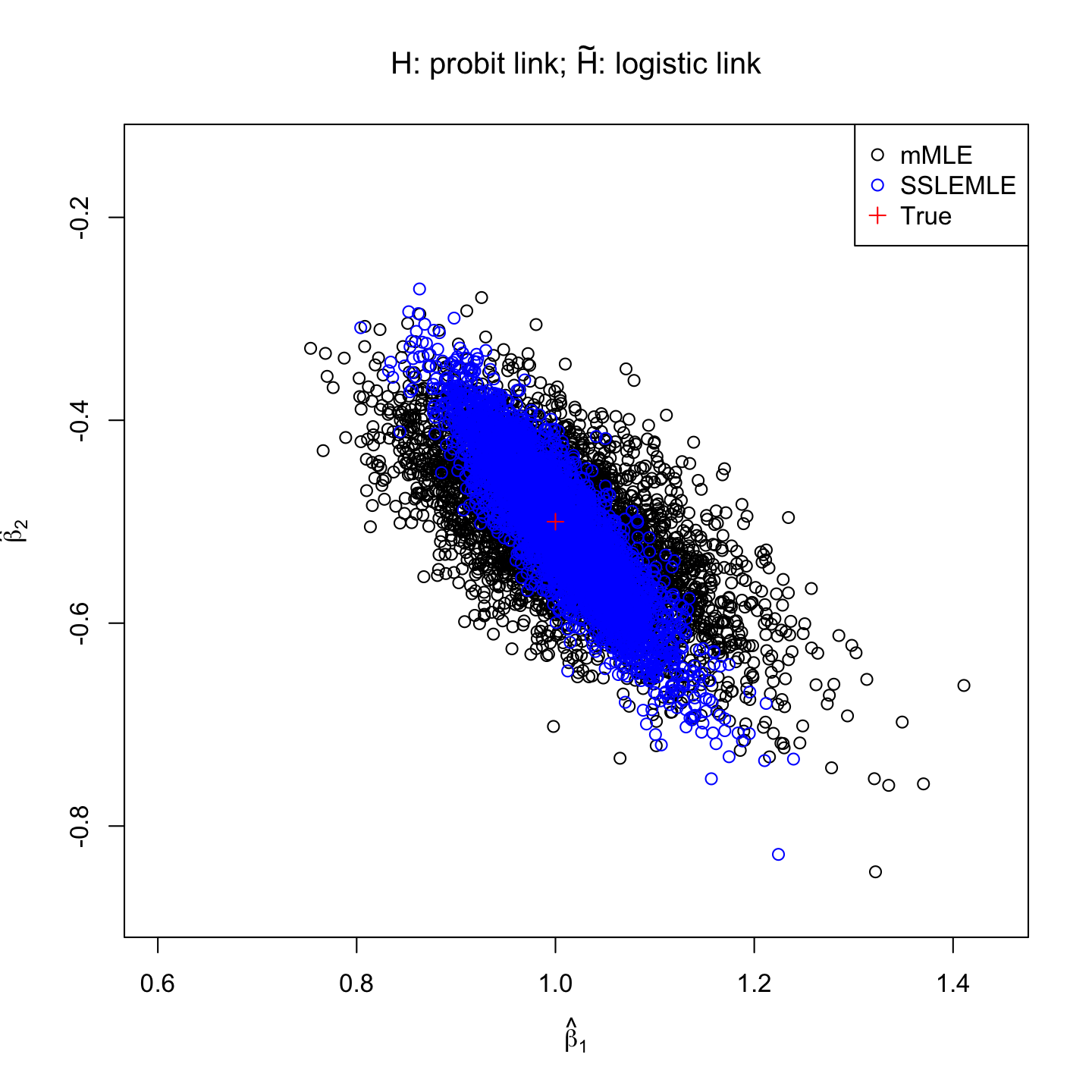}
    \caption{Matched MLEs and SSLEMLEs with $ \beta_0= \begin{pmatrix}
    1\\
    -1/2
    \end{pmatrix}$, $X \overset{d}{=} \tX$ and $ X\sim\mathcal{N}(\mu,\Sigma)$, where $\mu= \begin{pmatrix}
    1\\
    1/2
    \end{pmatrix}$ and $ \Sigma= \begin{pmatrix}
    1 & 1/2\\
    1/2 & 3/2 \end{pmatrix}$. The link function $H$ for matched data and the link function $\tH$ for unmatched data are different. We set $\lambda=0.04$ and $m=500$, so that $n=12500$.}
    \label{plot: Different H and tH}
\end{figure}

In the previous section, we considered the setting in which the matched and unmatched data share the same link function $H$. We now study the case where the link functions $H$ and $\tH$ belong to the same family but have different scale parameters. Here we still assume that $\mathbb{P}^{X} = \mathbb{P}^{\tX} $. Specifically,
\begin{align}\label{eq: model with different scale}
    &\mathbb{P}(Y = 1 \mid X ) = H\left(\beta_0^\top X\right) = H_1\left(\frac{\beta_0^\top X}{\sigma_M}\right), \notag \\
    &\mathbb{P}(\tY = 1 \mid \tX ) = \tH\left(\beta_0^\top X\right)  = H_1\left(\frac{\beta_0^\top X}{\sigma_U}\right)  = H\left(\beta_0^\top \left(\frac{\sigma_M}{\sigma_U}X\right)\right),
\end{align}
where $H_1$ denotes the link function with scale 1, and $\sigma_M,\sigma_U \in (0, \infty)$ are the scale parameters for the matched and unmatched data, respectively. We assume that $H_1$ (and hence $H$), $\sigma_M$ and $\sigma_U$ are known. (\ref{eq: model with different scale}) implies that this problem can be rewritten with $H = \tH$ and $\displaystyle \tX = \frac{\sigma_M}{\sigma_U} X$. We illustrate the statistical gain in this case with an example shown in Figure \ref{plot: parameter gain X Gaussian}. As shown in the plot, the gain tends to be larger when $\sigma_M/\sigma_U$ becomes larger, meaning that the unmatched link $\tH$ is sharper than the matched link $H$. Conversely, interchanging the two links tends to reduce the gain, although this gain remains greater or equal than 1. We expect a similar pattern in other settings where $H$ and $\tH$ are not equal. Indeed, Figure \ref{plot: Different H and tH} shows that the gain is larger when $H$ is the logistic link and $\tH$ is the probit link than in the reversed setting. This is expected and consistent with the fact that the probit link is steeper near the origin.

\bigskip

\section{Real data application}
We use the Dry Bean data set from the UCI Machine Learning Repository, originally introduced in \cite{KOKLU2020105507}. The data contain 13,611 bean samples from 7 dry bean varieties. Images were taken with a high-resolution camera, and 16 numerical features were extracted, including 12 size-related measurements and 4 shape-related features. For this data set, the original task is to classify the 7 bean varieties. In this paper, we focus on the first variety, Seker, and construct a binary outcome. Specifically, the response variable is defined as $Y=1$ if the bean is Seker and $Y=0$ otherwise. For simplicity, we use ``ShapeFactor1'' and ``ShapeFactor2'' as covariates. First, we used all 13,611 observations to fit a logistic regression model without an intercept, with the response variable defined as above. We multiplied both covariates by 10,000 so that the estimated parameters are on a smaller scale, which helps stabilize the numerical optimization. As shown in Table \ref{Table: real data reference model}, we used these results as a reference for the following analysis.
\begin{table}[ht]
\centering
\begin{tabular}{lrrrr}
\hline
Covariate & Estimate & Std. Error & $z$ value & $p$-value \\
\hline
10,000$\times$ShapeFactor1 & -0.3259 & 0.0077 & -42.42 & $<2\times10^{-16}$ \\
10,000$\times$ShapeFactor2 &  0.9299 & 0.0225 &  41.27 & $<2\times10^{-16}$ \\
\hline
\end{tabular}
\caption{Logistic regression results for the Dry Bean data. The residual deviance is 2653.7, the null deviance is 18868.9, the AIC is 2657.7, and the in-sample AUC is 0.9854.} \label{Table: real data reference model}
\end{table}

\begin{figure}[htbp]
    \centering
    \includegraphics[width=0.8\textwidth]{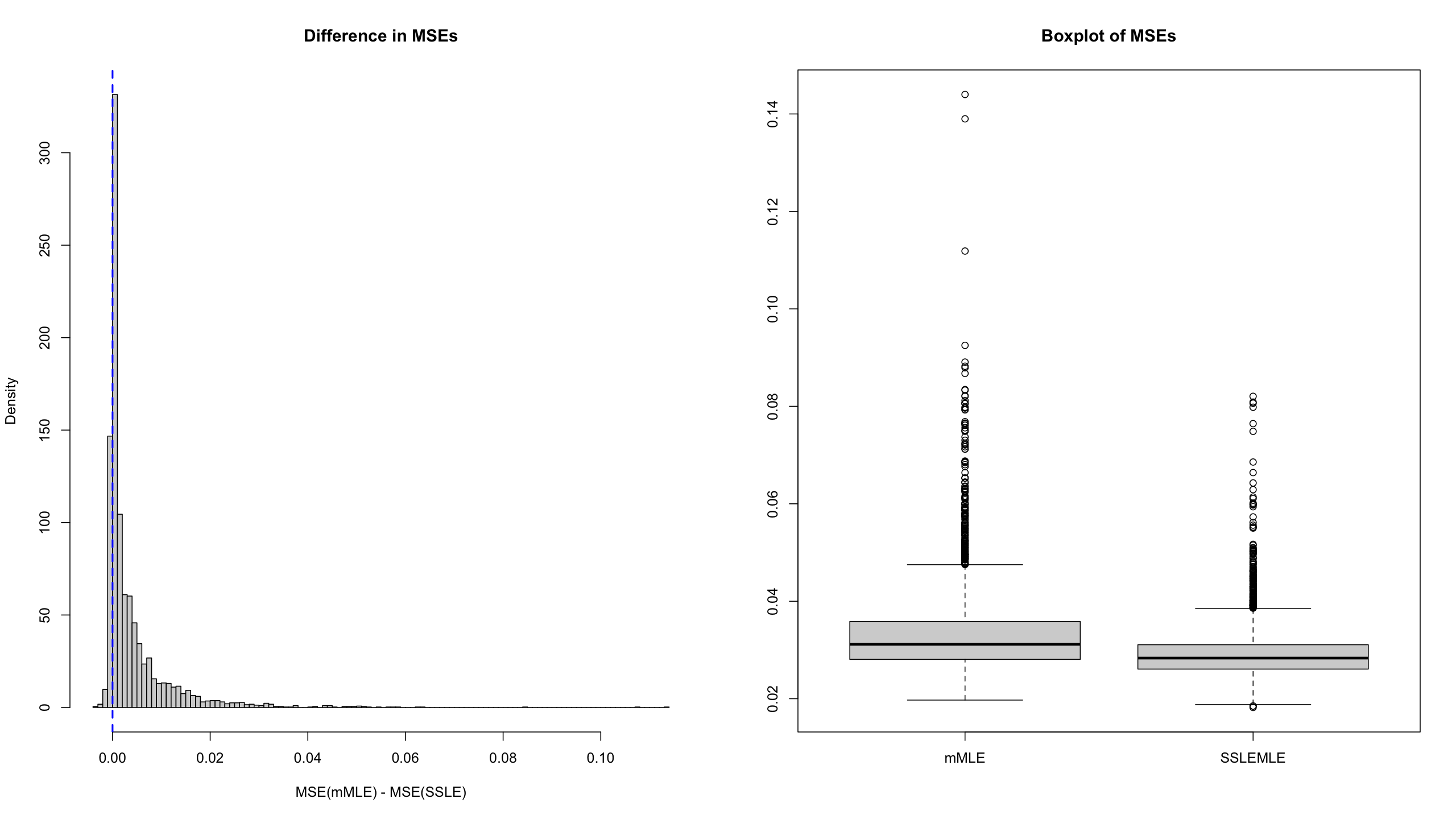}
    \includegraphics[width=0.8\textwidth]{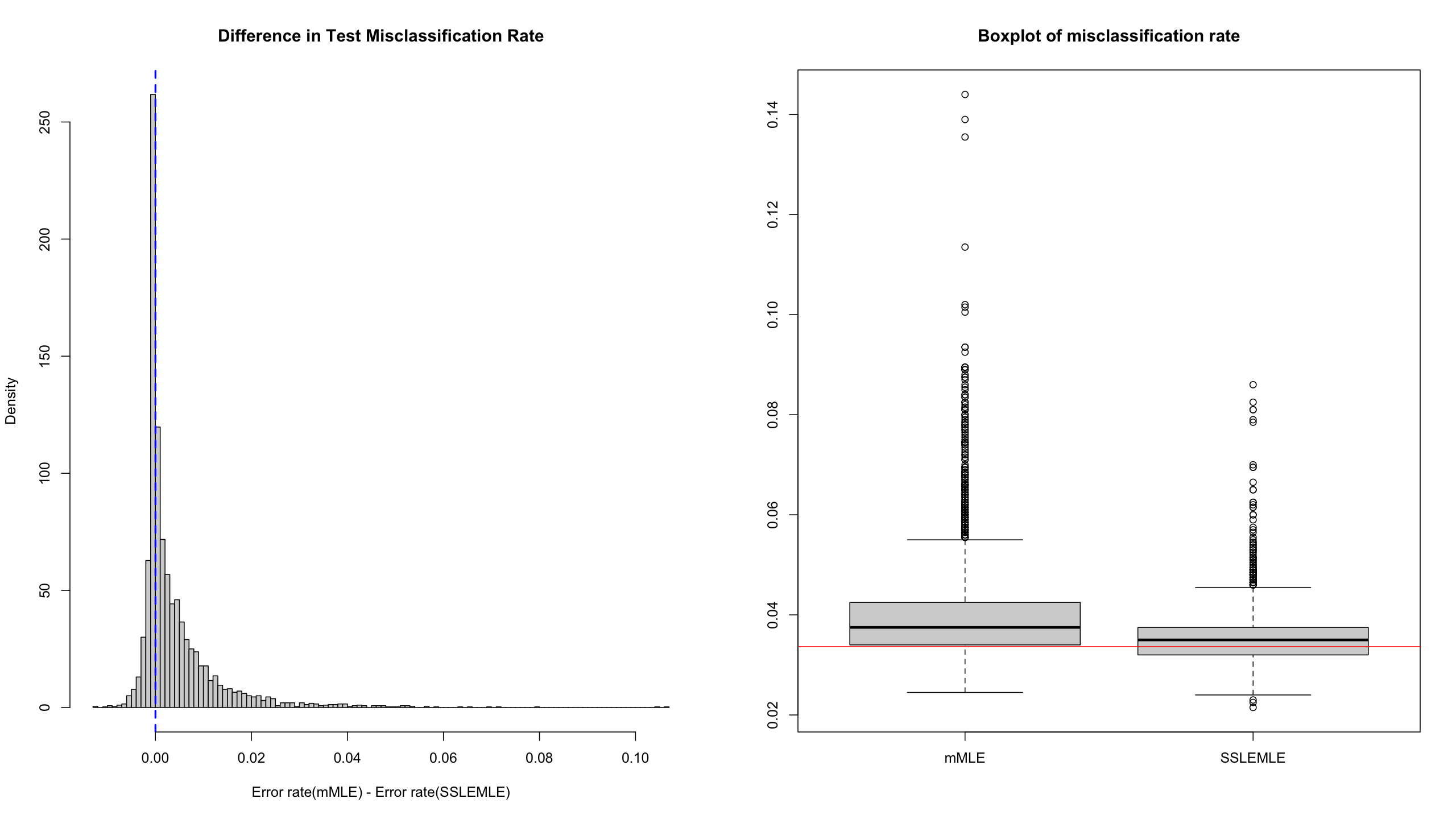}
    \caption{Top: Histogram of the MSE differences over 4,000 replications and boxplots of the MSEs for both estimators on the test data. Bottom: Histogram of the differences in misclassification rates and boxplots of the misclassification rates for both estimators, with the reference misclassification rate shown in red.}
    \label{plot: real data Misclassification and MSEs}
\end{figure}

We then proceeded as though only a small number of matched data is available, together with a large unmatched sample. In each replication, we first sampled 2,000 observations as our test data. From the remaining observations, we sampled $m=50$ matched observations and drew $n =10,000$ unmatched observations separately. We computed the matched MLE $\widehat{\beta}_m$ and the SSLEMLE $\widehat{\beta}_{n,m}$. This procedure was repeated 4,000 times. For each replication, we computed the mean squared error (MSE) and  misclassification rate on the test data. As shown in Figure \ref{plot: real data Misclassification and MSEs}, the SSLEMLE generally achieves lower MSEs and misclassification rates than the matched MLE.

\begin{figure}[htbp]
    \centering
    \includegraphics[width=0.8\textwidth]{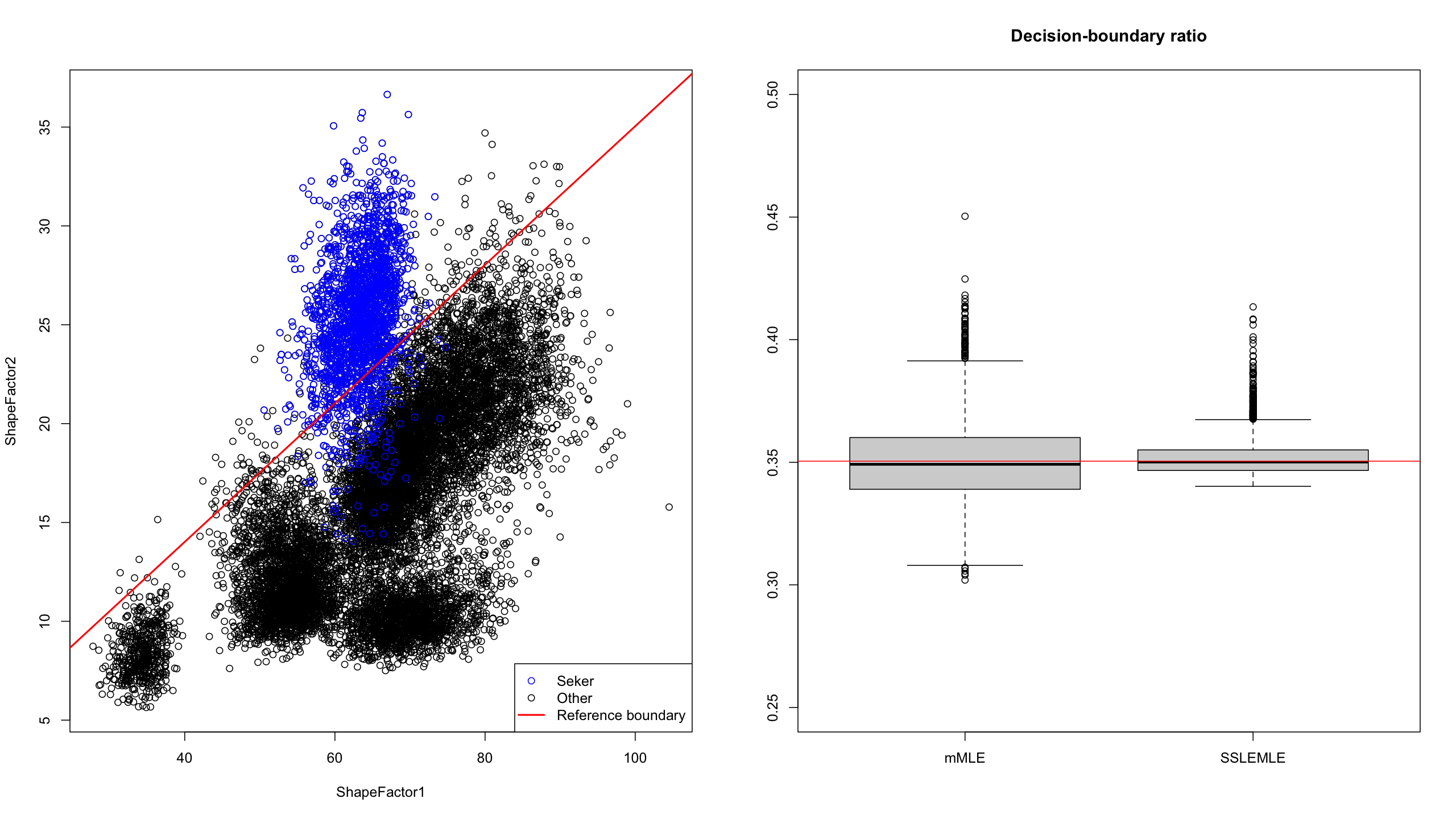}
    \caption{Left: All data points with the estimated reference value of $-\hat{\beta}_{1}/\hat{\beta}_{2}$ shown in red. Right: The ratios $-\widehat{\beta}_{m,1}/\widehat{\beta}_{m,2}$ and $-\widehat{\beta}_{n,m,1}/\widehat{\beta}_{n,m,2}$ over 4,000 replications. The mMLE has two extreme values that lie outside the shown range. The red line represents the ratio obtained from the reference estimator.}
    \label{plot: real data decision boundary}
\end{figure}

Note that in this real data example, if we classify a dry bean conditionally on the observed vector of features $X = x$ using 
$$
\mathbbm{1}_{\{H(\hat{\beta}_{1} x_1 + \hat{\beta}_{2}x_2) \geq 0.5\}}
$$
then, monotonicity of $H$  and the fact that $H(0)=0.5$ imply that this classification rule is equivalent to
\begin{align*}
    \mathbbm{1}_{\{\hat{\beta}_{1} x_1 + \hat{\beta}_{2}x_2 \geq 0\}}.
\end{align*}
If the signs of $\hat{\beta}_{0,1}$ and $\hat{\beta}_{0,2}$ are correctly estimated as in Table \ref{Table: real data reference model}, with
$\hat{\beta}_{1}<0$ and $\hat{\beta}_{2}>0$, then the classification rule can be re-written as 
\begin{align*}
     \mathbbm{1}_{\{ \frac{x_2}{x_1}\geq -\frac{\hat{\beta}_{1}}{\hat{\beta}_{2}} \}}.
\end{align*}
This means that, for the threshold $0.5$, only the ratio $-\hat{\beta}_{1}/\hat{\beta}_{2}$ determines the classification decision. As shown in Figure \ref{plot: real data decision boundary}, the SSLEMLEs provide a more concentrated estimation for the value $-\hat{\beta}_{1}/\hat{\beta}_{2}$ around the reference.

\section{Conclusions, discussion and future work}

This paper extends the SSLEMLE results obtained in \cite{balabdaoui2026SSLEMLE} in the framework of  linear regression to binary regression. Under some regularity assumptions, we were able to prove existence, consistency and asymptotic normality of the SSLEMLE. The asymptotic results were derived in the very general scheme of two different link functions and different distributions of the covariates in the matched and unmatched samples. The established theory allowed us further to quantify the asymptotic statistical gain of the SSLEMLE over the matched MLE and investigate this gain through synthetic and real datasets in several settings.

In Section \ref{Section: Matched and Unmatched Data with a Common Link Function and Covariate Distribution}, we assumed that $H=\tH$ is known. A natural extension would be to consider the case where  $H$ is only partially known or even unknown. For simplicity, we assume next that $\mathbb{P}^X = \mathbb{P}^{\tX}$. 

\paragraph{The case of  \texorpdfstring{$H$}{H} known only up to a scale}
Let us  assume that $H$ belongs to a scale family. This means that we now introduce an unknown scale parameter $\sigma_0 > 0$ such  $H =  H_1(\cdot/ \sigma_0)$, where $H_1$ is a known link function with scale parameter known to be equal to $1$. Then, the true model for a matched pair $(X, Y) $ is given by 
\begin{align*}
    \mathbb{P}\left(Y=\delta \mid X\right)  & = H \left(\beta_0^\top X\right)^\delta \,  \left(1- H \left(\beta_0^\top X\right)\right)^{1-\delta} \\
    &=  H_1\left(\frac{\beta_0^\top X}{\sigma_0}\right)^\delta \, \left( 1- H_1\left(\frac{\beta_0^\top X}{\sigma_0}\right) \right)^{1-\delta} . 
\end{align*}
for $\delta \in \{0, 1\}$. In the absence of pairing, the marginal distribution of $Y$ is 
\begin{eqnarray*}
\mathbb{P}(Y=\delta) =   \left \{\int H_1\left(\frac{\beta_0^\top x}{\sigma_0}\right) dF^X(x)\right\}^\delta \, \left \{1- \int H_1\left(\frac{\beta_0^\top x}{\sigma_0}\right) dF^X(x)\right \}^{1-\delta}.
\end{eqnarray*}
Thus, it is clear that the maximum likelihood approach studied above cannot be re-used in this case since one can replace simply $(\beta_0, \sigma_0)$ by another other pair $(\beta'_0, \sigma'_0)$ such that 
\begin{eqnarray}\label{noniden}
\frac{1}{\sigma_0} \beta_0  = \frac{1}{\sigma'_0} \beta'_0. 
\end{eqnarray}
Therefore, $\beta_0$ and $\sigma_0$ cannot be separately identified without an additional normalizing constraint. To make $\beta_0$ identifiable, we can for example impose the constraint that $\lVert\beta_0\rVert = 1$, whenever $\beta_0 \ne 0$. In this case, and by taking the Euclidean norm of both sides of the equality in (\ref{noniden}), we get that $\sigma'_0 = \sigma_0$. In other words, under the additional constraint that $\beta_0 \in \mathcal S^{p-1}$, the model becomes identifiable.  
Let $(\sigma, \beta) \mapsto \ell_{n,m}(\beta, \sigma)$ denote the natural extension of the log-likelihood $\beta \mapsto \ell_n(\beta)$ in (\ref{loglik}) to the more general case with 
\begin{eqnarray*}
(\beta, \sigma) \in \mathcal S \equiv \{\{0\} \times (0, \infty)\} \cup  \{\mathcal S^{p-1} \times (0, \infty) \}
\end{eqnarray*}
Then, it is not difficult to see that
\begin{eqnarray*}
\sup_{(\beta, \sigma) \in \mathcal S} \ell_{n,m}(\beta, \sigma)  = \sup_{v \in  \mathbb{R}^p} \ell_{n, m} (v). 
\end{eqnarray*}
Hence, we can mathematically return to the original optimization problem. In fact, if $\widehat v_n$ is the SSLEMLE with $H_0$ as the true link function, then
\begin{itemize}
\item if $\lVert \widehat v_n \rVert \ne 0$, $\widehat \beta_n  =  \widehat v_n / \Vert \widehat v_n \Vert$ and $\widehat \sigma_n =  1/\Vert \widehat v_n \Vert$,

\item if $\widehat v_n = 0$, then $\widehat \beta_n  = 0$ and $\widehat \sigma_n $ can be any real number in $(0, \infty)$.
\end{itemize}

Note that relaxing the assumption that $H$ is known might be very important in complex applications where it might be unclear which link function should be used.  This motivates us to  consider a more general semi-parametric model in which $H$ is only assumed to satisfy some qualitative constraint.


\paragraph{Maximum likelihood with  \texorpdfstring{$H$}{H} log-concave}

Let $\mathcal H_{\text{lc}}$ be the class of all distribution functions that are log-concave on $\mathbb R$. A distribution $H \in \mathcal H_{\text{lc}}$ if and only if there exists a function $\psi$ on $\mathbb R$ such that
\begin{itemize}
\item $\psi$ is a proper concave function on $\mathbb R$;
\item $\psi$ is non-decreasing;
\item $\lim_{t \to -\infty} \psi(t) = -\infty$ and $\lim_{t \to \infty} \psi(t) = 0$;
\item $H = e^{\psi}$.
\end{itemize}
In a future work, one can assume first that the true link function $H$ is only known to be an element of $ \mathcal H_{\text{lc}}$. Another related shape constraint is bi-log-concavity, that is, both $H$ and $1- H$ are log-concave. Such a constraint is satisfied by any CDF which admits a log-concave density.  
 Under either constraint,  the problem of maximizing the matched and pooled log-likelihoods, $\ell_m$ and $\ell_{n,m}$,  becomes now semi-parametric. To the best of our knowledge, such semi-parametric estimation problems have not been considered before in the literature. However, there are related works that are important to mention. The complexity of the problem considered here is reminiscent of that of estimation in the convex single index model studied by \cite{KuchibhotlaPatraSen2023}.  In that paper, the authors consider the model 
\begin{eqnarray*}
Y =m_0(\theta_0^{\top}X)+\varepsilon
\end{eqnarray*}
where $\mathbb{E}(\varepsilon\mid X)=0$, the index parameter $\theta_0\in\mathbb{R}^p$ is unknown, and the unknown regression function $m_0$ is convex. They propose a convex and Lipschitz constrained least-squares estimator for jointly estimating $m_0$ and $\theta_0$. Under mild moment and design assumptions, they establish consistency and convergence rates, showing in particular that the estimated regression function attains the minimax-optimal $n^{-2/5}$ rate in $L_2$ norm. Under additional smoothness, convexity, and identifiability conditions, they prove that the estimator of $\theta_0$ is $\sqrt{n}$-consistent and asymptotically normal. One of the main contributions of that work is to show that when the errors are homoscedastic, this estimator is semi-parametrically efficient. In the proof, the author make an important methodological contribution by presenting a geometric proof of efficiency based on projections onto the convex function class, which is designed to handle the non-smooth, shape-constrained nuisance estimator of $m_0$.  The other work related to semi-parametric estimation in binary regression under shape constraints is \cite{GroeneboomHendrickx2018} on the current status linear regression model.  The authors study the semi-parametric linear regression model
$$
Y = \beta_0^{\top}X + \varepsilon,
$$
when the response variable $Y$ is not observed directly. Instead, the available data consist of
$(T,X,\Delta)$, where $\Delta = \mathds{1}_{\{Y \leq T\}},
$ and $T$ is an \lq\lq inspection\rq\rq \ time, and the error $\varepsilon$ is independent of $(T, X)$ and has an unknown distribution function $F_0$. Thus, one observes only whether the response has occurred before the inspection time. The authors construct score based estimators of $\beta_0$ using the nonparametric maximum likelihood estimator of $F_0$, which acts as an infinite dimensional nuisance parameter. Their main result is that suitable estimators of $\beta_0$ are $\sqrt{n}$ consistent and asymptotically normal, even though the nuisance parameter is estimated non-parametrically and no smoothing of its maximum likelihood estimator is required. Under finiteness of the generalized Fisher information, a truncated version of the score yields estimators whose asymptotic variances can be made arbitrarily close to the semi-parametric efficiency bound.  Note that in estimating $F_0$, \cite{GroeneboomHendrickx2018} use isotonic regression since any distribution function is monotone. In a future work, one might add the smoother assumption of log-concavity of bi-log-concavity. If $T = 0$ with probability 1, then the current status linear regression model and our link-it regression model coincide. This naturally prompts the question on whether the asymptotic results obtained in this work can be easily extended to the setting where $T$ is degenerate, which clearly violates one of the assumptions made in \cite{GroeneboomHendrickx2018}. 

\bibliographystyle{plainnat}
\vskip 0.2in
\bibliography{unlinkedreg2}

\appendix
\section{Proof of existence} \label{Appendix: Existence}

\begin{proof}[Proof of Lemma \ref{lemma: Existence finite sample case}]
     Let $r_m$ be as above and denote also the number of successes in unmatched data by $r_n$. We have that 
     \begin{align*}
          \ell_{n,m} (0) &= \frac{r_n}{n+m} \log{\tH(0)} + \frac{n-r_n}{n+m} \log{\left(1-\tH(0)\right)} \\
          &+\frac{r_m}{n+m} \log{H(0)} + \frac{m-r_m}{n+m} \log{\left(1-H(0)\right)}.
     \end{align*}
     Define
     \begin{align}\label{M}
         M:= \min{\{\log{H(0)}, \, \log{\left(1-H(0)\right)}, \, \log{\tH(0)} ,\,  \log{\left(1-\tH(0)\right)}    \}} < 0.
     \end{align}
     Then, $\ell_{n,m} (0) \geq M$.   \\
     
  \par \noindent Consider the case $\mathbf{0 < r_m < m}$. We assume that the condition 
  $$0\in\operatorname{int}(\operatorname{conv}(X_1,\cdots,X_{r_m}, -X_{r_m+1}, -X_{r_m+2},\cdots, -X_m))$$
  is satisfied. We can write any given $\beta \in \mathbb{R}^p$ as $\beta = \lVert\beta\rVert u_\beta$ with $u_\beta \in\mathcal{S}^{p-1}$. Define \begin{align*}
        r_\beta^- :=  \sum_{k=1}^m \mathbbm{1}_{\{u_\beta^\top X_k<0\}}  \mathbbm{1}_{\{Y_k=1\}} =  \sum_{k=1}^{r_m} \mathbbm{1}_{\{u_\beta^\top X_k<0\}} 
    \end{align*}
    and 
    \begin{align*}
        r_\beta^+ :=  \sum_{k=1}^m \mathbbm{1}_{\{u_\beta^\top X_k>0\}}  \mathbbm{1}_{\{Y_k=0\}} =  \sum_{k=r_m+1}^{m} \mathbbm{1}_{\{u_\beta^\top X_k>0\}} .
    \end{align*}
Recall that we re-order the data so that the first $r_m$ observed matched pairs correspond to successes; i,e., $Y_k = 1, k =1, \dots, r_m$. This also means that $Y_k = 0$ for $k = r_m +1, \ldots, m $. Then,
\begin{align*}
&\sum_{k=1}^{m} \left( Y_k \log{H( \lVert\beta\rVert u_\beta^\top  X_k)} + \left(1-Y_k\right) \log{\left(1-H \left( \lVert\beta\rVert u_\beta^\top  X_k\right)\right)} \right) \\
& = \sum_{k=1}^{r_m}  \log{H( \lVert\beta\rVert u_\beta^\top  X_k)} +  \sum_{k = r_m+1}^m  \log{\left(1-H \left( \lVert\beta\rVert u_\beta^\top  
X_k\right)\right)} \\
& \le \sum_{k=1}^{r_m}  \log{H( \lVert\beta\rVert u_\beta^\top  X_k)} \mathbbm{1}_{\{u_\beta^\top X_k < 0\}}+  \sum_{k = r_m+1}^{m}  \log{\left(1-H \left( \lVert\beta\rVert u_\beta^\top  
X_k\right)\right)} \mathbbm{1}_{\{u_\beta^\top X_k>0\}}
\end{align*}
using the fact that $ \log(H) \le 0$ and $\log(1- H) \le 0$ on $\mathbb R$.
By Lemma \ref{lemma: bi log concavity left derivative} we have  
\begin{align*}
    \log{H}(t) \leq \frac{H_-'(0)}{H(0)} t, \quad \text{ and } \quad  \log{(1-H(t))}  \leq -\frac{H_-'(0)}{1-H(0)} t.
\end{align*}
It follows that
\begin{align*}
     &\sum_{k=1}^{r_m} \log{H(\lVert\beta\rVert u_\beta^\top X_k)} \mathbbm{1}_{\{u_\beta^\top X_k<0\}} + \sum_{k=r_m+1}^{m}  \log{\left(1-H( \lVert\beta\rVert u_\beta^\top X_k)\right)} \mathbbm{1}_{\{u_\beta^\top X_k>0\}}  \\
     &\leq  \frac{H_-'(0)}{H(0)} \sum_{k=1}^{r_m} \lVert\beta\rVert u_\beta^\top X_k  \mathbbm{1}_{\{u_\beta^\top X_k<0\}} -\frac{H_-'(0)}{1-H(0)}  \sum_{k=r_m+1}^{m}   \lVert\beta\rVert u_\beta^\top X_k\mathbbm{1}_{\{u_\beta^\top X_k>0\}} \\
     &= \lVert \beta \rVert \left( \frac{H_-'(0)}{H(0)}  \sum_{k=1}^{r_m} u_\beta^\top X_k  \mathbbm{1}_{\{u_\beta^\top X_k<0\}} + \frac{H_-'(0)}{1-H(0)}  \sum_{k=r_m+1}^{m}  - u_\beta^\top X_k\mathbbm{1}_{\{u_\beta^\top X_k>0\}} \right) \\
     &\leq \lVert \beta \rVert \frac{H_-'(0)}{L} \left(  \sum_{k=1}^{r_m} u_\beta^\top X_k  \mathbbm{1}_{\{u_\beta^\top X_k<0\}} +  \sum_{k=r_m+1}^{m}  - u_\beta^\top X_k\mathbbm{1}_{\{u_\beta^\top X_k>0\}} \right) \\
     &\leq \lVert \beta \rVert  \frac{H_-'(0)}{L}  \cdot \sup_{v\in \mathcal{S}^{p-1}} N_m(v) 
\end{align*}
with
\begin{align*}
    N_m(v) & =    \sum_{k=1}^{r_m} v^\top X_k  \mathbbm{1}_{\{v^\top X_k<0\}} + \sum_{k=r_m+1}^{m}  -v^\top X_k\mathbbm{1}_{\{v^\top X_k>0\}} \\
    &=   \sum_{k=1}^{r_m} \min{\{v^\top X_k,0\}} + \sum_{k=r_m+1}^{m}  -\max{\{v^\top X_k,0\}},  \ \ \text{and} \\
    L & := \max{\{H(0),1-H(0)\}}.
\end{align*}
Since $v \mapsto N_m(v)$ is continuous and $\mathcal{S}^{p-1}$ is compact, there exists some $v^\ast \in  \mathcal{S}^{p-1}$ such that $\sup_{v\in \mathcal{S}^{p-1}} N_m(v) = N_m(v^\ast)$. Hence,
\begin{align*}
&\sum_{k=1}^{r_m} \log{H(\lVert\beta\rVert u_\beta^\top X_k)} \mathbbm{1}_{\{u_\beta^\top X_k<0\}} + \sum_{k=r_m+1}^{m}  \log{\left(1-H( \lVert\beta\rVert u_\beta^\top X_k)\right)} \mathbbm{1}_{\{u_\beta^\top X_k>0\}} \\
& \leq \frac{H_-'(0)}{L}   N_m(v^\ast)   \lVert \beta \rVert.
\end{align*}
It is easy to see that $N_m(v) \leq 0$ for all $v \in \mathcal{S}^{p-1}$. We now prove that the inequality is strict; i.e., $N_m(v) < 0$ for all $v \in \mathcal{S}^{p-1}$. We show this by contradiction.  Suppose that there exists $w \in \mathcal{S}^{p-1}$ with $N_m(w) =0$. This can only occur if $w^\top X_k \geq 0$ for all $k = 1,\dots,r_m$ and $w^\top X_k \leq 0$ for all $k = r_m+1,\dots,m$. In other words, if $N_m(w) =0$ then  
\begin{align*}
    \sum_{k=1}^{r_m} \mathbbm{1}_{\{w^\top X_k<0\}} + \sum_{k=r_m+1}^{m} \mathbbm{1}_{\{w^\top X_k>0\}} =r^-_w  +  r^+_w = 0,
\end{align*}
which is not true. In fact, since $0\in\operatorname{int}(\operatorname{conv}(X_1,\cdots,X_{r_m}, -X_{r_m+1}, \cdots, -X_m))$, then Lemma \ref{lemma: existence-auxiliary} applied to 
$$  d=m, \ \text{and} \   (z_1, \ldots, z_d) = (X_1,\cdots,X_{r_m}, -X_{r_m+1},\cdots, -X_m)$$
implies that 
    \begin{align*}
       & \inf_{v\in \mathcal{S}^{p-1}}\frac{1}{m} \left(\sum_{k=1}^{r_m} \mathbbm{1}_{\{v^\top X_k<0\}}   +  \sum_{k=r_m+1}^m \mathbbm{1}_{\{-v^\top X_k < 0\}} \right)  \\
        &= \inf_{v\in \mathcal{S}^{p-1}}\left(\frac{1}{m} \sum_{k=1}^{r_m} \mathbbm{1}_{\{v^\top X_k<0\}}   + \frac{1}{m} \sum_{k=r_m+1}^m \mathbbm{1}_{\{v^\top X_k > 0\}} \right)
        > 0.
    \end{align*}
    Therefore, for all $\beta \in \mathbb R^p$
    \begin{eqnarray*}
     r^-_\beta  +  r^+_\beta  \ge m  \inf_{v\in \mathcal{S}^{p-1}}\left(\frac{1}{m} \sum_{k=1}^{r_m} \mathbbm{1}_{\{v^\top X_k<0\}}   + \frac{1}{m} \sum_{k=r_m+1}^m \mathbbm{1}_{\{v^\top X_k > 0\}} \right)  > 0,
    \end{eqnarray*}
and hence $r^-_w  +  r^+_w $ cannot be equal to $0$.  Using the fact that $\ell_n(\beta) \le 0$ for all $\beta \in \mathbb R^p$, it follows that 
\begin{eqnarray*}
\ell_{n, m}(\beta)  \le \frac{m}{n+m} \ell_m(\beta).
\end{eqnarray*}
The calculations above imply that the term on the right side is strictly smaller than $M$ provided that 
$$
\frac{m}{n+m} \Vert \beta \Vert \frac{H_-'(0)}{L} \sup_{v\in \mathcal{S}^{p-1}}\frac{1}{m} N_m(v)  < M
$$
where $M$ is the same constant in (\ref{M}). The preceding condition is equivalently to 
$$\beta \in \{v: \lVert v\rVert > R_{n,m} \}$$ 
where
\begin{align*}
    R_{n,m} = \frac{ML}{H_-'(0)} \cdot\frac{m+n}{m} \frac{1}{\sup_{v\in \mathcal{S}^{p-1}}\frac{1}{m} N_m(v)}.
\end{align*}
This implies that  $\ell_{n,m}(\beta) < \ell_{n,m}(0)$ for all $\beta \in \mathbb R^p$.

\bigskip

\par \noindent Consider now the case $\mathbf{r_m \in \{0,m\}}$. Without loss of generality, we can swap the definitions of success and failure, so it suffices to consider the case $r_m = m$. Similarly, as in the previous case, define 
\begin{align*}
        r_\beta^- := \sum_{k=1}^m \mathbbm{1}_{\{u_\beta^\top X_k<0\}}.
\end{align*}
Then, we have that 
\begin{align*}
    \ell_{n,m}(\beta) &\leq \frac{1}{n+m}\sum_{k=1}^m  \log H(\beta^\top X_k)  \\
&\leq  \frac{1}{n+m} \sum_{k=1}^{m} \log{ H(\lVert\beta\rVert u_\beta^\top X_k)}  \mathbbm{1}_{\{u_\beta^\top X_k<0\}}, \ \text{since $\log H \le 0$}  \\
    &\leq   \lVert\beta\rVert \frac{H_-'(0)}{H(0)} \frac{m}{n+m}\frac{1}{m}\sum_{k=1}^{m}u_\beta^\top X_k  \mathbbm{1}_{\{u_\beta^\top X_k<0\}}, \ \text{by (A1')} \\
    & = \lVert\beta\rVert \frac{H_-'(0)}{H(0)} \frac{m}{n+m}\frac{1}{m}\sum_{k=1}^{m} \min{\{u_\beta^\top X_k, 0\}} \\
    & \leq  \lVert\beta\rVert \frac{H_-'(0)}{H(0)} \frac{m}{n+m} \sup_{v\in \mathcal{S}^{p-1}} \frac{1}{m} \sum_{k=1}^{m} \min{\{v^\top X_k,0\}}.
\end{align*}
Since $0\in\operatorname{int}(\operatorname{conv}(X_1,\cdots,X_m))$, we have $r_\beta^- > 0$ for all $\beta \in \mathbb{R}^p$. Using continuity of the function $v \mapsto \frac{1}{m} \sum_{k=1}^{m} \min{\{v^\top X_k,0\}}$  and compactness of $\mathcal{S}^{p-1}$, it follows from arguments similar to those in the previous case that $\ell_{n,m}(\beta) < \ell_{n,m}(0)$ for all $\beta \in \{v: \lVert v \rVert > R_{n,m}\}$ with
\begin{align*}
    R_{n,m} = \frac{H(0)M}{H_-'(0)}\cdot \frac{m+n}{n}\frac{1}{\sup_{v\in \mathcal{S}^{p-1}} \frac{1}{m} \sum_{k=1}^{m} \min{\{v^\top X_k,0\}} }.
\end{align*}
Now, consider the closed ball  $\overline{\mathcal{B}}(0,R_{n,m}) := \{v\in\mathbb{R}^p: \left\lVert v\right\rVert \leq R_{n,m}\}$. Since $\beta \mapsto \ell_{n,m}(\beta)$ is continuous, this function attains its supremum on $\overline{\mathcal{B}}(0,R_{n,m})$ at some vector $\beta^\ast \in \overline{\mathcal{B}}(0,R_{n,m})$. Thus,  
\[
\ell_{n,m}(\beta^\ast) \geq  \ell_{n,m}(\beta) \quad \text{for all $\beta \in \overline{\mathcal{B}}(0,R_{n,m})$},
\]
which holds in particular for  $\beta = 0 \in \overline{\mathcal{B}}(0,R_{n,m})$. Hence,
\[
\ell_{n,m}(\beta^\ast) \geq  \ell_{n,m}(0).
\]
Since for all $\beta \in \{v \in\mathbb{R}^p: \left\lVert v \right\rVert >R_{n,m}\}$
\[
\ell_{n,m}(\beta^\ast) \geq  \ell_{n,m}(0) > \ell_{n,m}(\beta)
\]
it follows that $\beta^\ast$ is a maximizer of $\beta \mapsto \ell_{n,m}(\beta)$ over $\mathbb{R}^p$.

\end{proof}

\section{Proof of consistency} \label{Appendix: Consistency}

\begin{proof}[Proof of Proposition \ref{unifcons}.]
We first re-write $\ell$ as 
\begin{eqnarray*}
\ell & = &  \left(\frac{1}{\lambda +1} - \frac{n}{n + m} \right) \ell_U  + \left(\frac{\lambda}{\lambda +1} - \frac{m}{n+ m} \right) \ell_M  +  \frac{n}{n + m} \ell_U  + \frac{m}{n+ m}  \ell_M.
\end{eqnarray*}
We have that 
\begin{eqnarray*}
\left \vert \left(\frac{1}{\lambda +1} - \frac{n}{n + m} \right) \ell_U(\beta)  + \left(\frac{\lambda}{\lambda +1} - \frac{m}{n+ m} \right) \ell_M(\beta) \right \vert  \le \frac{\vert m/n - \lambda \vert}{(1+ \lambda)(1+ m/n)} \left(\vert \ell_U(\beta) \vert  +     \vert \ell_M(\beta) \vert \right).
\end{eqnarray*}
We show now that $\sup_{\beta \in \overline{\mathcal{B}}(0, R)} \vert \ell_U(\beta) \vert \vee \sup_{\beta \in \overline{\mathcal{B}}(0, R)} \vert \ell_U(\beta) \vert < \infty$. We have that 
\begin{eqnarray*}
\vert \ell_U(\beta) \vert & \le  &   \left \vert \log \left( \int \tH(\beta^\top x) d F^{\tX}(x) \right)  \right \vert +  \left \vert \log \left( 1- \int \tH(\beta^\top x) d F^{\tX}(x) \right)  \right \vert  \\
& \le &  \vert \log (\widetilde{c}_R) \vert  +  \vert \log (1-\widetilde{C}_R) \vert 
\end{eqnarray*}
for all $\beta \in \overline{\mathcal B}(0, R)$ by Assumption (A5-2).  Also,  we have that 
\begin{eqnarray*}
\vert \ell_M(\beta) \vert & \le  &  \int \vert \log H(\beta^\top x) \vert dF^X(x)  +  \int \vert \log (1- H(\beta^\top x)) \vert dF^X(x)  \\
    & \leq &  \left( \int \vert \log H(\beta^\top x) \vert^2 dF^X(x) \right)^{\frac{1}{2}}+ \left(\int \vert \log (1- H(\beta^\top x)) \vert^2 dF^X(x) \right)^{\frac{1}{2}}\\ 
    & \leq & \left( \int ( \log H(-R\lVert x\rVert) )^2 dF^X(x) \right)^{\frac{1}{2}}+ \left(\int ( \log (1- H(R \lVert x\rVert)) )^2 dF^X(x) \right)^{\frac{1}{2}} < \infty 
\end{eqnarray*}
using Assumptions (A1) and (A3).  We conclude that there exists some constant $\tilde{D}_R > 0$ such that for all $\beta \in \overline{\mathcal B}(0, R)$
\begin{eqnarray*}
\left \vert \left(\frac{1}{\lambda +1} - \frac{n}{n + m} \right) \ell_U(\beta)  + \left(\frac{\lambda}{\lambda +1} - \frac{m}{n+ m} \right) \ell_M(\beta) \right \vert  & \le  &  \tilde{D}_R \frac{\vert m/n - \lambda \vert}{(1+ \lambda)(1+ m/n)} \\
& \le & \frac{\tilde{D}_R}{1+ \lambda} \left \vert \frac{m}{n} - \lambda \right \vert  \to 0
\end{eqnarray*}
as  $m, n \to \infty $.   On the other hand, we have that 
\begin{eqnarray*}
&&\ell_{n, m}(\beta)  -  \frac{n}{n+m} \ell_U(\beta)  -  \frac{m}{n+m} \ell_M(\beta)\\
&&= \frac{n}{n+m} \bigg \{\overline{\tY}_n \log\left(  \int \tH(\beta^\top x) d\mathbb F^{\tX}_n(x)\right) - \mathbb E(\tY) \log\left( \int \tH(\beta^\top x) dF^{\tX}(x)  \right) \\
&& \hspace{1.2cm} + \ (1-\overline{\tY}_n) \log\left(1- \int \tH(\beta^\top x) d\mathbb F^{\tX}_n(x)\right)  - (1-\mathbb E(\tY)) \log\left( 1- \int \tH(\beta^\top x) dF^{\tX}(x) \right)\bigg \} \\
&& \ + \frac{m}{n + m} \int \left\{ y \log(H(\beta^\top x)) + (1-y) \log (1- H(\beta^\top x)) \right \} d(\mathbb P_m- \mathbb P)(x, y) 
\end{eqnarray*}
where 
\begin{eqnarray*}
&& \overline{\tY}_n \log\left(  \int \tH(\beta^\top x) d\mathbb F^{\tX}_n(x)\right) - \mathbb E(\tY) \log\left( \int \tH(\beta^\top x) dF^{\tX}(x)  \right)  \\
&&=   \overline{\tY}_n \log\left( \frac{\int \tH(\beta^\top x) d\mathbb F^{\tX}_n(x)}{\int \tH(\beta^\top x) d F^{\tX}(x)} \right)  + (\overline{\tY}_n- \mathbb E(\tY))  \log\left( \int \tH(\beta^\top x) dF^{\tX}(x)  \right) \\
&&:= A_n(\beta)
\end{eqnarray*}
and 
\begin{eqnarray*}
 && (1-\overline{\tY}_n) \log\left(1- \int \tH(\beta^\top x) d\mathbb F^{\tX}_n(x)\right)  - (1-\mathbb E(\tY)) \log\left( 1- \int \tH(\beta^\top x) dF^{\tX}(x) \right) \\
 && = (1-\overline{\tY}_n) \log\left( \frac{\int \left(1- \tH(\beta^\top x)\right) d\mathbb F^{\tX}_n(x)}{\int \left(1- \tH(\beta^\top x) \right) d F^{\tX}(x) } \right)  + (\mathbb E(\tY) - \overline{\tY}_n)  \log\left( \int \left(1- \tH(\beta^\top x) \right) dF^{\tX}(x)  \right) \\
 &&: = B_n(\beta).
\end{eqnarray*}
Denote also
\begin{eqnarray*}
C_m(\beta): =  \int \left\{ y \log(H(\beta^\top x)) + (1-y) \log (1- H(\beta^\top x)) \right \} d(\mathbb P_m- \mathbb P)(x, y).
\end{eqnarray*}
We will show next that 
\begin{eqnarray*}
\sup_{\beta \in \overline{\mathcal{B}}(0,R)} \lvert A_n(\beta) \rvert =  O_{\mathbb P^{\tX} \otimes \mathbb P^{\tY}}(n^{-1/2}), \ \ \sup_{\beta \in \overline{\mathcal{B}}(0,R)} \lvert B_n(\beta) \rvert =  O_{\mathbb P^{\tX} \otimes \mathbb P^{\tY}}(n^{-1/2}),
\end{eqnarray*}
and 
\begin{eqnarray*}
\sup_{\beta \in \overline{\mathcal{B}}(0,R)} \lvert C_m(\beta) \rvert =  O_{\mathbb P}(m^{-1/2}).
\end{eqnarray*}
We start with the latter statement. Denote $h_+(t) := \log{H(t)}$. By Assumption (A1), $H: \mathbb{R} \to (0,1)$ is strictly monotone increasing for all $t \in \mathbb{R}$. Hence, since $\log{z}$ is also strictly monotone increasing for all $z \in (0,1)$, the function $h_+: \mathbb{R} \to (-\infty,0)$ is strictly monotone increasing. Similarly, if $h_{-}(t) := \log{(1-H(t))}$, then it is easy to see that $h_{-}: \mathbb{R} \to (-\infty,0)$ is strictly monotone decreasing. Consider now the class of functions 
\begin{align*}
    \mathcal{L} = \{(x,y)\mapsto l_{\beta}(x) = \beta^\top x, \, \beta\in\overline{\mathcal{B}}(0,R) \}.
\end{align*}
The class $\mathcal{L}$ is indexed by $\beta$ and hence is a subset of the finite dimensional vector space
\begin{eqnarray*}
\{(x,y)\mapsto  l_\beta(x) = \beta^\top x, \ \beta \in \mathbb R^p \}.   
\end{eqnarray*}
From \cite[Lemma 2.6.16]{aadbookE2} it follows that $\mathcal L$ is a VC-subgraph of dimension smaller than or equal to $ d+1$. Using \cite[Lemma 2.6.20 - (viii)]{aadbookE2}, it follows that the class $\mathcal{H}_+:=h_+ \circ \mathcal L$ is a VC-subgraph. Furthermore, if we define $g_+(x,y) = y$, then it follows by \cite[Lemma 2.6.20 - (vi)]{aadbookE2} that
\begin{align*}
    \mathcal{H}^{\sharp}_+ := \mathcal{H}_+ \cdot g_+ &= \{ (x,y) \mapsto f(x,y) g_+(x,y): f\in \mathcal{H}_+\} \\
    &= \{ (x,y) \mapsto y\cdot \log{H(\beta^\top x),  \, \beta\in\overline{\mathcal{B}}(0,R) }\}
\end{align*}
is VC-subgraph. Similarly, if $g_-(x,y) = 1-y$, then by similar arguments we can show that
\begin{align*}
    \mathcal{H}^{\sharp}_- := \mathcal{H}_- \cdot g_- &= \{ (x,y) \mapsto \left(\left( h_- \circ \mathcal{L}\right)g_- \right) (x,y)\} \\
    &= \{ (x,y) \mapsto (1-y)\cdot \log{(1-H(\beta^\top x)),  \, \beta\in\overline{\mathcal{B}}(0,R) }\}
\end{align*}
is also VC-subgraph. We have
\begin{align*}
    \sup_{\beta\in\overline{\mathcal{B}}(0,R)} \lvert C_m(\beta)\rvert \leq \sup_{f^{\sharp}_+ \in \mathcal{H}^{\sharp}_+} \left \vert \int f^{\sharp}_+ \ d(\mathbb P_m- \mathbb P)(x, y) \right \vert +\sup_{f^{\sharp}_- \in \mathcal{H}^{\sharp}_-} \left \vert \int f^{\sharp}_- \ d(\mathbb P_m- \mathbb P)(x, y) \right \vert.
\end{align*}
Since $y \in \{0,1\}$, then
\begin{align*}
    F(x) := -\log{H(-R\lVert x\rVert)} \geq \lvert y\log{H(\beta^\top x)} \rvert
\end{align*}
 is an envelop of $\mathcal{H}^{\sharp}_+$. It follows from \cite[Theorem 2.6.7]{aadbookE2} that for $\delta \in (0,1)$ and any probability measure $Q$ such that $\Vert F \Vert^2_Q = \int F(x)^2  dQ(x) > 0$
\begin{eqnarray}\label{proof: consis C covering}
N(\delta \Vert F \Vert_Q, \mathcal{H}^{\sharp}_+, L_2(Q)) \le D \left(\frac{1}{\delta}\right)^{2 V(\mathcal{H}^{\sharp}_+)} 
\end{eqnarray}
for some universal constant $D > 0$, which can be taken to be larger than $ 1/e$ without loss of generality. Define now the uniform entropy integral 
\begin{eqnarray*}
J(\eta, \mathcal{H}^{\sharp}_+)  = \sup_{Q}   \int_0^\eta  \sqrt{1 + \log N(\delta \Vert F \Vert_Q, \mathcal{H}^{\sharp}_+, L_2(Q)) } d\delta
\end{eqnarray*}
for $\eta > 0$. In the above definition, the supremum is taken again over $Q$ such that $\Vert F \Vert_Q > 0$. Using the bound given in (\ref{proof: consis C covering}), it follows that
\begin{eqnarray*}
J(1, \mathcal{H}^{\sharp}_+) & \le  & \int_0^1  \left(\sqrt{1 + \log D} +  \sqrt{2 V(\mathcal{H}^{\sharp}_+)} \sqrt{\log\left(\frac{1}{\delta} \right)} \right) d\delta  \\
& \le &  \sqrt{1 + \log D}  +  \sqrt{2 V(\mathcal{H}^{\sharp}_+)} \int_0^1 \frac{1}{\sqrt{\delta}} d \delta \\
& =  & \sqrt{1 + \log D}  +  2\sqrt{2 V(\mathcal{H}^{\sharp}_+)}.
\end{eqnarray*}
By Assumption (A3), 
$$\Vert F \Vert_{\mathbb P^X}: = \left(\int F(x)^2 d\mathbb P^X(x)\right)^{1/2}$$
is finite. By \cite[Theorem 2.14.1]{aadbookE2}, it follows that
\begin{eqnarray*}
\mathbb E[\Vert \mathbb G_m \Vert^2_{\mathcal{H}^{\sharp}_+}]^{1/2} \lesssim J(1, \mathcal{H}^{\sharp}_+) \Vert F \Vert_{\mathbb P^X} 
\end{eqnarray*}
where $\Vert \mathbb G_m \Vert_{\mathcal{H}^{\sharp}_+} =  \sup_{g \in \mathcal{H}^{\sharp}_+} \vert \mathbb G_m g \vert$.  Using the Markov's inequality, it follows that 
\begin{eqnarray*}
\Vert \mathbb G_m \Vert_{\mathcal{H}^{\sharp}_+} = O_{\mathbb P}(1)
\end{eqnarray*}
which in turn implies that 
\begin{eqnarray*}
\sup_{f^{\sharp}_+ \in \mathcal{H}^{\sharp}_+} \left \vert \int f^{\sharp}_+ \ d(\mathbb P_m- \mathbb P)(x, y) \right \vert  = O_{\mathbb P}(m^{-1/2}) = o_{\mathbb P}(1).
\end{eqnarray*}
Since a similar reasoning can be used for the class of functions $\mathcal{H}^{\sharp}_-$ we can show that 
\begin{eqnarray*}
\sup_{f^{\sharp}_- \in \mathcal{H}^{\sharp}_-} \left \vert \int f^{\sharp}_- \ d(\mathbb P_m- \mathbb P)(x, y) \right \vert  = O_{\mathbb P}(m^{-1/2}) = o_{\mathbb P}(1)
\end{eqnarray*}
to conclude that $\sup_{\beta \in \overline{\mathcal{B}}(0,R)} \lvert C_m(\beta) \rvert =  O_{\mathbb P}(m^{-1/2})$. 

\medskip

\par \noindent To upper bound $\vert A_n(\beta) \vert$, we start by noting that if $\int \tH(\beta^\top x) d\mathbb{F}_n^{\tX}(x) > \int \tH(\beta^\top x) dF^{\tX}(x)$, then 
\begin{eqnarray*}
0 < \log{\left( \frac{\int \tH(\beta^\top x) d\mathbb{F}_n^{\tX}(x)}{\int \tH(\beta^\top x) dF^{\tX}(x)}\right)} &= & \log{\left( 1 + \frac{\int \tH(\beta^\top x) \left( d\mathbb{F}_n^{\tX}(x) - dF^{\tX}(x)\right)}{\int \tH(\beta^\top x) dF^{\tX}(x)}\right)} \\
& \leq  &  \frac{\int \tH(\beta^\top x) d(\mathbb{P}_n^{\tX} - \mathbb{P}^{\tX})(x)}{\int \tH(\beta^\top x) dF^{\tX}(x)} \\
& = &   \frac{1}{\int \tH(\beta^\top z) dF^{\tX}(z)}  \frac{1}{\sqrt n} \int \tH(\beta^\top x) \ d\mathbb G^{\tX}_n(x).
\end{eqnarray*}
If $\int \tH(\beta^\top x) d\mathbb{F}_n^{\tX}(x) < \int \tH(\beta^\top x) d\mathbb{F}^{\tX}(x)$, then  
\begin{align*}
    0& > \log{\left( \frac{\int \tH(\beta^\top x) d\mathbb{F}_n^{\tX}(x)}{\int \tH(\beta^\top x) d\mathbb{F}^{\tX}(x)}\right)}  =  -\log{\left( \frac{\int \tH(\beta^\top x) dF^{\tX}(x)}{\int \tH(\beta^\top x) d\mathbb{F}_n^{\tX}(x)}\right)} \geq \frac{\int \tH(\beta^\top x)\ d(\mathbb{P}_n^{\tX} - \mathbb{P}^{\tX})(x)}{\int \tH(\beta^\top x) d\mathbb{F}_n^{\tX}(x)} \\
    &=\frac{\int \tH(\beta^\top x) d\mathbb G^{\tX}_n(x)}{\int \tH(\beta^\top z) d\mathbb G^{\tX}_n(z) + \sqrt{n} \int \tH(\beta^\top z) dF^{\tX}(z)}.
\end{align*}
Combining this with Assumption (A5-2), we can write
\begin{align*}
    \lvert A_n(\beta)\rvert &\leq \left \vert \log{\left( \frac{\int \tH(\beta^\top x) d\mathbb{F}_n^{\tX}(x)}{\int \tH(\beta^\top x) d{F}^{\tX}(x)}\right)}\right \vert + \lvert \log{\widetilde{c}_R}\rvert \cdot \frac{1}{\sqrt{n}}\left \vert \int y\  d\mathbb{G}_n^{\tY}(y)  \right \vert \\
    &\leq  \frac{1}{\widetilde{c}_R} \frac{1}{\sqrt{n}}\left \vert \int 
    \tH(\beta^\top x)  \ d\mathbb G^{\tX}_n(x)\right \vert + \left \vert \frac{\frac{1}{\sqrt n} \int \tH(\beta^\top x) d\mathbb G^{\tX}_n(x)}{\frac{1}{\sqrt n} \int \tH(\beta^\top z) d\mathbb G^{\tX}_n(z) + \int \tH(\beta^\top z) d F^{\tX}(z)}\right \vert  \\
    & \ \  \ \ +  \lvert \log{\widetilde{c}_R}\rvert \cdot \frac{1}{\sqrt{n}}\left \vert \int y\  d\mathbb{G}_n^{\tY}(y)  \right \vert \\
    &\leq  \frac{1}{\widetilde{c}_R}\frac{1}{\sqrt{n}}  \left \vert \int 
    \tH(\beta^\top x) \ d\mathbb G^{\tX}_n(x)\right \vert + \left \vert \frac{\frac{1}{\sqrt n} \int \tH(\beta^\top x) d\mathbb G^{\tX}_n(x)}{\frac{1}{\sqrt n} \int \tH(\beta^\top z) d\mathbb G^{\tX}_n(z) + \int \tH(\beta^\top z) d F^{\tX}(z)}\right \vert \\
    & \ \  \ \ +   \lvert \log{\widetilde{c}_R}\rvert \cdot \frac{1}{\sqrt{n}}\left \vert \int y\  d\mathbb{G}_n^{\tY}(y)  \right \vert. 
\end{align*}
Since $ \tH\circ \mathcal{L}$ is a VC-subgraph and $F(x) := 1$ is its envelope, then  similar arguments as those used for $C_m$ yield 
\begin{eqnarray}\label{Knscaled}
\sup_{\beta \in \overline{\mathcal{B}}(0,R)}\frac{1}{\sqrt{n}}  \left \vert \int \tH(\beta^\top x) \ d\mathbb G^{\tX}_n(x)\right \vert   =  O_{\mathbb{P}^{\tX}}(n^{-1/2}) = o_{\mathbb{P}^{\tX}}(1).
\end{eqnarray}
Also, since $\tY \sim \text{Bernoulli}(p_0)$ has a finite variance, we have by the CLT that $  \int y\  d\mathbb{G}_n^{\tY}(y) =  O_{\mathbb{P}^{\tY}}(1)$  and hence
$$
\frac{1}{\sqrt{n}}\left \vert \int y\  d\mathbb{G}_n^{\tY}(y)  \right \vert  = O_{\mathbb{P}^{\tY}}(n^{-1/2}) = o_{\mathbb{P}^{\tY}}(1).
$$
Put
\begin{eqnarray*}
K_n:= \frac{1}{\sqrt n} \sup_{\beta \in \overline{\mathcal{B}}(0,R)} \left \vert \int \tH(\beta^\top x) \ d\mathbb G^{\tX}_n(x) \right \vert.
\end{eqnarray*}
Recall the constant $\widetilde{c}_R$ from Assumption (A5-2). By (\ref{Knscaled}), it follows that $K_n \le \widetilde{c}_R/2 $ occurs with probability tending to $1$ as $n \to \infty$.   This implies that 
\begin{eqnarray*}
&&\inf_{\beta \in \overline{\mathcal{B}}(0,R)}\left \vert \frac{1}{\sqrt n} \int \tH(\beta^\top z) d\mathbb G^{\tX}_n(z) + \int \tH(\beta^\top z) d F^{\tX}(z)  \right \vert \\
&& \ge \inf_{\beta \in \overline{\mathcal{B}}(0,R)} \left \{\int \tH(\beta^\top z) d F^{\tX}(z)  -   \frac{1}{\sqrt n} \left \vert \int \tH(\beta^\top z) d\mathbb G^{\tX}_n(z)  \right \vert   \right \}  \\
&& \ge \widetilde{c}_R  -  K_n  \\
&& \ge \widetilde{c}_R/2
\end{eqnarray*}
with probability tending to 1 as $n\to \infty$. Thus,
\begin{eqnarray*}
&& \sup_{\beta \in \overline{\mathcal{B}}(0,R)} \left \vert \frac{\frac{1}{\sqrt n} \int \tH(\beta^\top x) d\mathbb G^{\tX}_n(x)}{\frac{1}{\sqrt n} \int \tH(\beta^\top z) d\mathbb G^{\tX}_n(z) + \int \tH(\beta^\top z) d F^{\tX}(z)}\right \vert  \\
&& \le   \frac{K_n}{\inf_{\beta \in \overline{\mathcal{B}}(0,R)}\left \vert \frac{1}{\sqrt n} \int \tH(\beta^\top z) d\mathbb G^{\tX}_n(z) + \int \tH(\beta^\top z) d F^{\tX}(z)  \right \vert} \\
&& \le  \frac{2}{\widetilde{c}_R}  K_n 
\end{eqnarray*}
with probability tending to $1$. Then, we can conclude from (\ref{Knscaled}) that $\sup_{\beta \in \overline{\mathcal{B}}(0,R)} \lvert A_n(\beta) \rvert =  O_{\mathbb P^{\tX} \otimes \mathbb P^{\tY}}(n^{-1/2})= o_{\mathbb{P}^{\tX} \otimes \mathbb P^{\tY}}(1)$.  \\

To bound $\vert B_n(\beta) \vert$,  we write 
\begin{align*}
    \overline{\tH}(t) := 1 - \tH(t).
\end{align*}
Note that Assumption (A5-2) holds also for $\overline{\tH}$; i.e.
\begin{align*}
    0< 1- \widetilde{C}_R \leq \int \overline{\tH}(\beta^\top x) dF^{\tX}(x) \leq 1- \widetilde{c}_R < 1
\end{align*}
for all $\beta \in \overline{\mathcal{B}} (0,R)$. Since $\overline{\tH}$ is monotone decreasing, $\overline{\tH} \circ \mathcal{L}$ is a VC-subgraph. Since it admits $F(x):= 1$ as envelope, we can use similar arguments as for bounding $\vert A_n(\beta) \vert$ to show that $\sup_{\beta \in \overline{\mathcal{B}}(0,R)} \lvert B_n(\beta) \rvert =  O_{\mathbb P^{\tX} \otimes \mathbb P^{\tY}}(n^{-1/2})$. Thus we conclude that for any given $R > 0$, 
\begin{eqnarray*}
\sup_{\beta \in \overline{\mathcal{B}}(0, R)}  \vert \ell_{n, m}(\beta) - \ell(\beta) \vert   = o_{\mathbb P \otimes \mathbb P^{\tX} \otimes \mathbb P^{\tY}}(1).
\end{eqnarray*}

Next, we need to show that the regression model based on combining the matched and unmatched variables is identifiable. First, we show $\ell(\beta) \leq \ell(\beta_0)$ for all $\beta \in \mathbb{R}^p$. Define the random variable 
\begin{align}\label{W}
   W  =  \left(\frac{H(\beta^\top X)}{H(\beta_0^\top X)}\right)^Y \left(\frac{1-
    H(\beta^\top X)}{1-H(\beta_0^\top X)}\right)^{1-Y}.
\end{align}
For $\beta \in \mathbb R^p$, we can write
\begin{align*}
    &\ell_{M}(\beta) -  \ell_{M}(\beta_0) \\
    &= \int  \left \{ H( \beta_0^\top x) \log {\left(\frac{H(\beta^\top x)}{H(\beta_0^\top x)}\right)} + (1- H(\beta_0^\top x)) \log{\left( \frac{1 -H(\beta^\top x)}{1 -H(\beta_0^\top x)}\right) } \right \} dF^X(x) \\
    &= \mathbb{E}\left[\log{W}\right]\\
    &\leq \log{\mathbb{E}[W]}, \ \text{by Jensen's inequality} \\
    &= \log{\left(  \int H( \beta_0^\top x) \frac{H(\beta^\top x)}{H(\beta_0^\top x)} +(1- H(\beta_0^\top x))  \frac{1 -H(\beta^\top x)}{1 -H(\beta_0^\top x)}  dF^X(x)\right)} =0.
\end{align*}
For the unmatched part, we have for $\beta \in \mathbb R^p$
\begin{align*}
    &\ell_{U}(\beta) - \ell_{U}(\beta_0) \\
    &=  \left( \int  \tH(\beta_0^\top x) dF^{\tX}(x) \right) \log\left( \frac{\int \tH(\beta^\top x) dF^{\tX}(x)}{\int \tH(\beta_0^\top x) dF^{\tX}(x)}  \right) \\
    & \ \ + \left( 1- \int \tH(\beta_0^\top x) dF^{\tX}(x) \right) \log\left(\frac{ 1- \int \tH(\beta^\top x) dF^{\tX}(x)}{1- \int \tH(\beta_0^\top x) dF^{\tX}(x)}  \right) \\
    & \leq \log\Bigg(\int  \tH(\beta_0^\top x) dF^{\tX}(x) \frac{\int \tH(\beta^\top x) dF^{\tX}(x)}{\int \tH(\beta_0^\top x) dF^{\tX}(x)} \\
    & \ \ \ \ \ \ \ \ \ \  + \left( 1- \int \tH(\beta_0^\top x) dF^{\tX}(x) \right) \frac{ 1- \int \tH(\beta^\top x) dF^{\tX}(x)}{1- \int \tH(\beta_0^\top x) dF^{\tX}(x)}   \Bigg), \text{by concavity of $t  \mapsto \log(t)$} \\
    & = 0.
\end{align*}
We conclude that 
\begin{align*}
    \ell(\beta) = \frac{1}{\lambda+1}\ell_{U}(\beta) + \frac{\lambda}{\lambda+1}\ell_{M}(\beta) \leq  \frac{1}{\lambda+1}\ell_{U}(\beta_0) + \frac{\lambda}{\lambda+1}\ell_{M}(\beta_0) = \ell(\beta_0)
\end{align*}
for all $\beta \in \mathbb{R}^p$.

Next, we show that the previous inequality is strict, that is $\ell(\beta) < \ell(\beta_0)$ for all $\beta \in \mathbb R^p \setminus \{\beta_0\}$. Let $\varphi$ be strictly convex function. Then, for the integrable random variable $W$ defined in (\ref{W}) and $\varphi$ such that $W$ belongs to the domain of $\varphi$ almost surely the Jensen's inequality 
$$
\mathbb E[\varphi(W)] \ge \varphi(\mathbb E(W))
$$
is an equality if and only if $\mathbb P(W = \mathbb E(W)) = 1$. 
Since $- \log $ is strictly convex on $(0, \infty)$, it holds that 
\begin{align*}
&\int \left\{ H( \beta_0^\top x) \log {\left(\frac{H(\beta^\top x)}{H(\beta_0^\top x)}\right)} + (1- H(\beta_0^\top x)) \log{\left( \frac{1 -H(\beta^\top x)}{1 -H(\beta_0^\top x)}\right) } \right \} dF^X(x) \\
&= \mathbb{E}\left[\log{W}\right]= 0
\end{align*}
if and only if
$$
\mathbb{P}\left(W = \mathbb E \left[W\right]  \right) = 1. 
$$
Since $\displaystyle \mathbb E \left[ W\right] =\int \left \{ H( \beta_0^\top x) \frac{H(\beta^\top x)}{H(\beta_0^\top x)} +(1- H(\beta_0^\top x))  \frac{1 -H(\beta^\top x)}{1 -H(\beta_0^\top x)} \right \} dF^X(x)=1$, we conclude that 
\begin{eqnarray*}
H(\beta^\top x)  = H(\beta_0^\top x)
\end{eqnarray*}
for almost every $x  \in \operatorname{supp}(X) $. Since $H$ is strictly monotone, it implies that
$$
(\beta_0 - \beta)^\top x =0
$$
for almost all $x \in \operatorname{supp}(X)$.  This means that
$$
\mathbb P( u^\top  X  =0)  =1
$$
with $u = \beta_0  - \beta$.   By Assumption (A4), this implies that $u =0$ and hence $\beta = \beta_0$. We conclude that $\ell(\beta) < \ell(\beta_0)$ for all $\beta \ne  \beta_0$.

Now, let $r > 0$. we show next that 
$$
\sup_{\beta \in \mathcal O_r} \ell(\beta) < \ell(\beta_0).
$$

As shown above, there exists $R^\ast$ so that with probability 1 the maximizer  of the empirical likelihood can be restricted to the ball $\overline{\mathcal{B}}(0, R^\ast)$ for $n$ and $m$ large enough.  Therefore, showing the result of the theorem is equivalent to showing that
$$
\sup_{\beta \in \mathcal O_r \cap \overline{\mathcal{B}}(0, R^\ast)} \ell(\beta) < \ell(\beta_0).
$$
To avoid trivialities, we can take $r \in (0,1)$ so that $\mathcal O_r \cap \overline{\mathcal{B}}(0, R^\ast) \ne \emptyset$. In fact, if theorem is shown for $r$ small, then it will continue to hold for larger values of $r$. Hence, we assume in what follows that $r$. Using the fact that $\mathcal{O}_r \subseteq \overline{\mathcal{O}}_{r/2} =  \{ \beta \in \mathbb R^p:  \Vert \beta - \beta_0 \Vert \ge r/2 \}$, it follows that
\begin{eqnarray*}
\sup_{\beta \in \mathcal O_r \cap \overline{\mathcal{B}}(0, R^\ast)} \ell(\beta)  \le \sup_{\beta \in \overline{\mathcal O}_{r/2} \cap \overline{\mathcal{B}}(0, R^\ast)} \ell(\beta).
\end{eqnarray*}
The set $\overline{\mathcal O}_{r/2} \cap \overline{\mathcal{B}}(0, R^\ast)$ is compact as it is the intersection of a closed set and a compact set. Since the function $\beta \mapsto \ell(\beta)$ is continuous, it follows that $\beta \mapsto \ell(\beta)$ attains its supremum over this set at some $\beta^\ast \ne \beta_0$. Using the result obtained above, it holds that 
$$
\sup_{\beta \in \overline{\mathcal O}_{r/2} \cap \overline{\mathcal{B}}(0, R^\ast)} \ell(\beta) = \ell(\beta^\ast) < \ell(\beta_0)
$$
By \cite[Corollary 3.2.3]{aadbookE2}, it follows that $\widehat \beta_{n,m}  \to_{\mathbb P \otimes \mathbb P^{\tX} \otimes \mathbb P^{\tY}} \beta_0$. This finishes the proof.
\end{proof}

\section{Proofs of convergence} \label{Appendix: Convergence}

\begin{proof}[Proof of Lemma \ref{lemma: dot ell convergence in distribution}]
We have 
\begin{align*}
    \dot{\ell}_{n,m}(\beta_0) = &\frac{n}{n+m} \left(\frac{\bar{\tY}_n}{\int \tH(\beta_0^\top x) d\mathbb{P}_n^{\tX}(x)}-\frac{1-\bar{\tY}_n}{\int 1-\tH(\beta_0^\top x) d\mathbb{P}_n^{\tX}(x)}\right) \int \tH'(\beta_0^\top x) x \, d\mathbb{P}_n^{\tX}(x) \\
    &+ \frac{m}{n+m} \int \left(\frac{y}{H(\beta_0^\top x)} - \frac{1-y}{1-H(\beta_0^\top x)} \right) H'(\beta_0^\top x)x \, d\mathbb{P}_m(x,y).
\end{align*}
We start with the second term. Since 
\begin{align*}
    &\int \left(\frac{y}{H(\beta_0^\top x)} - \frac{1-y}{1-H(\beta_0^\top x)} \right) H'(\beta_0^\top x)x \, d\mathbb{P}(x,y) \\
    &=   \int \frac{H'(\beta_0^\top x)}{H(\beta_0^\top x)\left(1-H(\beta_0^\top x)\right)}  \left(y-H(\beta_0^\top x)\right)x \, d\mathbb{P}(x,y) \\
    &= \mathbb{E}_X\left[X\frac{H'(\beta_0^\top X)}{H(\beta_0^\top X)\left(1-H(\beta_0^\top X)\right)}\mathbb{E}_Y[ Y-H(\beta_0^\top X)\big|X]\right] \\
    &= \mathbb{E}_X\left[X\frac{H'(\beta_0^\top X)}{H(\beta_0^\top X)\left(1-H(\beta_0^\top X)\right)}\cdot 0\right] = 0
\end{align*}
and 
\begin{align*} 
    & \int \left(\frac{H'(\beta_0^\top x)}{H(\beta_0^\top x)\left(1-H(\beta_0^\top x)\right)} \right)^2 \left(y-H(\beta_0^\top x)\right)^2 \left(xx^\top\right) \, d\mathbb{P}(x,y)   \\
    &= \mathbb{E}_X\left[XX^\top \left(\frac{H'(\beta_0^\top X)}{H(\beta_0^\top X)\left(1-H(\beta_0^\top X)\right)}\right)^2\mathbb{E}_Y[ \left(Y-H(\beta_0^\top X)\right)^2\big|X]\right]    \\
    &= \int xx^\top \frac{\left(H'(\beta_0^\top x)\right)^2}{H(\beta_0^\top x)\left(1-H(\beta_0^\top x)\right)} \, d\mathbb{P}^X(x) = \Sigma_2,
\end{align*}
by the CLT, we have 
\begin{align*}
    \sqrt{m} \int \left(\frac{y}{H(\beta_0^\top x)} - \frac{1-y}{1-H(\beta_0^\top x)} \right) H'(\beta_0^\top x)x \, d\mathbb{P}_m(x,y) \overset{d}{\to} \mathcal{N}(0,\Sigma_2)
\end{align*}
and hence
\begin{align*}
    \sqrt{n+m}\frac{m}{n+m} \int \left(\frac{y}{H(\beta_0^\top x)} - \frac{1-y}{1-H(\beta_0^\top x)} \right) H'(\beta_0^\top x)x \, d\mathbb{P}_m(x,y) \overset{d}{\to} \mathcal{N}(0,\frac{\lambda}{1+\lambda}\Sigma_2).
\end{align*}
For the first term we have
\begin{align*}
     &\left(\frac{\bar{\tY}_n}{\int \tH(\beta_0^\top x) d\mathbb{P}_n^{\tX}(x)}-\frac{1-\bar{\tY}_n}{\int 1-\tH(\beta_0^\top x) d\mathbb{P}_n^{\tX}(x)}\right) \int \tH'(\beta_0^\top x) x \, d\mathbb{P}_n^{\tX}(x) \\
     &= \frac{\int \tH'(\beta_0^\top x) x \, d\mathbb{P}_n^{\tX}(x) }{\int \tH(\beta_0^\top x) d\mathbb{P}_n^{\tX}(x) \int 1- \tH(\beta_0^\top x) d\mathbb{P}_n^{\tX}(x)}\left(\int y \,d\mathbb{P}_n^{\tY}(y) - \int \tH(\beta_0^\top x) d\mathbb{P}_n^{\tX}(x)\right).
\end{align*}
Define 
\begin{align*}
    &W_n := \int \tH'(\beta_0^\top x) x \, d\mathbb{P}_n^{\tX}(x) \quad  &\text{and}  \quad  &W := \int \tH'(\beta_0^\top x) x \, d\mathbb{P}^{\tX}(x) \\
    &Q_n := \int \tH(\beta_0^\top x) \, d\mathbb{P}_n^{\tX}(x)  \quad  &\text{and}  \quad &Q := \int \tH(\beta_0^\top x)\, d\mathbb{P}^{\tX}(x) = \int  y \,d\mathbb{P}^{\tY}(y)
\end{align*}
We have 
\begin{align*}
    &\frac{\int \tH'(\beta_0^\top x) x \, d\mathbb{P}_n^{\tX}(x) }{\int \tH(\beta_0^\top x) d\mathbb{P}_n^{\tX}(x) \int 1-\tH(\beta_0^\top x) d\mathbb{P}_n^{\tX}(x)}\left(\int y \,d\mathbb{P}_n^{\tY}(y) - \int \tH(\beta_0^\top x) d\mathbb{P}_n^{\tX}(x)\right) \\
    &= \frac{W_n}{Q_n(1-Q_n)} \left( \int y \,d\mathbb{P}_n^{\tY}(y) - Q_n\right) \\
    &= \frac{1}{\sqrt{n}} \left( \frac{W_n}{Q_n(1-Q_n)} - \frac{W}{Q(1-Q)} \right) \left(\int y \,d\mathbb{G}_n^{\tY}(y) - \int \tH(\beta_0^\top x) \,d\mathbb{G}_n^{\tX}(x) \right) \\
    & \ \ \ \ + \frac{1}{\sqrt{n}} \frac{W}{Q(Q-1)} \int \tH(\beta_0^\top x) \,d\mathbb{G}_n^{\tX}(x) + \frac{1}{\sqrt{n}} \frac{W}{Q(1-Q)} \int y \,d\mathbb{G}_n^{\tY}(y). 
\end{align*}
By CLT we have
\begin{align*}
    \frac{W}{Q(1-Q)} \int y \,d\mathbb{G}_n^{\tY}(y) =\frac{\int \tH'(\beta_0^\top x) x \, d\mathbb{P}^{\tX}(x) }{\int \tH(\beta_0^\top x) d\mathbb{P}^{\tX}(x) \int 1-\tH(\beta_0^\top x) d\mathbb{P}^{\tX}(x)} \int y \,d\mathbb{G}_n^{\tY} \overset{d}{\to} \mathcal{N}(0,\Gamma_1)
\end{align*}
with 
\begin{align*}
    \Gamma_1 = \frac{\left(\int \tH'(\beta_0^\top x) x \, d\mathbb{P}^{\tX}(x) \right)\left(\int \tH'(\beta_0^\top x) x^\top \, d\mathbb{P}^{\tX}(x) \right)}{\int \tH(\beta_0^\top x) d\mathbb{P}^{\tX}(x) (1-\int \tH(\beta_0^\top x) d\mathbb{P}^{\tX}(x))}.
\end{align*}
In addition, we have also by CLT 
\begin{align*}
     &\frac{W}{Q(Q-1)} \int \tH(\beta_0^\top x) \,d\mathbb{G}_n^{\tX}(x)\\
     & = \frac{ \int \tH'(\beta_0^\top x) x \, d\mathbb{P}^{\tX}(x)}{\int \tH(\beta_0^\top x) d\mathbb{P}^{\tX}(x) \left( \int \tH(\beta_0^\top x) d\mathbb{P}^{\tX}(x)- 1\right)}  \int H(\beta_0^\top x) \,d\mathbb{G}_n^{\tX}(x) \overset{d}{\to} \mathcal{N}(0,\Gamma_2)
\end{align*}
with
\begin{align*} 
    \Gamma_2 &= \left(\int \tH'(\beta_0^\top x) x \, d\mathbb{P}^{\tX}(x) \right)\left(\int \tH'(\beta_0^\top x) x^\top \, d\mathbb{P}^{\tX}(x) \right)   \\
    &\phantom{=} \times\frac{ \int \tH(\beta_0^\top x)^2\, d\mathbb{P}^{\tX}(x)  - \left(\int \tH(\beta_0^\top x)\, d\mathbb{P}^{\tX}(x)\right)^2 }{\left(\int \tH(\beta_0^\top x) d\mathbb{P}^{\tX}(x)\right)^2 \left( 1- \int \tH(\beta_0^\top x) d\mathbb{P}^{\tX}(x)\right)^2} .
\end{align*}
As shown in the following Lemma \ref{lemma: res = o1} we have
\begin{align*}
    \left( \frac{W_n}{Q_n(1-Q_n)} - \frac{W}{Q(1-Q)} \right) \left(\int y \,d\mathbb{G}_n^{\tY}(y) - \int \tH(\beta_0^\top x) \,d\mathbb{G}_n^{\tX}(x) \right) = O_{\mathbb{P}^{\tX} \otimes \mathbb{P}^{\tY}}(\frac{1}{\sqrt{n}}).
\end{align*}
Since $\mathcal{N}(0,\Gamma_1)$ and $\mathcal{N}(0,\Gamma_2)$ are two independent Gaussian random variables, we have
\begin{align*}
    &\sqrt{n+m}\frac{n}{n+m} \left(\frac{\bar{\tY}_n}{\int \tH(\beta_0^\top x) d\mathbb{P}_n^{\tX}(x)}-\frac{1-\bar{\tY}_n}{\int 1-\tH(\beta_0^\top x) d\mathbb{P}_n^{\tX}(x)}\right) \int \tH'(\beta_0^\top x) x \, d\mathbb{P}_n^{\tX}(x) \\
    &\overset{d}{\to} \mathcal{N}(0,\frac{1}{1+\lambda} (\Gamma_1+\Gamma_2)).
\end{align*}
It follows that
\begin{align*}
    \sqrt{m+n}\ \dot{\ell}_{n,m}(\beta_0) \overset{d}{\to} \mathcal{N}(0,\,\frac{1}{1+\lambda} (\Gamma_1+\Gamma_2) +\frac{\lambda}{1+\lambda}\Sigma_2).
\end{align*}
The claim follows from Lemmas \ref{lemma: Sigma2 finite} and \ref{lemma: Gamma1 and Gamma2 finite}.
\end{proof}

\begin{lemma} \label{lemma: res = o1}
    Under the conditions and notation of Lemma \ref{lemma: dot ell convergence in distribution}, we have
    $$\left( \frac{W_n}{Q_n(1-Q_n)} - \frac{W}{Q(1-Q)} \right) \left(\int y \,d\mathbb{G}_n^{\tY}(y) - \int \tH(\beta_0^\top x) \,d\mathbb{G}_n^{\tX}(x) \right) = O_{\mathbb{P}^{\tX} \otimes \mathbb{P}^{\tY}}(\frac{1}{\sqrt{n}}).$$
\end{lemma}

\medskip

\begin{proof}
    We have 
    \begin{align*}
         &\frac{W_n}{Q_n(1-Q_n)} - \frac{W}{Q(1-Q)} \\
         &= \left(W_n-W\right) \left(\frac{1}{Q_n(1-Q_n)} - \frac{1}{Q(1-Q)}\right) + \left(W_n-W\right) \frac{1}{Q(1-Q)}  \\
         & \phantom{=} + W\left(\frac{1}{Q_n(1-Q_n)} - \frac{1}{Q(1-Q)}\right).
    \end{align*}
    By CLT we have
    \begin{align*}
        &\left(W_n-W\right) \frac{1}{Q(1-Q)} \left(\int y \,d\mathbb{G}_n^{\tY}(y) - \int \tH(\beta_0^\top x) \,d\mathbb{G}_n^{\tX}(x) \right) \\
        &=\frac{1}{Q(1-Q)} \cdot  \frac{1}{\sqrt{n}} \int \tH'(\beta_0^\top x) x\, d\mathbb{G}_n^{\tX}(x) \left(\int y \,d\mathbb{G}_n^{\tY}(y) - \int \tH(\beta_0^\top x) \,d\mathbb{G}_n^{\tX}(x) \right) \\
        &= O_{\mathbb{P}^{\tX} \otimes \mathbb{P}^{\tY}}(\frac{1}{\sqrt{n}}).
    \end{align*}
    Define $h(t): = \frac{1}{t(1-t)}$. By applying $\delta$-technique, it follows
    \begin{align*}
        h(Q_n)- h(Q) &= h'(t)\big|_{t = Q} \left(Q_n - Q\right) + o_{\mathbb{P}^{\tX}}(\frac{1}{\sqrt{n}}) \\
        &= \frac{2Q-1}{Q^2(1-Q)^2} \left(Q_n - Q\right) + o_{\mathbb{P}^{\tX}}(\frac{1}{\sqrt{n}}) .
    \end{align*}
    By Assumption (A5-2) and monotonicity of $s \mapsto (2s-1)/(s(1-s))^2$ on $(0,1)$ we have
    \begin{align*}
        \left\lvert \frac{2Q-1}{Q^2(1-Q)^2}\right\rvert \leq \max{\left\{ \frac{\lvert 2 \widetilde{c}_{\lVert\beta_0\rVert}-1\rvert }{\widetilde{c}_{\lVert\beta_0\rVert}^2(1-\widetilde{c}_{\lVert\beta_0\rVert})^2}, \frac{\lvert 2 \widetilde{C}_{\lVert\beta_0\rVert}-1\rvert }{\widetilde{C}_{\lVert\beta_0\rVert}^2(1-\widetilde{C}_{\lVert\beta_0\rVert})^2} \right\}} < \infty.
    \end{align*}
    Then it follows by the CLT that 
    \begin{align*}
        &\left(W_n-W\right) \left(\frac{1}{Q_n(1-Q_n)} - \frac{1}{Q(1-Q)}\right) + W\left(\frac{1}{Q_n(1-Q_n)} - \frac{1}{Q(1-Q)}\right) \\
        &=\frac{2Q-1}{Q^2(1-Q)^2}\left(W_n-W\right)  \left(Q_n - Q\right) + o_{\mathbb{P}^{\tX}}(n^{-1}) + \frac{2Q-1}{Q^2(1-Q)^2}W  \left(Q_n - Q\right) \\
        &= O_{\mathbb{P}^{\tX}}(\frac{1}{\sqrt{n}}).
    \end{align*}
    As we also have
    \begin{align*}
       \int y \,d\mathbb{G}_n^{\tY}(y) - \int \tH(\beta_0^\top x) \,d\mathbb{G}_n^{\tX}(x) = O_{\mathbb{P}^{\tX} \otimes \mathbb{P}^{\tY}}(\frac{1}{\sqrt{n}}),
    \end{align*}
    by the CLT, the claim follows.
\end{proof}

\medskip

\begin{proof}[Proof of Lemma \ref{lemma: ddot ell convergence in probability}]
\begin{align*}
    \ddot{\ell}_{n,m}(\beta_0) &= \frac{n}{n+m} \Bigg\{ \left( \frac{\bar{\tY}_n}{\int \tH(\beta_0^\top x) \, d\mathbb{P}_n^{\tX}(x)}-\frac{1-\bar{\tY}_n}{1-\int \tH(\beta_0^\top x) \, d\mathbb{P}_n^{\tX}(x)}\right) \int \tH''(\beta_0^\top x) xx^\top \, d\mathbb{P}_n^{\tX}(x)\\
    &\phantom{= \frac{n}{n+m} \Bigg\{} -\left( \frac{\bar{\tY}_n}{\left(\int \tH(\beta_0^\top x) \, d\mathbb{P}_n^{\tX}(x)\right)^2}+\frac{1-\bar{\tY}_n}{\left(1-\int \tH(\beta_0^\top x) \, d\mathbb{P}_n^{\tX}(x)\right)^2}\right) \\
    &\phantom{= \frac{n}{n+m} \Bigg\{ -}\times \left(\int \tH'(\beta_0^\top x) x \, d\mathbb{P}_n^{\tX}(x)\right)\left(\int \tH'(\beta_0^\top x) x^\top \, d\mathbb{P}_n^{\tX}(x)\right) \Bigg\}    \\
    & \quad + \frac{m}{n+m} \int \Bigg(\left(\frac{y}{H(\beta_0^\top x)} - \frac{1-y}{1-H(\beta_0^\top x)}\right) H''(\beta_0^\top x)\\
    &\phantom{\quad + \frac{m}{n+m} \int \Bigg(} - \left(\frac{y}{H(\beta_0^\top x)^2} + \frac{1-y}{(1-H(\beta_0^\top x))^2}\right) \left(H'(\beta_0^\top x)\right)^2 \Bigg) xx^\top \,d\mathbb{P}_m(x,y) \\
    & := \frac{n}{n+m}  \mathrm{I}_n + \frac{m}{n+m} \mathrm{II}_m.
\end{align*}
We start with $\mathrm{II}_m$. Since SLLN implies that
\begin{align*}
    &\int xx^\top \left(\frac{y}{H(\beta_0^\top x)} - \frac{1-y}{1-H(\beta_0^\top x)}\right) H''(\beta_0^\top x) \,d\mathbb{P}_m(x,y) \\
    &\overset{\mathbb{P}}{\to} \mathbb{E}_X \left[ XX^\top  H''(\beta_0^\top X)\, \mathbb{E}_Y \left[\frac{Y}{H(\beta_0^\top X)} - \frac{1-Y}{1-H(\beta_0^\top X)} \big| X\right] \right] = 0,
\end{align*}
it follow that
\begin{align*}
    \mathrm{II}_m &= - \int \left(\frac{y}{H(\beta_0^\top x)^2} + \frac{1-y}{(1-H(\beta_0^\top x))^2}\right) \left(H'(\beta_0^\top x)\right)^2  xx^\top \,d\mathbb{P}_m(x,y) \\
    &\overset{\mathbb{P}}{\to} -\mathbb{E}_X \left[ XX^\top \left(H'(\beta_0^\top X)\right)^2 \mathbb{E}_Y \left[\frac{Y}{H(\beta_0^\top X)^2} + \frac{1-Y}{(1-H(\beta_0^\top X))^2}  \big| X\right] \right] \\
    &= -\mathbb{E}_X \left[ XX^\top  \frac{\left(H'(\beta_0^\top X)\right)^2 }{H(\beta_0^\top X) \left(1 - H(\beta_0^\top X)\right)} \right] = -\Sigma_2.
\end{align*}
Thus we have
\begin{align*}
    \frac{m}{n+m} \mathrm{II}_m \to_{\mathbb{P}} - \frac{\lambda}{1+\lambda} \Sigma_2.
\end{align*}
We continue with $ \mathrm{I}_n$.
By SLLN and continuous mapping theorem we have
\begin{align*}
    &\left( \frac{\bar{\tY}_n}{\int \tH(\beta_0^\top x) \, d\mathbb{P}_n^{\tX}(x)}-\frac{1-\bar{\tY}_n}{1-\int \tH(\beta_0^\top x) \, d\mathbb{P}_n^{\tX}(x)}\right) \int \tH''(\beta_0^\top x) xx^\top \, d\mathbb{P}_n^{\tX}(x) \\
    &= \left(\left(\frac{1}{U_n}- \frac{1}{U}\right) \left(V_n - V\right) + \frac{1}{U} (V_n-V) + \left(\frac{1}{U_n}- \frac{1}{U}\right)V + \frac{V}{U}\right) \int \tH''(\beta_0^\top x) xx^\top \, d\mathbb{P}_n^{\tX}(x) \\
    & \phantom{=} -\left(\left(\frac{1}{U_n-1}- \frac{1}{U-1}\right) \left(V_n - V\right) + \frac{1}{U-1} (V_n-V) + \left(\frac{1}{U_n-1}- \frac{1}{U-1}\right)(V-1) + \frac{V-1}{U-1}\right) \\
    &\phantom{\phantom{=} -}  \times  \int \tH''(\beta_0^\top x) xx^\top \, d\mathbb{P}_n^{\tX}(x) \\
& \overset{\mathbb{P}^{\tX} \otimes \mathbb{P}^{\tY} }{\to} \frac{V}{U} \int \tH''(\beta_0^\top x) xx^\top \, d\mathbb{P}^{\tX}(x) -  \frac{V-1}{U-1} \int \tH''(\beta_0^\top x) xx^\top \, d\mathbb{P}^{\tX}(x) 
 = 0
\end{align*}
with $U_n = \int \tH(\beta_0^\top x) \, d\mathbb{P}_n^{\tX}(x)$, $U = \int \tH(\beta_0^\top x) \, d\mathbb{P}^{\tX}(x)$, $V_n = \int y \, d\mathbb{P}_n^{\tY}(y)$, and $V = \int y \, d\mathbb{P}^{\tY}(y) = \int \tH(\beta_0^\top x) \, d\mathbb{P}^{\tX}(x) = U$. In addition, we can write the second part of $\mathrm{I}_n$ as following:
\begin{align*}
     &-\left( \frac{\bar{\tY}_n}{\left(\int \tH(\beta_0^\top x) \, d\mathbb{P}_n^{\tX}(x)\right)^2}+\frac{1-\bar{\tY}_n}{\left(1-\int \tH(\beta_0^\top x) \, d\mathbb{P}_n^{\tX}(x)\right)^2}\right) \\
     & \quad \times\left(\int \tH'(\beta_0^\top x) x \, d\mathbb{P}_n^{\tX}(x)\right)\left(\int \tH'(\beta_0^\top x) x^\top \, d\mathbb{P}_n^{\tX}(x)\right) \\
     & = -\left( \frac{\bar{\tY}_n}{\left(\int \tH(\beta_0^\top x) \, d\mathbb{P}^{\tX}(x)\right)^2}+\frac{1-\bar{\tY}_n}{\left(1-\int \tH(\beta_0^\top x) \, d\mathbb{P}^{\tX}(x)\right)^2}\right) \\
     & \phantom{= -} \times\left(\int \tH'(\beta_0^\top x) x \, d\mathbb{P}^{\tX}(x)\right)\left(\int \tH'(\beta_0^\top x) x^\top \, d\mathbb{P}^{\tX}(x)\right) + r_n \\
\end{align*}
with \begin{align*}
     r_n &= \frac{1}{\mathbb{E} \left[\tY\right] \left(1-\mathbb{E} \left[\tY\right]\right)} \left(BB^\top -B_nB_n^\top\right) \\
     & \phantom{= } + \frac{1}{(\mathbb{E} \left[\tY\right])^2} \left(BB^\top -B_nB_n^\top\right) \cdot \frac{1}{\sqrt{n}} \int y \,d\mathbb{G}_n^{\tY}(y) \\
     & \phantom{= } - \frac{1}{\left(1-\mathbb{E} \left[\tY\right]\right)^2} \left(BB^\top -B_nB_n^\top\right) \cdot \frac{1}{\sqrt{n}} \int y \,d\mathbb{G}_n^{\tY}(y) \\
     & \phantom{= }+  \left(\frac{1}{A_n^2} - \frac{1}{A^2}\right)\left(BB^\top -B_nB_n^\top\right) \cdot \frac{1}{\sqrt{n}} \int y \,d\mathbb{G}_n^{\tY}(y)\\
    & \phantom{= }+ \mathbb{E} \left[\tY\right] \left(BB^\top -B_nB_n^\top\right) \left(\frac{1}{A_n^2} - \frac{1}{A^2}\right)\\
    & \phantom{= } - BB^\top   \left(\frac{1}{A_n^2} - \frac{1}{A^2}\right)  \frac{1}{\sqrt{n}} \int y \,d\mathbb{G}_n^{\tY}(y) \\
     & \phantom{= }+     BB^\top   \mathbb{E} \left[\tY\right] \left(\frac{1}{A_n^2} - \frac{1}{A^2}\right)\\
    & \phantom{= } -\left(\frac{1}{(1-A_n)^2} - \frac{1}{(1-A)^2}\right)\left(BB^\top -B_nB_n^\top\right) \cdot \frac{1}{\sqrt{n}} \int y \,d\mathbb{G}_n^{\tY}(y) \\
    & \phantom{= } + \left(1- \mathbb{E}\left[\tY\right]\right)\left(\frac{1}{(1-A_n)^2} - \frac{1}{(1-A)^2}\right)\left(BB^\top -B_nB_n^\top\right) \\
    & \phantom{= } + \left(\frac{1}{(1-A_n)^2} - \frac{1}{(1-A)^2}\right)BB^\top \cdot \frac{1}{\sqrt{n}} \int y \,d\mathbb{G}_n^{\tY}(y) \\
     & \phantom{= } -  \left(1- \mathbb{E}\left[\tY\right]\right)  \left(\frac{1}{(1-A_n)^2} - \frac{1}{(1-A)^2}\right)BB^\top
\end{align*}
where 
\begin{align*}
    &A = \int \tH(\beta_0^\top x) d\mathbb{P}^{\tX}(x) \quad &\text{and} \quad & A_n = \int \tH(\beta_0^\top x) d\mathbb{P}_n^{\tX}(x) \\
    &B = \int \tH'(\beta_0^\top x) x \,  d\mathbb{P}^{\tX}(x) \quad &\text{and} \quad & B_n = \int \tH'(\beta_0^\top x) x \,  d\mathbb{P}_n^{\tX}(x)
\end{align*}
SLLN and continuous mapping theorem imply that $r_n \overset{\mathbb{P}^{\tX} \otimes \mathbb{P}^{\tY} }{\longrightarrow} 0$. Thus, it follows that
\begin{align*}
    &-\left( \frac{\bar{\tY}_n}{\left(\int \tH(\beta_0^\top x) \, d\mathbb{P}^{\tX}(x)\right)^2}+\frac{1-\bar{\tY}_n}{\left(1-\int \tH(\beta_0^\top x) \, d\mathbb{P}^{\tX}(x)\right)^2}\right) \\
     & \phantom{-} \times\left(\int \tH'(\beta_0^\top x) x \, d\mathbb{P}^{\tX}(x)\right)\left(\int \tH'(\beta_0^\top x) x^\top \, d\mathbb{P}^{\tX}(x)\right) \\
     & \longrightarrow_{\mathbb{P}^{\tY}}-\left( \frac{1}{\int \tH(\beta_0^\top x) \, d\mathbb{P}^{\tX}(x)}+\frac{1}{1-\int \tH(\beta_0^\top x) \, d\mathbb{P}^{\tX}(x)}\right) \\
     & \phantom{\longrightarrow_{\mathbb{P}^{\tY}}-} \times\left(\int \tH'(\beta_0^\top x) x \, d\mathbb{P}^{\tX}(x)\right)\left(\int \tH'(\beta_0^\top x) x^\top \, d\mathbb{P}^{\tX}(x)\right) = -\Gamma_1.
\end{align*}
Thus we have
\begin{align*}
      \ddot{\ell}_{n,m}(\beta_0) \to_{\mathbb{P}^{\tX} \otimes \mathbb{P}^{\tY}  \otimes \mathbb{P} } -\frac{1}{1+\lambda} \Gamma_1 - \frac{\lambda}{1+\lambda} \Sigma_2.
\end{align*}
The claim follows from Lemmas \ref{lemma: Sigma2 finite} and \ref{lemma: Gamma1 and Gamma2 finite}.
\end{proof}

\medskip

\section{Proof of example with Logistic link and a uniform covariate} \label{appendix: Uniform covariate example}

\begin{proof}[Proof of Lemma \ref{lemma: G_infinity X uniform H logis}]
By evenness, we can focus on the limit as $t \to \infty$. For $t$ large, we know that 
\begin{eqnarray*}
I_1(t)  & =  &   \frac{1}{t} \log\left( \frac{1+ e^t}{2} \right) \\
& = &  \frac{1}{t} ( \log(1+e^t)  - \log(2))  \\
& = & \frac{1}{t} ( t - \log(2) + \log(e^{-t}  + 1) ).
\end{eqnarray*}
Hence, 
\begin{eqnarray*}
I_1(t)=1-\frac{\log 2}{t}+o(t^{-1}), 
\end{eqnarray*}
and
\begin{eqnarray*}
I_1(t)(1-I_1(t)) =\frac{\log 2}{t}+o(t^{-1}).
\end{eqnarray*}
Also, for $t$ large
\begin{eqnarray*}
I_2(t)  & = &  \frac{1}{t} H(t)  - \frac{1}{t^2} \log\left( \frac{1+ e^t}{2} \right)  \\
& = &  \frac{1}{t}  \frac{1}{1+e^{-t}} - \frac{1}{t^2} \left( -\log 2 + t + \log(e^{-t} + 1)  \right)   \\
& = &  -\frac{1}{t} \frac{e^{-t}}{1+e^{-t}}   + \frac{\log 2}{t^2} -  \frac{1}{t^2} \log(e^{-t} + 1) \\
& = &  \frac{\log 2}{t^2}  + o(t^{-2}).
\end{eqnarray*}
For $I_3$, we have that
\begin{eqnarray*}
I_3(t) & = &  \frac{1}{t} \left(\log\left( \frac{1+ e^t}{2} \right) - H(t) + \frac{1}{2} \right) \\
& = &  \frac{1}{t} \left( - \log(2) + t + \log(e^{-t} +1) - \frac{1}{2}  + o(1) \right) \\
& = &  1  - \frac{\log 2 + 1/2}{t}  + o(t^{-1}).
\end{eqnarray*}
Lastly, using a change of variables, we can write
\begin{eqnarray*}
I_4(t) =  \frac{1}{t^3} \int_0^t  x^2 H'(x) dx.
\end{eqnarray*}
Using integrations by parts, we have
\begin{eqnarray*}
\int_0^\infty  x^2 H'(x) dx  & = &  - [x^2 (1-H(x))]_0^\infty + 2 \int_0^\infty x (1-H(x)) dx  \\
& = & 2 \int_0^\infty x (1-H(x)) dx \\
& =&  2 \int_0^\infty x  \frac{1}{1 + e^x} dx \\
& = & 2 \int_0^\infty x  \frac{e^{-x}}{1 + e^{-x}} dx.
\end{eqnarray*}
Using the expansion $e^{-x}/(1+e^{-x})  = e^{-x} \sum_{k=0}^\infty  (-1)^k e^{-kx}  = \sum_{k=1}^\infty  (-1)^{k-1} e^{-kx}$ we can write 
\begin{eqnarray*}
\int_0^\infty x  \frac{e^{-x}}{1 + e^{-x}} dx  &  =  &  \sum_{k=1}^\infty (-1)^{k-1} \int_0^\infty x e^{-kx} dx \\
& = &  \sum_{k=1}^\infty (-1)^{k-1} \frac{1}{k^2}  = \frac{\pi^2}{12}.
\end{eqnarray*}
Thus, 
\begin{eqnarray*}
I_4(t) =  \frac{1}{t^3} \frac{\pi^2}{6} + o(t^{-3}).
\end{eqnarray*}
Finally,
\begin{eqnarray*}
I_3(t) - I_1(t)^2 & = &  1  - \frac{\log 2 + 1/2}{t} - \left(1-\frac{\log 2}{t}+o(t^{-1}) \right)^2 + o(t^{-1}) \\
& = &    \frac{\log 2 - 1/2}{t}  + o(t^{-1}).
\end{eqnarray*}
Therefore, 
\begin{eqnarray*}
\tau(t)=\frac{I_2(t)^2}{I_1(t)\{1-I_1(t)\}I_4(t)} \to  \frac{6\log 2}{\pi^2},
\end{eqnarray*}
and
\begin{eqnarray*}
\tau (t)\left(2-\kappa(t)\right)
        &=&   \frac{I_2(t)^2}{I_4(t)}
        \left[
        \frac{1}{I_1(t)\{1-I_1(t)\}}
        +
        \frac{I_3(t)-I_1(t)^2}{I_1(t)^2\{1-I_1(t)\}^2}
        \right] \\
        &\to & \frac{12\log 2-3}{\pi^2}.
\end{eqnarray*}
We conclude that  $G(t, \infty) \to G_\infty$ as $\vert t \vert \to \infty$ where $G_\infty$ is from (\ref{Ginfty}).
\end{proof}

\section{Proofs of auxiliary lemmas used to establish existence} \label{Appendix: Existence auxiliary lemma}

\begin{proof}[Proof of Lemma \ref{Lemma: Case I holds a.s.}]
    Put 
     \begin{align*}
     p_0 := \int H(\beta_0^\top x) dF^X(x).
     \end{align*}
    By Assumption (A2) we have $0 < p_0 <1$. Note that for all $i \in \{1, \dots, m \}$
     \begin{align*}
         \mathbb{E}\left(Y_i\right) = \mathbb{P}\left(Y_i = 1\right) = \int \mathbb{P}\left(Y_i = 1 \big| X_i=x \right) dF^X(x) = \int H(\beta_0^\top x) dF^X(x) = p_0.
     \end{align*}
     By SLLN there exists almost surely $m' \geq 1$ such that for all $m \geq m'$
     \begin{align*}
         &0< \frac{p_0}{2}  \leq \frac{1}{m} \sum_{i=1}^m Y_i.
     \end{align*}
     Since $Y_i \in \{0,1\}$ for all $i=1,\cdots,m$, it means that there exists at least one $i \in \{1, \dots, m\}$ such that $Y_i=1$. 
     Similarly, there exists almost surely $m'' \geq 1$ such that for all $m \geq m''$ 
      \begin{align*}
         \frac{1}{m}\sum_{j=1}^m (1- Y_j) \geq  \frac{1-p_0}{2} > 0,
     \end{align*}
     implying that there exists at least one $j \in \{1, \dots, m \}$ such that $Y_j=0$. 
     Taking $\tilde{m} := \max{\{m',m''\}}$ the lemma follows.
\end{proof}

\medskip

\begin{proof}[Proof of Lemma \ref{lemma: Case I holds a.s. 2}]
 By Assumption (A4) there exists $ r > 0$ such that $\mathcal{B}(a,r) \subseteq \operatorname{supp}(X)$. Note that there exists some small $\epsilon > 0$ such that
 \begin{align*}
     \mathcal{B}\left(a,\epsilon \right) \subseteq a + \left [-\frac{r}{2\sqrt{p}}, \frac{r}{2\sqrt{p}} \right ]^{p} \subseteq \mathcal{B}(a,r)
 \end{align*}
 where $\displaystyle a + \left [-\frac{r}{2\sqrt{p}}, \frac{r}{2\sqrt{p}}\right]^{p}$ is a hypercube with center $a \in \mathbb{R}^p$. Let  $c_1,\cdots,c_{2^p}$ denote the $2^p$ corners of this hypercube. For each of these corners,  there exists a neighborhood $\displaystyle \mathcal{B}\left(c_k,\frac{r}{4\sqrt{p}} \right)$ for $k =1, \cdots, 2^p$, such that for every point $ \displaystyle d_k \in \mathcal{B}\left(c_k,\frac{r}{4\sqrt{p}}\right)$ we have
 \begin{align*}
   \mathcal{B}\left(a,\epsilon\right) \subseteq  \operatorname{int} (\operatorname{conv}(d_1,\dots,d_{2^p})). 
 \end{align*}
For example, we can choose $\epsilon = \frac{r}{8\sqrt{p}}$. Now, Assumption (A1) implies that  
\begin{align*}
    \mathbb{P}\left(X \in \mathcal{B}\left(c_k,\frac{r}{4\sqrt{p}}\right),\, Y = 1\right) 
    =  \int_{\mathcal{B}\left(c_k,\frac{r}{4\sqrt{p}}\right)} H(\beta_0^\top x) f^X(x)\, dx  > 0 
\end{align*}
for all $k=1, \cdots, 2^p$. For any fixed $k \in \{1,\cdots, 2^p\}$ the SLLN implies that , with probability 1, there exists $\tilde{m}$ such that for all $m \geq \tilde{m}$
\begin{align*}
    \frac{1}{m}\sum_{i=1}^m \mathbbm{1}_{\{X_i \in \mathcal{B}(c_k,\frac{r}{4\sqrt{p}}),\, Y_i = 1\}} > \frac{1}{2} \mathbb{P}\left(X \in \mathcal{B}(c_k,\frac{r}{4\sqrt{p}}),\, Y = 1\right) > 0.
\end{align*}
This implies that there exists at least one $j \in \{1,\cdots, r_m\}$ such that $\displaystyle X_j \in \mathcal{B}\left(c_k,\frac{r}{4\sqrt{p}}\right)$. By the same argument, across all $2^p$ neighborhoods of the corners, with probability 1, there exists $\tilde{m} \geq 1$ such that for all $m \geq \tilde{m}$, each neighborhood contains at least one successful response with covariate variable $X$ lying in it. We denote these covariates by $X^{c}_1, \cdots, X^{c}_{2^p}$. Then,
\begin{align*}
     \mathcal{B}\left(a, \epsilon \right) \subseteq  \operatorname{int} (\operatorname{conv}(X^{c}_1,\dots,X^{c}_{2^p})) \subseteq   \operatorname{conv}(X_1,\dots,X_{r_m}). 
\end{align*}
Similarly, we have
\begin{align*}
     \mathcal{B}\left(-a,\epsilon \right) \subseteq \operatorname{conv}(-X_{r_m+1},\dots,-X_{m}). 
\end{align*}
Then, it follows that
\begin{align*}
\operatorname{conv} \left( \mathcal{B}(a,\epsilon) \cup  \mathcal{B}(-a,\epsilon)\right) \subseteq \operatorname{conv}(X_1,\cdots,X_{r_m}, -X_{r_m+1},\cdots, -X_m)
\end{align*}
and hence
\begin{align*}
    0\in\operatorname{int}(\operatorname{conv}(X_1,\cdots,X_{r_m}, -X_{r_m+1},\cdots, -X_m)).
\end{align*}
\end{proof}

\medskip

\begin{proof} [Proof of Lemma \ref{lemma: ULLN existence asymptotic}]
By Lemma \ref{lemma: bi log concavity left derivative} we have
    \begin{align*}
    \log{H}(-R \lVert x \rVert) \leq  \frac{H_-'(0)}{H(0)} \left(-R  \lVert x \rVert\right) < 0,
    \end{align*}
    which implies that $\displaystyle \left( \log{H}(-R \lVert x \rVert)\right)^2 \geq \left(R\frac{H_-'(0)}{H(0)}\right)^2 \lVert x \rVert^2 $. By Assumption (A3) this implies that 
    \begin{align}\label{NormBounded}
        \int \lVert x \rVert^2 dF^X(x) \leq \frac{1}{\left(R\frac{H_-'(0)}{H(0)}\right)^2} \int \left( \log{H}(-R \lVert x \rVert)\right)^2 dF^X(x) < \infty.
    \end{align}
    Define the function
    \begin{align*}
        f_v(x,y) &= v^\top x \mathbbm{1}_{\{ v^\top x < 0\}} \mathbbm{1}_{\{y=1\}} - v^\top x \mathbbm{1}_{\{v^\top x > 0\}} \mathbbm{1}_{\{y=0\}} \\
        &= \min{\{v^\top x, 0\}} \mathbbm{1}_{\{y=1\}} - \max{\{ v^\top x,0 \}}\mathbbm{1}_{\{y=0\}},
    \end{align*}
 and let 
$$
\mathcal{F} :=  \{f_v: v\in \mathcal{S}^{p-1}\}.
$$
 Since the functions $t \mapsto \min(t, 0)$ and $t \to \max(t, 0)$ are 1-Lipchitz, it follows that for all $v,w \in \mathcal{S}^{p-1}$ 
    \begin{align*}
       \lvert f_v(x,y) - f_w(x,y) \rvert \leq  \lvert (v-w)^\top x\rvert \leq  \lVert v-w\rVert \lVert x\rVert.
    \end{align*}
   In addition, for all $v \in \mathcal S^{p-1}$,  $\vert f_v(x,y) \vert \le \Vert v \Vert \Vert x \Vert$, implying that $F(x, y) := \Vert x \Vert $ is an envelope of $\mathcal F$. Then, by \cite[Theorem 2.7.17]{aadbookE2} we have
    \begin{align*}
        N_{[\ ]}\left(2\epsilon \lVert F \rVert_{\mathbb{P},1} ,\mathcal{F}, L_1(\mathbb{P})  \right ) \leq N(\epsilon, \mathcal{S}^{p-1},\lVert\cdot\rVert).
    \end{align*}
    As stated in page 147 of \cite{aadbookE2}, $N(\epsilon, \mathcal{S}^{p-1},\lVert\cdot\rVert) \leq D(\epsilon, \mathcal{S}^{p-1},\lVert\cdot\rVert)$, where the latter denotes the packing number which is the maximum of  $\epsilon$-separated points in $\mathcal{S}^{p-1}$. By \cite[Problems and Complements (7) in page 143]{aadbookE2} we know that
    \begin{align*}
        D(\epsilon, \mathcal{S}^{p-1},\lVert\cdot\rVert) \leq D(\epsilon, \bar{\mathcal{B}}(0,1),\lVert\cdot\rVert) \leq \left(\frac{3}{\epsilon}\right)^p \quad \text{for all $0 < \epsilon \leq 1$}.
    \end{align*}
    Since $\lVert F \rVert_{\mathbb{P},1} = \int \lVert x\rVert \,d\mathbb{P}(x,y) = \mathbb{E}\left[\lVert X\rVert\right] \leq \sqrt{\mathbb{E}\left[\lVert X\rVert^2\right]} < \infty$ by (\ref{NormBounded}), it follows that
    \begin{align*}
         N_{[\ ]}\left(\epsilon  , \mathcal{F}, L_1(\mathbb{P})  \right ) \leq \left(\frac{6\mathbb{E}\left[\lVert X\rVert\right]} {\epsilon}\right)^p <\infty \quad \text{for all $0 < \epsilon \leq 2\mathbb{E}[\lVert X\rVert] $}.
    \end{align*}
    And for all $\epsilon >  2\mathbb{E}[\lVert X\rVert]$ we have 
    \begin{align*}
        N_{[\ ]}\left(\epsilon  , \mathcal{F}, L_1(\mathbb{P})  \right ) \leq  N_{[\ ]}\left( 2\mathbb{E}[\lVert X\rVert] , \mathcal{F}, L_1(\mathbb{P})  \right ) <\infty. 
    \end{align*}
    Thus, $ N_{[\ ]}\left(\epsilon  , \mathcal{F}, L_1(\mathbb{P})  \right ) < \infty$ for all $\epsilon>0$. By \cite[Theorem 2.4.1]{aadbookE2} we have $\mathcal{F} = \{f_v: v\in \mathcal{S}^{p-1}\}$ is Glivenko-Cantelli. Therefore, we have
    \begin{align*}
        &\left\lvert \sup_{v\in \mathcal{S}^{p-1}} \int f_v(x,y)\, d\mathbb{P}_m(x,y) - \sup_{v\in \mathcal{S}^{p-1}} \int f_v(x,y)\, d\mathbb{P}(x,y)  \right\rvert \\
        &\leq \sup_{v\in \mathcal{S}^{p-1}} \left\lvert \int f_v(x,y)\, d\mathbb{P}_m(x,y) -  \int f_v(x,y)\, d\mathbb{P}(x,y)  \right\rvert \overset{\text{a.s.}}{\longrightarrow} 0.
    \end{align*}
\end{proof}

\medskip

\begin{proof}[Proof of Lemma \ref{lemma: asymptotic non-zero denominator}]
    Define
    \begin{align*}
        f(v)&:= \mathbb{E}\left[v^\top X  \mathbbm{1}_{\{v^\top X<0\}}  \mathbbm{1}_{\{Y=1\}} \right] + \mathbb{E}\left[-v^\top X  \mathbbm{1}_{\{v^\top X>0\}}  \mathbbm{1}_{\{Y=0\}} \right] \\
        &=\mathbb{E}\left[\min{\{v^\top X,0\}}  \mathbbm{1}_{\{Y=1\}} \right] + \mathbb{E}\left[- \max{\{v^\top X,0\}}\mathbbm{1}_{\{Y=0\}} \right].
    \end{align*}
     We have shown in Lemma \ref{lemma: ULLN existence asymptotic} that $\mathbb{E}[\lVert X\rVert^2] < \infty$, and hence $\mathbb{E}[\lVert X\rVert] \le \sqrt{\mathbb{E}[\lVert X\rVert^2]} < \infty$. Thus $f$ is 1-Lipschitz. In fact, 
    \begin{align*}
        \lvert f(v) - f(w)\rvert &\leq \left\lvert \mathbb{E}\left[\left(\min{\{v^\top X,0\}}- \min{\{w^\top X,0\}}\right)  \mathbbm{1}_{\{Y=1\}} \right]\right\rvert \\
        & \ \ \ \  + \left\lvert \mathbb{E}\left[-\left(\max{\{v^\top X,0\}}- \max{\{w^\top X,0\}}\right)  \mathbbm{1}_{\{Y=0\}} \right] \right\rvert \\
        &\leq \mathbb{E}\left[\left\lvert \left(v -w\right)^\top X \right\rvert   \mathbbm{1}_{\{Y=1\}} \right] + \mathbb{E}\left[\left\lvert \left(v -w\right)^\top X \right\rvert   \mathbbm{1}_{\{Y=0\}} \right] \\
        &\leq \mathbb{E}\left[\left\lVert v-w\right\rVert \left\lVert X\right\rVert\right] \\
        &= \mathbb{E}\left[\lVert X\rVert\right]\lVert v-w\rVert. 
    \end{align*}
    Since $\mathcal{S}^{p-1}$ is compact, it follows that there exists $v^\ast \in \mathcal{S}^{p-1}$ such that  $\displaystyle  \sup_{v\in \mathcal{S}^{p-1}} f(v) = f(v^\ast)$. It is clear that $f(v) \le 0$, but we will show that $f(v) < 0$ for all $v\in \mathcal{S}^{p-1}$. If there exists $w \in \mathcal{S}^{p-1}$ such that $f(w) = 0$, then we have almost surely $\min{\{w^\top X,0\}}  \mathbbm{1}_{\{Y=1\}}  = 0$ and $\max{\{w^\top X,0\}}\mathbbm{1}_{\{Y=0\}} = 0$. Since $\min{\{w^\top X,0\}}  \mathbbm{1}_{\{Y=1\}}  = 0 $ is equivalent to $ \{w^\top X \ge 0 \} \cup  \{Y = 0 \}$, we must have
    \begin{align*}
        0&= \mathbb{P}\left(w^\top X< 0 , Y =1\right) \\
        &= \mathbb{E}\left[ \mathbbm{1}_{\{w^\top X < 0\}} \mathbbm{1}_{\{Y=1\}} \right]\\
         &= \mathbb{E} \left[\mathbb{E}\left[ \mathbbm{1}_{\{w^\top X < 0\}} \mathbbm{1}_{\{Y=1\}} \mid X\right]\right]\\
        &= \mathbb{E}\left[ \mathbbm{1}_{\{w^\top X < 0\}} \mathbb{P}\left(Y=1 \mid X\right) \right] \\
        &=  \mathbb{E}\left[ \mathbbm{1}_{\{w^\top X < 0\}} H(\beta_0^\top X) \right]
    \end{align*}
    By Assumption (A1) we have $H > 0$, thus we have almost surely $\mathbbm{1}_{\{w^\top X < 0\}} = 0$. It follows that we must almost surely have $w^\top X\geq 0$. Similarly, we can show that we must have  $w^\top X\leq 0$ almost surely. Thus, we conclude that $ \mathbb P(w^\top X= 0) = 1$. Since $w \ne 0$ (recall that $\Vert w \Vert =1$), this conclusion is in contradiction with Assumption (A4). Therefore, we must have $ \sup_{v\in \mathcal{S}^{p-1}} f(v) = f(v^\ast) < 0$. 
\end{proof}

\section{Proofs of auxiliary lemmas used to establish convergence} \label{Appendix: convergence auxialary lemma}
\begin{proof} [Proof of Lemma \ref{lemma: Sigma2 finite}]
    For any concave function $g$ we have for all $t \geq a$
    \begin{align*}
        g'(t) \leq \frac{g(t) - g(a)}{t-a}.
    \end{align*}
    By (A1'), we put $g = \log{H}$. Then, for all $a < \beta_0^\top x$, we have
    \begin{align*}
         \frac{H'(\beta_0^\top x)}{H(\beta_0^\top x)} = \frac{d}{dt}\log{H(t)}|_{t=\beta_0^\top x} &\leq \frac{\log{H}\left(\beta_0^\top x\right) - \log{H}\left(a\right)}{\beta_0^\top x - a} \\
         & \leq -\frac{ \log{H}\left(a\right)}{\beta_0^\top x - a} \quad \text{since $\log{H}\left(\beta_0^\top x\right) < 0$ by (A1).}
    \end{align*}
    By Cauchy-Schwarz inequality it holds that
    \begin{align*}
        \beta_0^\top x \geq -\lVert\beta_0\rVert \lVert x\rVert > -\left(\lVert\beta_0\rVert +1\right)\lVert x\rVert.
    \end{align*}
    We put $a =  -\left(\lVert\beta_0\rVert +1\right)\lVert x\rVert$, then $t - a \geq \lVert x\rVert$. It follows that
    \begin{align*}
        \frac{H'(\beta_0^\top x)}{H(\beta_0^\top x)} \leq -\frac{ \log{H}\left(-\left(\lVert\beta_0\rVert +1\right)\lVert x\rVert\right)}{\lVert x\rVert}.
    \end{align*}
    Similarly, we have
    \begin{align*}
        \frac{H'(\beta_0^\top x)}{1-H(\beta_0^\top x)} \leq -\frac{ \log{\left(1-H\left(\left(\lVert\beta_0\rVert +1\right)\lVert x\rVert\right)\right)}}{\lVert x\rVert}.
    \end{align*}
    Thus it follows by (A3) that
    \begin{align*}
         &\int \lVert x\rVert^2 \frac{\left(H'(\beta_0^\top x)\right)^2}{H(\beta_0^\top x)\left(1-H(\beta_0^\top x)\right)} \, d\mathbb{P}^X(x)  \\
         &\leq \int   \log{H}\left(-\left(\lVert\beta_0\rVert +1\right)\lVert x\rVert\right) \cdot \log{\left(1-H\left(\left(\lVert\beta_0\rVert +1\right)\lVert x\rVert\right)\right)}  \, d\mathbb{P}^X(x) \\
         &\leq \frac{1}{2} \int   \log{H}\left(-\left(\lVert\beta_0\rVert +1\right)\lVert x\rVert\right)^2  \, d\mathbb{P}^X(x)+ \frac{1}{2}\int  \log{\left(1-H\left(\left(\lVert\beta_0\rVert +1\right)\lVert x\rVert\right)\right)}^2  \, d\mathbb{P}^X(x)\\
         &< \infty.
    \end{align*}
\end{proof}

\medskip

\begin{proof}[Proof of Lemma \ref{lemma: Gamma1 and Gamma2 finite}]
    By (A1) we have $0 < H(t) < 1$ for all $t\in\mathbb{R}$, then it follows that
    \begin{align*}
        &0 < \int \tH(\beta_0^\top x) d\mathbb{P}^{\tX}(x) (1-\int \tH(\beta_0^\top x) d\mathbb{P}^{\tX}(x)) \leq \frac{1}{4} 
    \end{align*}
    and 
    \begin{align*}
        0 \leq \operatorname{Var}\left(\tH\left(\beta_0^\top \tX\right)\right) = \int \tH(\beta_0^\top x)^2\, d\mathbb{P}^{\tX}(x)  - \left(\int \tH(\beta_0^\top x)\, d\mathbb{P}^{\tX}(x)\right)^2 \leq 1.
    \end{align*}
    For any vector $v \in \mathbb{R}^p$ we have
    \begin{align*}
        &v^\top \left(\int \tH'(\beta_0^\top x) x \, d\mathbb{P}^{\tX}(x) \right)\left(\int \tH'(\beta_0^\top x) x^\top \, d\mathbb{P}^{\tX}(x) \right) v \\
        &= \left(\int \left( v^\top x \right) \tH'(\beta_0^\top x) \, d\mathbb{P}^{\tX}(x) \right)^2 \\
        &= \left(\int \left( v^\top x \right) \frac{\tH'(\beta_0^\top x)}{\sqrt{\tH(\beta_0^\top x)\left(1-\tH(\beta_0^\top x)\right)}}\sqrt{\tH(\beta_0^\top x)\left(1-\tH(\beta_0^\top x)\right)} \, d\mathbb{P}^{\tX}(x) \right)^2 \\
        &\leq \left(\int \left( v^\top x \right)^2 \frac{\left(\tH'(\beta_0^\top x)\right)^2}{\tH(\beta_0^\top x)\left(1-\tH(\beta_0^\top x)\right)}\, d\mathbb{P}^{\tX}(x) \right) \\
        &\phantom{\leq} \times \left(\int \tH(\beta_0^\top x)\left(1-\tH(\beta_0^\top x)\right)\, d\mathbb{P}^{\tX}(x) \right) \quad \text{by Cauchy-Schwarz inequality} \\
        &= v^\top \left(\int xx^\top \frac{\left(\tH'(\beta_0^\top x)\right)^2}{\tH(\beta_0^\top x)\left(1-\tH(\beta_0^\top x)\right)}\, d\mathbb{P}^{\tX}(x) \right) v\\
        &\phantom{\leq} \times \left(\int \tH(\beta_0^\top x)\left(1-\tH(\beta_0^\top x)\right)\, d\mathbb{P}^{\tX}(x) \right)
    \end{align*}
    Using similar arguments in Lemma \ref{lemma: Sigma2 finite}, we have here 
    \begin{align*}
        \left(\int xx^\top \frac{\left(\tH'(\beta_0^\top x)\right)^2}{\tH(\beta_0^\top x)\left(1-\tH(\beta_0^\top x)\right)}\, d\mathbb{P}^{\tX}(x) \right) < \infty.
    \end{align*}
    Since $ 0< \int \tH(\beta_0^\top x)\left(1-\tH(\beta_0^\top x)\right)\, d\mathbb{P}^{\tX}(x)  \leq 1/4$, we have $\Gamma_1 < \infty$ and $\Gamma_2 < \infty$.
\end{proof}

\medskip

\begin{proof} [Proof of Lemma \ref{lemma: Sigma2 positive definite}]
    For all $v \in \mathbb{R}^p \setminus \{0\}$, we have
    \begin{align*}
        v^\top \Sigma_2 v &=  \int \left(v^\top x\right)^2 \frac{\left(H'(\beta_0^\top x)\right)^2}{H(\beta_0^\top x)\left(1-H(\beta_0^\top x)\right)} \, d\mathbb{P}^X(x) \\
        &\geq 4 \int \left(v^\top x\right)^2 \left(H'(\beta_0^\top x)\right)^2 \, d\mathbb{P}^X(x).
    \end{align*}
    By (A1) and (A1'), $H'(\beta_0^\top X)>0$ almost surely. Hence the integrand is
    nonnegative and is strictly positive whenever $v^\top X\neq 0$. By (A4), $X$ admits a density with respect to the Lebesgue measure  and thus $\mathbb{P}\left(v^\top X = 0 \right) =0$. Therefore we have $ v^\top \Sigma_2 v  > 0$.
\end{proof}

\end{document}